\documentclass[12pt, a4paper]{article}
\usepackage{amssymb,latexsym, amsmath,amscd,amsfonts,amsthm}%,refcheck,mathbbol, mathrsfs}
\usepackage[hypertex]{hyperref} %xdvi in jussieu wants this line, no idea
\newtheorem{prop}{Proposition}[section]
\newtheorem{defi}[prop]{Definition}

\newtheorem{lem}[prop]{Lemma}
\newtheorem{thm}[prop]{Theorem}

\newtheorem{remar}[prop]{Remark}
\newtheorem{cor}[prop]{Corollary}

\newtheorem{ass}[prop]{Assumption}

\DeclareMathAlphabet{\mathpzc}{OT1}{pzc}{m}{it}
\DeclareMathOperator{\Aut}{Aut}
\DeclareMathOperator{\End}{End}
\DeclareMathOperator{\Hom}{Hom}
\DeclareMathOperator{\Ind}{Ind}
\DeclareMathOperator{\cInd}{c-Ind}

\DeclareMathOperator{\GL}{GL}

\DeclareMathOperator{\Ker}{Ker}

\DeclareMathOperator{\tr}{tr}

\DeclareMathOperator{\supp}{Supp}
\DeclareMathOperator{\id}{id}

\DeclareMathOperator{\vol}{vol}
\newcommand{\cIndu}[3]{\cInd_{#1}^{#2}{#3}}

\newcommand{\Indu}[3]{\Ind_{#1}^{#2}{#3}}

\newcommand{\curlya}{{\mathfrak A}}
\newcommand{\UU}{\mathbf U}

\newcommand{\urlyb}{\mathfrak B}

\newcommand{\curlyb}[1]{\mathfrak B_{#1}}

\newcommand{\ub}[2]{\mathbf U^{#1}(\curlyb{#2})}
\newcommand{\ubz}[1]{\mathbf U(\curlyb{#1})}
\newcommand{\LL}{\mathcal L}

\newcommand{\curlyp}{\mathfrak P}

\newcommand{\oF}{\mathfrak o_F}
\newcommand{\pF}{\mathfrak p_F}
\newcommand{\oE}{\mathfrak o_E}
\newcommand{\oK}{\mathfrak o_K}

\newcommand{\pL}{\mathfrak p_L}
\newcommand{\pK}{\mathfrak p_K}
\newcommand{\pE}{\mathfrak p_E}
\newcommand{\kF}{\mathfrak k_F}
\newcommand{\kE}{\mathfrak k_E}
\newcommand{\kA}{\mathfrak K(\curlya)}

\newcommand{\eBoe}[1]{e(\curlyb{#1}|\oE)}

\newcommand{\FF}{\mathcal F}
\newcommand{\Fq}{\mathbf{F}_q}
\newcommand{\ZZ}{\mathbb{Z}}
\newcommand{\trAF}{{\tr}_{A/F}}
\newcommand{\pif}{\varpi_F}
\newcommand{\mir}{\mathcal M} % mirabolic
\newcommand{\rim}{\mathcal M_{\curlya}} % rimabolic
\newcommand{\bes}{\mathcal J}
\newcommand{\MM}{\mathbb M}

\newcommand{\BB}{\mathcal B}
\newcommand{\ee}{\mathrm e}
\newcommand{\CC}{\mathbb C}
\newcommand{\rV}{\mathrm V}

\newcommand{\rw}{\mathrm w}
\newcommand{\rF}{\mathrm F}
\newcommand{\rG}{\mathrm G}
\newcommand{\rK}{\mathrm K}
\newcommand{\WW}{\mathcal W}
\newcommand{\VV}{\mathcal V}
\newcommand{\Eins}{\mathbf 1} 
\newcommand{\kB}{\mathfrak K(\urlyb)}
\newcommand{\KK}{\mathfrak K}
\newcommand{\GG}{\mathcal G}

\def\bdots{\mathinner{\mkern1mu\raise1pt\hbox{.}\mkern2mu\raise4pt\hbox{.}
           \mkern2mu\raise7pt\vbox{\kern7pt\hbox{.}}\mkern1mu}}
\def\Re{{\mathrm R\mathrm e}}
\def\volU{\vol_{U}}
\def\volF{\vol_{F}}

\def\labelenumi{{\rm(\roman{enumi})}}
\def\theenumi{\roman{enumi}}

\title{On the realisation of maximal simple types and %$\varepsilon$-factors
 epsilon factors of pairs}
\author{Vytautas Paskunas and Shaun Stevens}
\date{\today}
\begin{document} 
\maketitle

\begin{abstract}
Let~$G$ be the group of rational points of a general linear group over 
a non-archimedean local field~$F$. We show that certain
representations of open, compact-mod-centre subgroups of~$G$, (the maximal
simple types of Bushnell and Kutzko) can be realized 
as concrete spaces. In the level zero case our result is essentially
due to Gel$'$fand. This allows us, for a supercuspidal
representation~$\pi$ of~$G$, to compute a distinguished
matrix coefficient of~$\pi$. By integrating, we obtain an explicit
Whittaker function for~$\pi$. We use this to compute
the~$\varepsilon$-factor of pairs, for
supercuspidal representations~$\pi_1$,~$\pi_2$ of~$G$, when~$\pi_1$
and~ the contragredient of $\pi_2$ differ only at the `tame level' (more
precisely,~$\pi_1$ and~$\check{\pi}_2$ contain the same simple character).  
We do this by
computing both sides of the functional equation defining the epsilon
factor, using the definition of Jacquet, Piatetskii-Shapiro, Shalika. We
also investigate the behaviour of the~$\varepsilon$-factor under twisting 
of~$\pi_1$ by tamely ramified quasi-characters. Our
results generalise the special case~$\pi_1=\check{\pi}_2$ %and~$\pi_1$ 
totally wildly ramified, due to Bushnell and Henniart.
\end{abstract}

\tableofcontents

\newpage

%%%%%%%%%%%%%%%%%%%%%%%%%%%%%%%%%%%%%%%%%%%%%%%%%%%%%%%%
%%%%%%%%%%%%%%%%%%%%%% Introduction %%%%%%%%%%%%%%%%%%%%
%%%%%%%%%%%%%%%%%%%%%%%%%%%%%%%%%%%%%%%%%%%%%%%%%%%%%%%%

\section*{Introduction}
\addcontentsline{toc}{section}{Introduction}%\numberline {}Introduction}
\def\theprop{\Alph{prop}}

Let~$F$ be a non-archimedean local field and fix an additive 
character~$\psi_F$ of~$F$, with conductor~$\pF$ (the maximal ideal 
of the ring of integers~$\oF$ of~$F$). Let~$V$ be
an~$N$-dimensional~$F$-vector space and let~$G=\Aut_F(V)$. 
This paper concerns the supercuspidal representations of~$G$, so we
adopt the notation of~\cite{bk}, where these are classified in terms
of \emph{maximal simple types}~(\cite{bk}\S6). 

A maximal simple type is a pair~$(J,\lambda)$, consisting of a (rather
special) compact open subgroup~$J$ of~$G$ and an irreducible
representation~$\lambda$ of~$J$. It is constructed from a 
\emph{simple stratum}~$[\curlya,n,0,\beta]$, where~$\curlya$ is a 
principal 
hereditary~$\oF$-order in~$A$ and $\beta\in A$ satisfy various 
properties (see~\cite{bk}~(1.5.5)
for the full definition). The algebra~$E=F[\beta]$ is a field,
and~$E^{\times}$ normalises~$J$. If we set~$\mathbf J=E^{\times} J$
and let~$\Lambda$ be a representation of~$\mathbf J$ 
such that~$\Lambda|_{J}\cong \lambda$, then~$\cIndu{\mathbf
  J}{G}{\Lambda}$ is an irreducible supercuspidal representation
of~$G$. Conversely, any irreducible supercuspidal
representation of~$G$ arises this way. 

It follows from the fact
that~$\cIndu{\mathbf J}{G} {\Lambda}$ has a Whittaker model, that
there exist a maximal unipotent subgroup~$U$ of~$G$ and a
non-de\-ge\-ne\-rate character $\psi_{\alpha}$ of~$U$ such that~$\Hom_{U\cap
  \mathbf J} (\psi_{\alpha}, \Lambda)\neq 0$. Moreover, the uniqueness of the
Whittaker model implies that the pair~$(U,\psi_{\alpha})$ is
determined up to conjugation by~$\mathbf J$ (see~\cite{bh}). 

Associated to the maximal simple type, we have two other groups
$H^1\subset J^1\subset J$ (see~\cite{bk}\S3.1). They are normal
subgroups of~$J$ and~$\lambda|_{H^1}$ is a multiple of a 
\emph{simple character}~$\theta$. 
Hence,~$\psi_{\alpha}|_{U\cap H^1}= \theta|_{U\cap H^1}$ and 
we may define a character $\Psi$ of~$(J\cap U) H^1$, by  
$$ 
\Psi(uh)= \psi_{\alpha}(u)\theta(h), \quad \forall u\in J\cap U, \quad
\forall h\in H^1.
$$          
Let~$0\subset V_1\subset\cdots  \subset V_N= V$ be the maximal flag
corresponding to $U$, and set~$\mir=\{g\in G: (g-1) V\subseteq
V_{N-1}\}$, so that~$\mir$ is a mirabolic subgroup of~$G$. 
Our first main result, in~\S\ref{super}, is

\begin{thm}\label{Ione}
Let~$U$ be a maximal unipotent subgroup
of~$G$ and let~$\psi_{\alpha}$ be a non-degenerate character of~$U$ such
that~$\Hom_{U\cap \mathbf J}(\psi_{\alpha}, \Lambda)\neq 0$.
Then~$\Lambda|_{(\mir\cap \mathbf J) J^1}$ is irreducible and 
$$
\Lambda|_{(\mir\cap \mathbf J) J^1}\cong \Indu{(U\cap \mathbf
  J)H^1}{(\mir\cap \mathbf J) J^1}{\Psi}.
$$
Moreover, the same result holds if we replace:~$\mathbf J$ by~$\kA$, 
the~$G$-normaliser of~$\curlya$; $\Lambda$ 
by~$\rho=\Indu{\mathbf J}{\kA}{\Lambda}$; 
and~$(\mir\cap \mathbf J)J^1$ by~$(\mir\cap
\kA)\UU^1(\curlya)$, where~$\UU^1(\curlya)$ is the group of principal 
units of~$\curlya$. 
\end{thm}

We show in~\S\ref{charind} that this property in fact
characterises the representations of the form~$\Indu{\mathbf
  J}{\kA}{\Lambda}$. More precisely, 
writing~$\rim= (\mir\cap \kA)\UU^1(\curlya)$ we have

\begin{prop}\label{Itwo}
Let~$\tau$ be a representation of~$\kA$ such that 
$$ 
\tau|_{\rim}\cong 
\Indu{(J\cap U) H^1}{\rim} \Psi.
$$
Then 
$$ 
\tau\cong \Indu{\mathbf J}{\kA}{\Lambda'},
$$
for some representation~$\Lambda'$ of~$\mathbf J$, such
that~$(J,\Lambda'|_J)$ is a maximal simple type containing~$\theta$.
\end{prop}

This is an analogue of Gel$'$fand's characterisation of the  cuspidal
representations of~$\GL_N(\Fq)$.

\smallskip

Before continuing with the applications of Theorem~\ref{Ione}, we 
will say a few words about its proof.
The strategy is as follows: We first 
construct a special pair~$(U, \psi_{\alpha})$, such that~$\Hom_{U\cap
  \mathbf J}(\psi_{\alpha}, \Lambda)\neq 0$, by carefully picking a
basis of~$V$, and then letting~$U$ be the group of upper-triangular
matrices with respect to this basis, and~$\psi_{\alpha}$ be the
'standard' nondegenerate character. The choice of basis is made so
that we can control the restriction of~$\psi_{\alpha}$ to~$U\cap
B^{\times}$, where~$B$ is the centraliser of~$\beta$. We then prove
Theorem~\ref{Ione} for this particular
pair~$(U,\psi_{\alpha})$. The general case
follows from the fact that any other such pair~$(U',
\psi_{\alpha'})$ is conjugate to our particular choice by some~$g\in
\mathbf J$. 

We remark also that Theorem~\ref{Ione} should follow easily
from~\cite{bh}~Theorem~2.9, but there are problems with the proof of
that result: the unipotent group~$U$ used in the proof is not (in
general) the group required by the statement of the Theorem; moreover,
there is a gap in the proof of~\cite{bh}~Lemma~2.10 which,
so far as we know, nobody has been able to fix. In the course of the proof of 
Theorem \ref{Ione}, we end up proving an analogue of \cite{bh} Theorem 2.9, see 
Theorem \ref{wavingit}. We get around the problem of \cite{bh} Lemma 2.10 by 
using the case when~$E$ is maximal in~$A$ as a `black
box'. If~\cite{bh}~Lemma~2.10 were true then the basis of~$V$ referred
to above could be written explicitly in terms of~$\beta$ and some of
our proofs would simplify.

\medskip

The interest in Theorem~\ref{Ione} is that, 
since~$F^{\times} \Ker \Lambda$ is of finite
index in~$\mathbf J$, it allows us to apply a
very general Theorem of Alperin and James~\cite{alperin} for finite groups. 
Following~\cite{alperin}, we define the Bessel function~$\bes_\Lambda$
of~$\Lambda$ by 
$$ 
\bes_\Lambda(g)= Q^{-1} 
\sum_{u\in (U\cap\mathbf J)H^1/\UU^{n+1}(\curlya)} 
     \Psi(u)\tr_{\Lambda}(gu^{-1}).
$$ 
where~$Q=((U\cap\mathbf J)H^1:\UU^{n+1}(\curlya))$. Now, \cite{alperin} implies:
\begin{thm}\label{Ibes} Let $\mathcal S$ be the space of functions $f:(\mir\cap\mathbf J) J^1\rightarrow \CC$, such that
$$
f(ug)= \Psi(u) f(g), \quad \forall u\in (U\cap\mathbf J)H^1, \quad \forall g\in(\mir\cap \mathbf J)J^1,
$$ 
and, for all~$g\in\mathbf J$, let~$L(g)\in \End_{\CC} (\mathcal S)$ be the operator: 
$$ 
[L(g)f](m)= 
\sum_{m_1\in{(\mir\cap \mathbf J)J^1}/(U\cap\mathbf J)H^1} 
\bes_\Lambda(m g m_1) f(m_1^{-1}),
$$ 
Then~$L$ defines a representation of~$\mathbf J$ on $\mathcal S$,
which is isomorphic to~$\Lambda$.
\end{thm}
The analogous result for $\GL_N(\Fq)$ (or, equivalently, for level
zero supercuspidal representations of $G$) was first observed by
Gel$'$fand~\cite{gelfand}.

\medskip

Now we proceed as in~\cite{bh}~\S3 and
construct a Whittaker function for~$\pi=\cIndu{\mathbf
  J}{G}{\Lambda}$. However, the observation that~$\Lambda$ can be
realised via Bessel functions results in
explicit formulae and extra information about this Whittaker
function. Our main result concerning Whittaker functions (see~\S\ref{xwf})
is as follows:
 
\begin{thm}\label{Ithree} 
Define~$\WW\in
\Indu{U}{G}{\psi_{\alpha}}$ by~$\supp \WW\subseteq U\mathbf J$ and  
$$ 
\WW(ug)=\psi_{\alpha}(u)\bes_{\Lambda}(g), \quad \forall u\in U, \quad
\forall g\in \mathbf J,
$$
then~$\WW$ is a Whittaker function for~$\pi=\cIndu{\mathbf
  J}{G}{\Lambda}$. Moreover,~$(\supp \WW)\cap \mir= U( H^1\cap \mir)$
and  
$$ 
\WW(u h)= \psi_{\alpha}(u)\theta(h), \quad \forall u\in U, \quad
\forall h\in H^1\cap \mir.
$$
\end{thm} 

It is the second part of Theorem~\ref{Ithree} that really requires
Theorem~\ref{Ione} and 
the explicit realization in terms of Bessel functions. The first part of 
Theorem~\ref{Ithree} has also
been obtained by Roberto Johnson~\cite{johnson}, in the special case
when~$\pi$ is a  Carayol representation.

\medskip

Finally, in~\S\ref{epsi} we use our Whittaker functions to
compute~$\varepsilon$-factors of pairs in the following situation. We 
fix a simple stratum as above and 
consider two supercuspidal representations~$\pi_1=\cIndu{\mathbf
  J}{G}{\Lambda_1}$ and~$\pi_2=\cIndu{\mathbf J}{G}{\Lambda_2}$, such
that~$\Lambda_1|_{H^1}$ and~$\Lambda_2|_{H^1}$ are multiples of 
\emph{the same simple character~$\theta$}. Then, for~$i=1,2$, we can 
write~$\Lambda_i|_J\cong \kappa\otimes \sigma_i$, where~$\kappa$ is
a~$\beta$-{\emph extension} (see~\cite{bk}~(5.2.1)) and~$\sigma_i$ 
is the lift of a cuspidal representation of~$J/J^1\cong \GL_r(\kE)$, 
where~$\kE$ is the residue field of~$E$ and $r=\dim_E(V)$. 

One may show that, for~$i=1,2$, we have~$\Lambda_i\cong 
\tilde{\kappa}\otimes \Sigma_i$, where~$\tilde{\kappa}$ and~$\Sigma_i$ 
are representations of~$\mathbf J$ which restrict to~$\kappa$ 
and~$\sigma_i$ respectively. Moreover, we may think of~$\Sigma_i$ as 
a representation of $\kA\cap B^\times\cong E^{\times}\GL_r(\oE)$ 
(where, we recall,~$B$ is the centraliser of~$E$). We 
set~$\tau_i=\cIndu{\kB}{B^{\times}}{\Sigma_i}$; 
then~$\tau_i$ is a supercuspidal level zero representation
of~$B^{\times}\cong \GL_r(E)$. 

In Theorem~\ref{apl} we
relate~$\varepsilon(\pi\times\check{\pi}_2,s, \psi_F)$
and~$\varepsilon(\tau_1\times \check{\tau}_2, s,\psi_E)$, 
where~$\psi_E$ is an additive character of~$E$ with conductor~$\pE$ 
which extends~$\psi_F$. We obtain:

\begin{thm}\label{Ifour}
$$ 
\varepsilon(\pi_1\times \check{\pi}_2, s , \psi_F)= 
\zeta\omega_{\tau_1}(\nu^{-r})
\omega_{\tau_2}(\nu^{r})q^{(s-1/2)r v_E(\nu)N/e} 
\varepsilon(\tau_1\times \check{\tau}_2, s,\psi_E),
$$
where:~$\zeta= \omega_{\tau_2}(-1)^{r-1} \omega_{\pi_2}(-1)^{N-1}$;%
~$q=q_F$ is the cardinality of~$\kF$;~$v_E$ is the additive valuation 
on~$E$ with image~$\mathbb Z$;
and~$\nu=\nu(\theta, \psi_F,\psi_E)\in E^{\times} /(1+\pE)$ is 
an invariant which we define in~\S\ref{numerical}.
\end{thm}

We prove this by computing both sides of the functional equation 
for the epsilon factor, using the definition of Jacquet, 
Piatetskii-Shapiro, Shalika~\cite{jpss}, with the Whittaker functions
of Theorem~\ref{Ithree}. We are able to do the calculation because 
the fact that the operator~$L$ in Theorem \ref{Ibes}
 defines a group action imposes various identities 
on the Bessel function~$\bes_{\Lambda}$. 

Theorem~\ref{Ifour} implies:

\begin{cor}
 Let~$\chi:F^{\times}\rightarrow \CC^{\times}$ be a tamely 
ramified quasi-character and put~$\chi_E=\chi\circ \mathrm{N}_{E/F}$; then  
$$
\frac{\varepsilon(\pi_1\chi\times \check{\pi}_2,s,\psi_F)}
{\varepsilon(\pi_1\times \check{\pi}_2,s,\psi_F)}  =
\chi(\mathrm{N}_{E/F}(\nu^{-r^2})) 
\frac{\varepsilon(\tau_1\chi_E\times \check{\tau}_2,s,\psi_E)}
{\varepsilon(\tau_1\times \check{\tau}_2,s, \psi_E)},
$$
where~$\nu=\nu(\theta_F,\psi_F,\psi_E)$.
\end{cor}

We recover a
result of Bushnell and Henniart~\cite{can}, concerning the effect
on~$\varepsilon(\pi\times\check{\pi},s,\psi_F)$ of twisting~$\pi$ by 
tamely ramified quasi-characters, when~$\pi$ is totally wildly 
ramified, as a special case of the above Corollary.  Moreover, we show 
in~\S\ref{behaviour} that the
invariant~$\nu$ behaves well under the tame lifting operation for
simple characters of Bushnell and Henniart~\cite{ihes}, which 
implies \cite{can} Theorem 7.1.
\medskip

We end the introduction with a brief summary of the contents of each 
section. We begin in~\S\ref{notation} with notation and some elementary 
results about nondegenerate characters and induced supercuspidal 
representations. In~\S\ref{noteonbeta} we begin the groundwork for the 
proof of Theorem~\ref{Ione}, proving a similar result for 
a~$\beta$-extension~$\kappa$. In~\S\ref{Udef} we define the particular 
unipotent subgroup, and the basis, used in the proof of 
Theorem~\ref{Ione}; this proof appears in~\S\ref{super}, along with the proof 
of Proposition~\ref{Itwo}. In~\S\ref{realization}, we apply the Theorem 
of Alperin and James to define Bessel functions, and construct our 
explicit Whittaker function from Theorem~\ref{Ithree}. 
In~\S\ref{numerical}, we define the numerical invariant~$\nu$ which 
appears in Theorem~\ref{Ifour}. Finally, the proof of 
Theorem~\ref{Ifour}, and its application to twisting by tamely ramified 
quasi-characters, appears in~\S\ref{epsi}.

\medskip

\textit{Acknowledgements.} This work has grown out of discussions
which started at the LMS Symposium ``$L$-functions and Galois
representations'' in Durham, July~2004. The authors would like to
thank the organisers of the symposium, and also Richard Hill for some
productive conversations there. Part of this work was done 
while the first-named author was a CNRS post-doc at Jussieu; the
first-named author would like to thank Marie-France Vigneras, the
Institut de Math\'ematiques de Jussieu and the CNRS for providing
excellent working conditions in a stimulating environment. The
second-named author would like to thank the SFB~701 ``Spektrale
Strukturen und Topologische Methoden in der Mathematik'', Bielefeld,
for their invitation and hospitality at a crucial stage in the writing
of the paper.

%%%%%%%%%%%%%%%%%%%%%%%%%%%%%%%%%%%%%%%%%%%%%%%%%%%%%%%%%%%%%%%%%%%%%%
%%%%%%%%%%%%%%%%%%%%% Notation and Preliminaries %%%%%%%%%%%%%%%%%%%%%
%%%%%%%%%%%%%%%%%%%%%%%%%%%%%%%%%%%%%%%%%%%%%%%%%%%%%%%%%%%%%%%%%%%%%%

\section{Notation and Preliminaries}\label{notation}
\def\theprop{\arabic{section}.\arabic{prop}}

Let~$F$ be a locally compact non-archimedean local field, with ring of
integers~$\oF$, maximal ideal~$\pF$, and residue field~$\kF=\oF/\pF$
with~$q_F=p^f$ elements,~$p$ prime. We fix~$\varpi_F$
a uniformizing element of~$F$ and  let~$v_F$ denote the
additive valuation of~$F$, normalised so that $v_F(\varpi_F)=1$. We use similar
notation for any field extension of~$F$.

Let~$V$ be an~$N$-dimensional~$F$-vector space,~$A=\End_F(V)$
and~$G=\Aut_F(V)$ so, after choosing a basis for~$V$, we have 
$$
A\cong \mathbb M_N(F),\qquad G\cong \GL_N(F).
$$

%For~$r$ a real number, we write~$\lceil r\rceil$ for the smallest
%integer greater than or equal to~$r$, and~$\lfloor r\rfloor$ for the
%largest integer less than or equal to~$r$.

%%%%%%%%%%%%%%%%% Unipotent subgroups and characters %%%%%%%%%%%%%%%%%

\subsection{Unipotent subgroups and characters}\label{UandC}

The results of this section are stated without proof, since these
proofs are straightforward. One way to prove them would be to choose
a suitable basis for~$V$ with respect to which the unipotent subgroups
considered consist of matrices which are upper triangular.

We fix, once and for all, an additive character~$\psi_F:F\to\mathbb C$
which is trivial on~$\pF$, non-trivial on~$\oF$. For any~$a\in A$, we
define a function~$\psi_a:A\to\mathbb C$
by
$$
\psi_a(x)\ =\ (\psi_F\circ\trAF)(a(x-1)),\qquad\hbox{for }x\in A,
$$
where~$\trAF$ denotes the matrix trace. We use the same notation for
the restriction of~$\psi_a$ to various subsets of~$A$.

Let~$\FF$ be an~$F$-flag in~$V$,
$$
\FF:\qquad 0=V_0 \subset V_1\subset V_2\subset \cdots \subset V_s=V,
$$
Let~$P=P_\FF$ be the~$G$-stabiliser of~$\FF$, a parabolic subgroup of~$G$, and
let~$U=U_\FF$ be its unipotent radical. We also put
\begin{eqnarray*}
X_{\FF} &=&\{x\in A:xV_i\subseteq V_{i+1},\ 0\le i\le s-1\}, \\
X_{\FF}^+ &=&\{x\in X_{\FF}:xV_i\not\subset V_i,\ 0\le i\le s-1\}, \\
X_{\FF}^- &=&\{x\in X_{\FF}:xV_i\subseteq V_{i},\ 0\le i\le s-1\}.
\end{eqnarray*}

\begin{lem}\label{charX} Let~$a\in A$. The function~$\psi_a$ defines a
linear character of~$U$ if and only if~$a\in X_{\FF}$. 
Moreover,~$\psi_a$ is trivial on~$U$ if and only if~$a\in X_{\FF}^-$. 
\end{lem}

%\begin{proof} Does this need a proof? We could always just say that
%one way to prove it would be to choose a basis for~$V$ which conforms
%to the flag~$\FF$.
%\end{proof}

Now suppose that~$\FF$ is a maximal~$F$-flag so that~$s=N$ and~$\dim_F
V_i=i$, for~$0\le i\le N$. A smooth linear character~$\chi$ of~$U$ is
said to be {\it nondegenerate\/} if its~$G$-normaliser is~$F^\times
U$. We can describe this more concretely by choosing a
basis~$v_1,...,v_N$ for~$V$ such that~$V_i=\bigoplus_{j=1}^i Fv_j$,
for~$1\le i\le N$. Then~$U$ is identified with the upper triangular
unipotent matrices in~$\GL_N(F)$ and the smooth characters~$\chi$
of~$U$ are given by
$$
\chi(u)=\psi_F\left(\sum_{i=1}^{N-1} \mu_i u_{i,i+1}\right),
\qquad\hbox{for }u=\left(u_{ij}\right)\in U,
$$
where~$\mu_i\in F$,~$1\le i\le N-1$, are fixed scalars. It is easy
to see that~$\chi$ is nondegenerate if and only if all~$\mu_i\ne
0$; or, equivalently, if and only if, for all~$1\le j\le N-1$,
there exists~$u^{(j)}\in U$ with~$u^{(j)}_{i,i+1}=0$, for~$i\ne j$,
such that~$\chi(u^{(j)})\ne 1$. %We also have:

\begin{lem} Let~$\FF$ be a maximal~$F$-flag and~$a\in
X_{\FF}$. Then~$\psi_a$ is nondegenerate if and only if~$a\in
X_{\FF}^+$.
\label{nondeg}
\end{lem}

%\begin{proof} Does this need a proof?
%\end{proof}

%%%%%%%%%%%%%%%%%%%%%%%% Induced supercuspidals  %%%%%%%%%%%%%%%%%%%%%

\subsection{Induced supercuspidals}

Since the supercuspidal representations of~$G$ are all obtained by
irreducible induction from compact-mod-centre subgroups, the following
Proposition (which is mostly taken from~\cite{bh}~\S1) will be useful.

\begin{prop}\label{conju} Let~$\KK$ be an open, compact-mod-centre
subgroup of~$G$, and suppose that~$\rho$ is a representation of~$\KK$,
such that~$\pi=\cIndu{\KK}{G}{\rho}$ is an irreducible supercuspidal
representation of~$G$. Let~$U$ be a maximal unipotent subgroup of~$G$,
and let~$\chi$ be a smooth character of~$U$. %Then the following hold: 
\begin{enumerate}
\item If~$\Hom_{U\cap \KK}( \rho, \chi)\neq 0$
  then~$\chi$ is non-degenerate. 
\item If~$\chi$ is non-degenerate then there
exists~$g\in G$ such that~$\Hom_{U\cap \KK^g}( \chi, \rho^g)\neq 0$. 
\item If~$\Hom_{U\cap \KK}( \chi, \rho)\neq 0$
and~$\Hom_{U\cap \KK^g}( \chi, \rho^g)\neq 0$, for some~$g\in G$, then
there exists~$u\in U$ such that~$\KK^u= \KK^g$ 
and~$\rho^u\cong\rho^g$.
\end{enumerate} 
\end{prop}

\begin{proof} (i) Let~$\phi \in \Hom_{U\cap \KK}( \rho, \chi)\neq 0$ be such
that~$\phi\neq 0$, and fix a Haar measure~$du$ on~$U$. Then this gives a
non-zero~$\Phi \in \Hom_U(\pi, \chi)$, by  
$$ 
\Phi(f)= \int_U \chi(u^{-1}) \phi(f(u)) du, \quad \forall f \in
\cIndu{\KK}{G}{\rho}.
$$ 
If~$\chi$ is degenerate then there exists a unipotent radical~$U'$ of
some proper parabolic subgroup of~$G$, such that the restriction
of~$\chi$ to~$U'$ is trivial, but this implies~$\Hom_{U'}(\pi,
\Eins)\neq 0$. However this may not happen as~$\pi$ is
supercuspidal. 

Parts (ii) and (iii) follow from~\cite{bh}~Proposition~1.6 and~(1.8).
\end{proof}

%%%%%%%%%%%%%%%%%%%%%%%%%%%%%%%%%%%%%%%%%%%%%%%%%%%%%%%%%%%%%%%%%%%%%%
%%%%%%%%%%%%%%%%%%%%% A note on~$\beta$-extensions %%%%%%%%%%%%%%%%%%%
%%%%%%%%%%%%%%%%%%%%%%%%%%%%%%%%%%%%%%%%%%%%%%%%%%%%%%%%%%%%%%%%%%%%%%

\section{A note on~$\boldsymbol\beta$-extensions}\label{noteonbeta}

The main result of this section is Theorem~\ref{Psikappa}, which
asserts that the restriction of a~$\beta$-extension~$\kappa$ to a
certain subgroup of~$J$ is isomorphic to a representation induced from
a linear character. This result will be  used in~\S\ref{super}. In
section~\ref{Iwahori} we recall some results on Iwahori
decompositions. We will use the definitions and notations of~\cite{bk}
with little introduction.

Let~$[\curlya,n,0,\beta]$ be a principal simple stratum in~$A$
(see~\cite{bk}~(1.5.5)). In particular,~$\curlya$ is a hereditary,
principal~$\oF$-order in~$A$, with Jacobson
radical~$\curlyp$, and~$\beta\in\curlyp^{-n}\setminus\curlyp^{1-n}$
is such that~$E=F[\beta]$ is a field with~$E^\times$ 
normalising~$\curlya$. We denote by~$B$ the~$A$-centraliser of~$E$ and
put~$\urlyb=\curlya\cap B$. 

Let~$\LL=\{L_k:k\in\mathbb Z\}$ be the~$\oF$-lattice chain in~$V$
associated to~$\curlya$, see~\cite{bk}~(1.1.2). Since~$E^\times$
normalises~$\curlya$ we may consider~$\LL$ also as an~$\oE$-lattice
chain. Let~$e=e(\urlyb|\oE)$ be the~$\oE$-period of~$\LL$.  We fix
an~$E$-basis~$\{ w_1, \ldots, w_r\}$ of~$V$, where~$r[E:F]=N$, such
that 
$$
L_0= \oE w_1 + \oE w_2+\cdots +\oE w_r
$$
$$
L_i=\oE w_1+ \cdots +\oE w_{\frac{r}{e}(e-i)} + \pE
w_{\frac{r}{e}(e-i)+1}+ \cdots+\pE w_{r}
$$
for~$0< i <e$. The choice of this basis identifies~$B$ with~$\mathbb
M_r(E)$ and~$\urlyb$ with a subring of~$\mathbb M_r(\oE)$, which
consists of block upper-triangular matrices modulo~$\pE$, such that
each block on the diagonal is of the size~$\frac{r}{e}\times
\frac{r}{e}$.    

Let~$\FF_0$ be the~$F$-flag in~$V$, given by
$$ 
\FF_0: 0 \subset E w_1\subset \cdots \subset\oplus_{i=1}^j E
w_i\subset \cdots\subset\oplus_{i=1}^{r} E w_i=V.
$$ 
Let~$P_0$ be the~$G$-stabiliser of~$\FF_0$. Moreover,
put~$G_i=\Aut_F(Ew_i)$, for~$1\le i\le r$, and
$$
M_0\ =\ \prod_{i=1}^r G_i,
$$
a Levi component of the parabolic subgroup~$P_0$ of~$G$. Let~$U_0$ be
the unipotent radical of~$P_0$, so that~$P_0=M_0U_0$. We also denote
by~$U_0^-$ the unipotent radical of the parabolic subgroup opposite
to~$P_0$ relative to~$M_0$. 

We note that~$\FF_0$ is also a maximal~$E$-flag in~$V$. This yields:

\begin{lem}\label{refine} Let~$\FF$ be an~$F$-flag refining~$\FF_0$
and let~$U$ be the unipotent radical of the~$G$-stabiliser of~$\FF$,
then 
$$
U\cap B^{\times}= U_0 \cap B^{\times}.
$$
\end{lem}

%%%%%%%%%%%%%%%%%%%%% Iwahori decompositions %%%%%%%%%%%%%%%%%%%%%%%%%

\subsection{Iwahori decompositions}\label{Iwahori}

Put~$J=J(\beta, \curlya)$,~$J^1=J^1(\beta, \curlya)$
and~$H^1=H^1(\beta, \curlya)$ (see~\cite{bk}~\S3 for the definitions
of these groups). By~\cite{ihes}~Example~10.9,~$J^1$ and~$H^1$ have
Iwahori decompositions with respect to~$(M_0,P_0)$: 
\begin{eqnarray*}
&J^1= (J^1\cap U_0^{-}) (J^1 \cap M_0) (J^1 \cap U_0)& \\
&H^1= (H^1\cap U_0^{-}) (H^1 \cap M_0) (H^1 \cap U_0)&
\end{eqnarray*}
It will also be useful for us to form the group
$$ 
(J^1\cap P_0)H^1= (H^1\cap U_0^{-}) (J^1 \cap M_0) (J^1 \cap U_0).
$$
Now let~$\FF$ be a maximal~$F$-flag refining~$\FF_0$ and~$U$ the
corresponding unipotent subgroup. We want to understand the
group~$(J\cap U)H^1$. For~$1\le i \le r$, let~$\FF_i$ be the
maximal~$F$-flag in~$Ew_i$ given by intersection of~$\FF$
with~$Ew_i$. Let~$U_i$ be the unipotent radical of
the~$G_i$-stabiliser of~$\FF_i$; then~$U_i= U\cap G_i$ and
$$
U\cap M_0 =\prod_{i=1}^r U_i.
$$ 

For~$1\le i\le r$, we denote by~$\curlya_i$ the hereditary~$\oE$-order
in~$A_i=\End_F(Ew_i)$ given by the lattice chain~$\LL_i=\{L_k\cap
Ew_i:k\in\mathbb Z\}$. 

\begin{lem} \label{J1UH}~$(J^1\cap U)H^1$ has an Iwahori
decomposition with respect~$(M_0,P_0)$ and
$$
((J^1\cap U)H^1)\cap M_0 = 
\prod_{i=1}^r (J^1(\beta, \curlya_i)\cap U_i)H^1(\beta, \curlya_i).
$$
\end{lem}

\begin{proof} Since~$J^1$ and~$U$ have Iwahori decompositions with
respect to~$(M_0,P_0)$, so does~$J^1\cap U$. Since this group
normalises~$H^1$, which also has an Iwahori decomposition, we see
that~$(J^1\cap U)H^1$ has an Iwahori decomposition with
respect~$(M_0,P_0)$ and, in particular,
$$ 
((J^1\cap U)H^1)\cap M_0= (J^1\cap U \cap M_0) (H^1 \cap M_0).
$$
The lemma now follows from the decompositions of~$J^1\cap M_0$
and~$H^1\cap M_0$ (see~\cite{ihes}~\S10) and the decompositions
above.
\end{proof}

Since~$\dim_E (E w_i) =1$, the algebra~$E=F[\beta]$ is a maximal
subfield of~$A_i$. We will use the lemma above as a reduction step to
the case when~$E$ is a maximal subfield of~$A$.

\begin{lem}\label{likeit} Let~$U'$ be a unipotent radical of some
parabolic subgroup of~$G$, then the image of 
$$
U'\cap J \rightarrow J/J^1 \cong \ubz{}/\ub{1}{} 
%\cong \prod_{i=1}^{\eBoe{}}\Aut_{\kE}(L_{i-1}/ L_{i})
$$
is contained in the unipotent radical of some Borel subgroup
of~$\ubz{}/\ub{1}{}$.
\end{lem} 

\begin{cor}\label{dec} We have
$$ 
J\cap U = (\ubz{}\cap U_0)(J^{1}\cap U).
$$
Moreover,~$(J\cap U)H^1$ has an Iwahori decomposition with respect
to~$(M_0,P_0)$.
\end{cor} 

\begin{proof} The~$E$-basis~$\{w_1, \ldots, w_r\}$ of~$V$
identifies~$B^{\times}$ with~$\GL_r(E)$ and~$\ubz{}\cap U_0$ with a
subgroup of unipotent upper-triangular matrices with entries
in~$\oE$. This implies that the image of~$\ubz{}\cap U_0$
in~$\ubz{}/\ub{1}{}$ is the unipotent radical of some Borel subgroup
of~$\ubz{}/\ub{1}{}$. Lemma~\ref{likeit} implies that  
$$ 
(J\cap U)J^1=  (\ubz{}\cap U_0)J^1.
$$
Intersecting both sides with~$U$ gives the first part of the
lemma. The second follows immediately from Lemma~\ref{J1UH}
\end{proof}

%%%%%%%%%%%%%%%%%%%%%%%%%% $\beta$-extensions  %%%%%%%%%%%%%%%%%%%%%%%

\subsection{$\beta$-extensions}

Let~$\mathcal C (\curlya, 0, \beta)$ be the set of simple characters
of~$H^1$, in the sense of~\cite{bk}~(3.2.3). Let~$\theta\in \mathcal C
(\curlya, 0, \beta)$; then~$\theta$ is a linear character and, there
exists a unique irreducible representation~$\eta$ of~$J^1$
containing~$\theta$,~\cite{bk}~(5.1.1).

\begin{lem}\label{etaindone} Let~$\mathbf{k}_{\theta}$ be the
nondegenerate alternating form on~$J^1/H^1$ given by 
$$
\mathbf{k}_{\theta}(x,y)=\theta([x,y]),\qquad\hbox{for }x,y\in J^1,
$$
introduced in~\cite{bk}~\S3.4. Let~$\mathcal U$ be a subgroup of~$J^1$
containing~$H^1$ and let~$\overline{\mathcal U}$ be the image
of~$\mathcal U$ in~$J^1/H^1$. Suppose there exists a linear
character~$\chi$ of~$\mathcal U$ such that~$\chi|_{H^1}= \theta$. Then
the following are
equivalent:
\begin{enumerate}
\item[(i)]~$\chi$ occurs in~$\eta$ with multiplicity one;
\item[(ii)]~$\overline{\mathcal U}$ is a maximal totally isotropic
subspace of~$J^1/H^1$ for the form~$\mathbf{k}_{\theta}$; 
\item[(iii)]~$\eta\cong \Indu{\mathcal U}{J^1}{\chi}$.
\end{enumerate}
\end{lem}  

\begin{proof} Since~$H^1$ is normal in~$J^1$
and~$J^1/H^1$ is abelian (by~\cite{bk}~(3.1.15)),~$\mathcal U$ is a
normal subgroup of~$J^1$. Moreover, since~$\chi|_{H^1} =\theta$
and~$J^1/H^1$ is abelian, the commutator subgroup of~$\mathcal U$ will
lie in the kernel of~$\theta$ and hence~$\overline{\mathcal U}$ is a
totally isotropic subspace of~$J^1/H^1$ for the
form~$\mathbf{k}_{\theta}$.

(i)$\Rightarrow$(ii) Let~$\overline{\mathcal U}_{max}$ be a maximal
totally isotropic subspace of~$J^1/H^1$ containing~$\overline{\mathcal
U}$ and let~$\mathcal U_{max}$ be its inverse image in~$J^1$, so
that~$\mathcal U_{max}/\Ker(\theta)$ is a maximal abelian subgroup
of~$J^1/\Ker(\theta)$. The character~$\chi$ admits extension to a
linear character of~${\mathcal U}_{max}$ in exactly~$(\mathcal
U_{max}:\mathcal U)$ ways and every one of these extensions occurs
in~$\eta$. Then~$\chi$ occurs in~$\eta$ with multiplicity at least
this index so that~$\mathcal U_{max}=\mathcal U$.

(ii)$\Rightarrow$(iii) Suppose that~$j\in J^1$ intertwines~$\chi$ with
itself. Let~$\mathcal U'$ be the subgroup of~$J^1$ generated by~$j$
and~$\mathcal U$ and let~$\overline {\mathcal U'}$ be the image
of~$\mathcal U'$ in~$J^1/H^1$. A typical element of~$\overline
{\mathcal U'}$ is a coset~$j^{a} x H^1$, where~$a$ is an integer
and~$x\in \mathcal U$. Since
$$
\chi([j^a x, j^b y])= \chi(j^a x j^{-a}) \chi(j^{a+b} y x^{-1}
j^{-a-b}) \chi(j^b y^{-1} j^{-b})= \chi([x,y])=1
$$
for all~$x, y \in \mathcal U$, the subspace~$\overline {\mathcal U'}$
is totally isotropic for the
form~$\mathbf{k}_{\theta}$. Hence~$\overline {\mathcal U'}=\overline
{\mathcal U}$ and~$j\in\mathcal U$. In particular,~$\Indu{\mathcal
  U}{J^1}{\chi}$ is irreducible. Since~$\eta$ is the unique
irreducible representation of~$J^1$ containing~$\theta$ and~$\chi$
contains~$\theta$, we have~$\eta\cong\Indu{\mathcal U}{J^1}{\chi}$.

(iii)$\Rightarrow$(i) Since~$\eta$ is irreducible, this is just
Frobenius reciprocity.
\end{proof}

Now let~$\kappa$ be a representation of~$J$, such
that~$\kappa|_{J^1}\cong \eta$ and~$\kappa$ is intertwined by the
whole of~$B^{\times}_{\beta}$, that is, a~$\beta$-extension of~$\eta$
in the sense of~\cite{bk}~(5.2.1).

\begin{thm}\label{Psikappa} Let~$[\curlya, n, 0, \beta]$ be a simple
stratum. Let~$\theta\in \mathcal C (\curlya, 0, \beta)$ and
let~$\kappa$ be a~$\beta$-extension of~$\eta$ as above. Let~$\FF$ be
a maximal~$F$-flag in~$V$, let~$U$ be the unipotent radical of
the~$G$-stabiliser of~$\FF$ and let~$\chi$ be a smooth character
of~$U$, such that  
\begin{enumerate}
 \item $\FF_0 \subseteq \FF$,
 \item $\theta|_{U\cap H^1}= \chi|_{U\cap H^1},$
 \item $\chi$ is trivial on~$U_0$.
\end{enumerate}
Let~$\Theta$ be the linear character of~$(J\cap U)H^1$ 
defined by 
$$
\Theta( u h) =\chi(u) \theta(h), \quad \forall u\in U\cap J, \quad
\forall h\in H^1.
$$
Then
$$
\kappa|_{(J\cap U)J^1}\cong \Indu{(J\cap U)H^1}{(J\cap U)J^1}{\Theta}.
$$
\end{thm}

Before proving Theorem~\ref{Psikappa}, we remark that it is not clear
that there exist a flag~$\FF$ and a character~$\chi$ satisfying the
hypotheses. This would follow from \cite{bh}~Lemma~2.10, if we could
fix the proof of that result. Instead, we will have to wait for
Theorem~\ref{wavingit} to see that there are indeed such a flag and
character.

\begin{proof} We begin by proving
\begin{equation} \label{Psieta}
\Indu{(J^1\cap U)H^1}{J^1}{\Psi}\cong \eta.
\end{equation}
\textit{Step 1.} We will prove~\eqref{Psieta} in  the special case
when~$E=F[\beta]$ is a maximal subfield of~$A$.

Suppose that~$E$ is a maximal subfield of~$A$ so that~$B=E$
and~$J=\oE^{\times}J^1$. The pair~$(J,\kappa)$ is a simple type in the
sense of~\cite{bk}~(5.5.10) and, since~$\curlyb{}=\oE$, we
have~$\eBoe{}=1$. Now~\cite{bk}~(6.2.2) and~(6.2.3) imply that there
exists a representation~$\Lambda$ of~$E^{\times}J$ such
that~$\Lambda|_{J} \cong \kappa$
and~$\pi=\cIndu{E^{\times}J}{G}{\Lambda}$ is an irreducible
supercuspidal representation of~$G$. If~$u\in U$, then~$\det_A(u)=1$
and this implies that
$$
(E^{\times}J)\cap U = J\cap U.
$$
Since~$\dim_E V=1$ the unipotent radical~$U_0$ is trivial and hence
Corollary~\ref{dec} implies that 
$$
J\cap U = J^1\cap U.
$$
Hence~$\Lambda|_{(E^{\times}J)\cap U}\cong \eta|_{J^1\cap
U}$. Since~$\eta$ is the unique irreducible representation
containing~$\theta$ and~$\Theta|_{H^1} = \theta$, we obtain
that~$\Theta$ occurs in~$\eta|_{(J^1\cap U)H^1}$. 
Since~$\Theta|_{U\cap J}= \chi|_{U\cap J}$, we obtain 
\begin{eqnarray*}
 1\ \le\ \dim\Hom_{(J^1\cap U)H^1}(\Theta,\eta) &\le &  
\dim\Hom_{J^1\cap U}(\chi,\eta) \\   &=  &
\dim\Hom_{(E^\times J^1)\cap U}(\chi,\Lambda)\ \le\ 1,
\end{eqnarray*}
where the last inequality follows from~\cite{bh}~Proposition~1.6(iii). Hence 
$$
\dim\Hom_{(J^1\cap U) H^1}(\Theta,\eta)=1.
$$  
The equivalence~\eqref{Psieta} now follows immediately from
Lemma~\ref{etaindone} applied to~$\mathcal U=(J^1\cap U)H^1$
and~$\chi=\Theta$.

\textit{Step 2.} We will prove~\eqref{Psieta} in  the general case by
reducing to Step 1.

 For~$1\le j \le r$ we have an equality of sets: 
$$ 
\{L_i\cap E w_j : L_i \in \LL \}= \{\pE^k w_j : k \in \mathbb Z \}.
$$
This follows from the explicit description of lattices in~$\LL$, in
terms of the~$E$-basis~$\{w_1,\ldots w_r\}$
of~$V$. Hence~\cite{ihes}~Example~10.9 implies that the character~$\theta$
is trivial on its restrictions to~$H^1\cap U_0$ and~$H^1\cap
U_0^-$. Moreover, by~\cite{ihes}~page~167, the restriction of~$\theta$
to~${H^1\cap G_i}=H^1(\beta,\curlya_i)$ is the simple
character~$\theta_i$ in~$\mathcal C(\curlya_i,0,\beta)$ corresponding
to~$\theta$ under the canonical
bijection~$\tau_{\curlya,\curlya_i,\beta}$ of~\cite{bk}~\S3.6,
where~$\curlya_i$ is the hereditary~$\oF$-order corresponding to the
lattice chain $\{\pE^k w_i:k\in\ZZ\}$ in~$Ew_i$. 
%%
%%  we cannot refer to \cite{bk} (7.2.3) directly since (7.1.11) does not hold.
%% since for (7.1.11) you need (7.1.1) and the ``t in (7.1.1)= our
%% r''. Now the comment straight after (7.1.1) implies that~$r$ must
%% divide~$e(\urlyb | \oE)$, which can happen if and only
%% if $r=e(\urlyb|\oE)$, i.e. $\urlyb=\urlyb_m$. But the minimal case
%% is also a special case of \cite{ihes} (10.9). 
%%
In particular,~\cite{ihes}~Example~10.9 implies that the analogue
of~\cite{bk}~(7.2.3) holds in our situation:
\begin{enumerate}
 \item\label{one} the subspaces~$(J^1\cap U_0)/(H^1\cap U_0)$
and~$(J^1\cap U_0^-)/(H^1\cap U_0^-)$ of~$J^1/H^1$ are both totally
isotropic for the form~$\mathbf{k}_{\theta}$, and orthogonal to the
subspace~$(J^1\cap M_0)/(H^1\cap M_0)$;  
 \item\label{two} the restriction of~$\mathbf{k}_{\theta}$ to the group
$$
(J^1\cap M_0)/(H^1\cap M_0) = \prod_{i=1}^r
J^1(\beta,\curlya_i)/H^1(\beta,\curlya_i) 
$$
is the orthogonal sum of the pairings~$\mathbf{k}_{\theta_i}$; 
\item\label{three} we have an orthogonal sum decomposition
$$
\frac{J^1}{H^1} = \frac{J^1\cap M_0}{H^1\cap M_0} \perp 
\left(\frac{J^1\cap U_0^-}{H^1\cap U_0^-} \times 
\frac{J^1\cap U_0}{H^1\cap U_0}\right).
$$
In particular, the restriction of~$\mathbf{k}_{\theta}$ to the
group~$(J^1\cap U_0^-) / (H^1\cap U_0^-) \times (J^1\cap U_0) /
(H^1\cap U_0)$ is non-degenerate.
\end{enumerate}

Let~$\overline{((J^1\cap U)H^1)\cap M_0}$ be the image of the natural
homomorphism 
$$
((J^1\cap U)H^1)\cap M_0\rightarrow (J^1\cap M_0)/(H^1\cap M_0).
$$
Lemma~\ref{J1UH} and~\eqref{two} above imply that 
$$
\overline{((J^1\cap U)H^1)\cap M_0}=\prod_{i=1}^r \overline{J^1(\beta,
  \curlya_i)\cap U_i}
$$
where~$\overline{J^1(\beta, \curlya_i)\cap U_i}$ is the image of the
natural homomorphism  
$$
J^1(\beta, \curlya_i)\cap U_i \rightarrow
J^1(\beta,\curlya_i)/H^1(\beta,\curlya_i).
$$ 
Since~$E$ is a maximal subfield of~$A_i$ and we have
proved~\eqref{Psieta} when~$E$ is maximal, $\overline{J^1(\beta,
\curlya_i)\cap U_i}$ is a maximal isotropic subspace
in~$J^1(\beta,\curlya_i)/H^1(\beta,\curlya_i)$ for the
form~$\mathbf{k}_{\theta_i}$. 

Now~\eqref{two} implies that~$\overline{((J^1\cap U)H^1)\cap M_0}$ is a
maximal isotropic subspace of~$(J^1\cap M_0)/(H^1\cap M_0)$
for~$\mathbf{k}_{\theta}$. Moreover,~\eqref{one} and~\eqref{three}
imply that~$(J^1\cap U_0)/(H^1 \cap U_0)$ is a maximal isotropic
subspace of~$(J^1\cap U_0^-) / (H^1\cap U_0^-) \times (J^1\cap U_0) /
(H^1\cap U_0)$. 

It follows from the orthogonal sum decomposition in~\eqref{three} that 
$$ 
\overline{((J^1\cap U)H^1)\cap M_0}\times (J^1\cap U_0) / (H^1\cap
U_0)
$$
is a maximal isotropic subspace in~$J^1/H^1$ for the
form~$\mathbf{k}_{\theta}$.

Since~$(J^1\cap U)H^1$ contains~$J^1\cap U_0$, the image of~$(J^1\cap
U)H^1$ in~$J^1/H^1$ contains a maximal totally isotropic subspace
of~$J^1/H^1$ described above. Since there exists a linear
character~$\Theta$ of~$(J^1\cap U)H^1$ extending $\theta$, we see that
this image is itself isotropic and hence must be a maximal totally
isotropic subspace. The equivalence~\eqref{Psieta} now follows from
Lemma~\ref{etaindone}. 

\textit{Step 3.} Finally, we will deduce Theorem~\ref{Psikappa}
from~\eqref{Psieta} by examining the construction of~$\kappa$
in~\cite{bk}~\S5.

Let~$\LL_m$ be the~$\oE$-lattice chain in~$V$ given by
$$ 
\LL_m = \{ \varpi_{E}^k( \oE w_1 + \cdots + \oE w_j + \pE
w_{j+1}+\cdots+ \pE w_r) : k\in \mathbb Z, 1\le j \le r \}
$$  
%\marginpar{\tiny Should~$nr$ be~$nr/e(\urlyb|\oE)$ here? Or are we
%  already only in the maximal case?}
Let~$\curlyb{m}=\End_{\oE}^0(\LL_m)$ so that~$\curlyb{m}$ is a
minimal~$\oE$-order in~$B$. Similarly,
let~$\curlya=\End_{\oF}^{0}(\LL_m)$ so that~$\curlya_m$ is  the unique
hereditary~$\oF$-order in~$A$ normalised by~$E^\times$ such
that~$\curlya_m\cap B=\curlyb{m}$. Moreover,~$[\curlya_m,nr/e,0,\beta]$, where 
$e=e(\urlyb|\oE)$, 
is a simple stratum in~$A$ and the groups
$$
H_m^1\ =\ H^1(\beta,\curlya_m),\qquad J_m^1\ =\ J^1(\beta,\curlya_m)
$$
have Iwahori decompositions with respect to~$(M_0,P_0)$. 

We denote by~$\theta_m$ the simple character in~$\mathcal
C(\curlya_m,0,\beta)$ corresponding to~$\theta$ via the canonical
bijection~$\tau_{\curlya,\curlya_m,\beta}$ of
\cite{bk}~\S3.6. Then~$\theta_m$ is trivial on~$H_m^1\cap U_0$
and~$H_m^1\cap U_0^-$ (see~\cite{ihes} (10.9)). We also denote
by~$\eta_m$ the unique irreducible representation of~$J_m^1$ which
contains~$\theta_m$.

Since~$\LL \subseteq \LL_m$ we have~$\curlya_m \subseteq \curlya$
and~$\curlyb{m}\subseteq \urlyb$. Moreover,  
$$
\ubz{}\cap U_0=\ubz{m}\cap U_0
$$ 
since with respect to the~$E$-basis~$\{w_1, \ldots, w_r\}$ of~$V$,
both groups are identified with unipotent upper-triangular matrices,
with entries in~$\oE$.

We define 
$$
\tilde\eta\ :=\ \Indu{(J\cap U)H^1}{(J\cap U)J^1}{\Theta}.
$$
Note first that~$\tilde\eta|_{J^1}=\Indu{(J^1\cap
U)H^1}{J^1}{\Theta}\cong \eta$ so~$\tilde\eta$ is certainly
irreducible. Moreover,~\cite{bk}~(5.1.1) implies that
the~$\mathbf{U}^1(\curlya_m)$-intertwining of~$\tilde\eta$ is
contained in
$$
I_{\mathbf{U}^1(\curlya_m)}(\eta) = (J^1 B^\times J^1)\cap
\mathbf{U}^1(\curlya_m) = \mathbf{U}^1(\curlyb{m})J^1 
= (J\cap U)J^1
$$
where the last equality follows from
Corollary~\ref{dec}. Hence~$\Indu{(J\cap
  U)J^1}{\mathbf{U}^1(\curlya_m)}{\tilde\eta}$ is irreducible.

Now we claim that
$$
\Indu{(J\cap U)J^1}{\mathbf{U}^1(\curlya_m)}{\tilde\eta} \cong
\Indu{J^1_m}{\mathbf{U}^1(\curlya_m)}{\eta_m}.
$$
It is enough to show that 
$$
\Hom_{\mathbf{U}^1(\curlya_m)}\left(\Indu{(J\cap
  U)J^1}{\mathbf{U}^1(\curlya_m)}{\tilde\eta}, 
\Indu{H^1_m}{\mathbf{U}^1(\curlya_m)}{\theta_m}\right) \ne 0,
$$
since the latter is a multiple
of~$\Indu{J^1_m}{\mathbf{U}^1(\curlya_m)}{\eta_m}$, which is
irreducible by the same argument as above. By Mackey Theory, it is
enough to show that
$$
\Hom_{((J\cap U)H^1)\cap H_m^1}(\Theta,\theta_m)\ne 0.
$$
But~$((J\cap U)H^1)\cap H_m^1$ has an Iwahori decomposition with respect
to~$(M_0,P_0)$, since both~$(J\cap U)H^1$ and~$H_m^1$
do. Moreover,~$\Theta$ is trivial on~$((J\cap U)H^1)\cap U_0^-$,
since~$((J\cap U)H^1)\cap U_0^-= H^1\cap U_0^-$ and~$\theta$ is
trivial on~$H^1 \cap U_0^-$. Since~$\chi$ is trivial on~$U_0$ we get
that~$\Theta$ is trivial  on~$((J\cap U)H^1)\cap U_0$. Hence
both~$\Theta$ and~$\theta_m$ are trivial on the subgroups~$((J\cap U)
H^1)\cap H^1_m \cap U_0^{-}$ and~$((J\cap U) H^1)\cap H^1_m \cap
U_0$. Finally, by~\cite{ihes}~\S10,
$$
H_m^1\cap M_0\ =\ \prod_{i=1}^r H^1(\beta,\curlya_i)\ =\ H^1\cap M_0
$$
so
$$
(J\cap U)H^1\cap H_m^1\cap M_0\ =\ H^1\cap M_0\ \subset\ H^1\cap H^1_m,
$$
where~$\Theta=\theta$ and~$\theta_m$ agree by~\cite{bk}~(3.6.1).

In particular, we have shown that, upto equivalence,~$\tilde\eta$
satisfies the conditions of the representation (also denoted
$\tilde\eta$) in~\cite{bk}~(5.1.15). Since these conditions uniquely
determine~$\tilde\eta$, we conclude that our~$\tilde\eta$ is, up to
equivalence, the same one used in the construction of~$\kappa$ in
\cite{bk}~(5.2.4) and hence 
$$
\kappa|_{(J\cap U)J^1} \cong \tilde\eta \cong \Indu{(J\cap U)H^1}{(J\cap
U)J^1}{\Theta} 
$$
as required.
\end{proof}

%%%%%%%%%%%%%%%%%%%%%%%%%%%%%%%%%%%%%%%%%%%%%%%%%%%%%%%%%%%%%%%%%%%%%%
%%%%%%%%%%%%%%%%%%% A particular unipotent subgroup %%%%%%%%%%%%%%%%%%
%%%%%%%%%%%%%%%%%%%%%%%%%%%%%%%%%%%%%%%%%%%%%%%%%%%%%%%%%%%%%%%%%%%%%%

\section{A particular unipotent subgroup}\label{Udef}

The main result of this section says that we may choose a
maximal~$F$-flag~$\FF$ in~$V$ and a character~$\chi$ of
the~$G$-stabiliser~$U$ of~$\FF$ which satisfy the hypotheses of
Theorem~\ref{Psikappa} and, moreover, such that we can control the
restriction of certain characters of~$U$ to~$U\cap
\UU(\urlyb)$ (see Theorem~\ref{wavingit} for details). We will
continue with the notation of the previous section but, since we will
consider only supercuspidal representations in the applications, we will
assume that~$e(\urlyb|\oE)=1$.

%%%%%%%%%%%%%%%%%%%%%%%%%% An~$F$-basis of~$E$ %%%%%%%%%%%%%%%%%%%%%%%

\subsection{An~$F$-basis of~$E$}

Although we would like to work in a basis-free way, we are forced to 
choose one, since the zeta function in \S\ref{epsi}, whose functional equation 
defines $\varepsilon$-factors of pairs, is defined on matrices. In this
subsection we  consider the special case when~$E$ is maximal, so we may
identify~$V=E$. We show that there exists an~$F$-basis $\BB$ of~$E$ with 'nice' 
properties. The basis is chosen  to ease the pain of calculations in \S7, 
see also Corollary \ref{goodbasis}. The conditions imposed on $\BB$ imply certain 
uniqueness result, which will be used to define a 
numerical invariant in \S\ref{numerical}. 

Let~$\FF=\{V_i: 1\le i \le d\}$, where~$d=[E:F]$, be a maximal~$F$-flag
in~$E$, let~$U$ be the unipotent radical of the~$\Aut_F(E)$-stabiliser
of~$\FF$, and let~$\chi$ be a smooth, non-degenerate character of~$U$. 
Let $\psi_E$ be an additive character of~$E$, trivial on~$\pE$, and such that 
$$ 
\psi_E(x)= \psi_F(x), \quad \forall x\in F.
$$ 

\begin{defi}\label{balba} An~$F$-basis~$\BB=\{ x_1, \ldots, x_d\}$
of~$E$ is~\emph{$(U, \chi, \psi_E)$-balanced} if the  following hold: 
\begin{enumerate}
\item $V_i= Fx_1+\cdots +F x_i$,~$1\le i \le d$;
\item there exists functions~$a_{i}: \ZZ\rightarrow \ZZ$, 
for~$1\le i \le d$ such that
$$ 
\pE^k= \sum_{i=1}^d \pF^{a_{i}(k)} x_{i}, \quad \forall k\in \ZZ;
$$ 
\item if~$u\in U$, and~$(u_{ij})$ is the matrix of~$u$ with respect
to~$\BB$, then  
$$
\chi(u)=\psi_F(\sum_{i=1}^{d-1} u_{i, i+1});
$$
\item if we embed naturally~$E\hookrightarrow \End_F(E)$ and,
for~$\xi\in E$, we let $(\xi_{ij})$ be the matrix of~$\xi$ with
respect to~$\BB$, then
$$
\psi_E(\xi)=\psi_F(\xi_{dd}),\quad  \forall \xi\in E.
$$
\end{enumerate}
\end{defi} 

\begin{prop}\label{balance} There exists an~$(U,\chi,\psi_E)$-balanced
basis. Moreover, if~$\BB=\{x_1, \ldots x_d\}$
is~$(U,\chi,\psi_E)$-balanced and~$\BB'=\{y_1, \ldots, y_d\}$ is a
basis of~$E$ over~$F$ which satisfies
Definition~\ref{balba}~(i),(iii) and~(iv), then~$x_d x_1^{-1}= y_d
y_1^{-1}$.
\end{prop}

\begin{proof} 
Proposition~II-3 in~\cite{weil} applied to~$\FF$
and~$\{\pE^i:i\in\ZZ\}$ implies that there exists
an~$F$-basis~$\BB=\{x_{1}, \ldots, x_{d}\}$ of~$E$ which satisfies~(i)
and~(ii). Since~$\chi$ is non-degenerate, after replacing~$x_i$ by
some~$\lambda_i x_i$, where~$\lambda_i \in F^{\times}$, we may ensure
that~$\BB$ satisfies~(iii).

For~$\xi\in E$, let~$(\xi_{ij})$ be a matrix of~$\xi$ with respect
to~$\BB$. Consider the function~$\phi_{\BB}: E\rightarrow
\CC^{\times}$ given by~$\phi_{\BB}(\xi)=\psi_F(\xi_{dd})$. It is clear
that~$\phi_{\BB}$ is an additive character. If~$x\in F$,
then~$\phi_{\BB}(x)=\psi_F(x)$, hence~$\phi_{\BB}$ is non-trivial
on~$\oE$. Since~$x_d\in \pE^{v_E(x_d)}$ and~$x_d\not\in
\pE^{v_{E}(x_d)+1}$ we have~$a_d(v_{E}(x_d))=0$
and~$a_d(v_{E}(x_d)+1)=1$, hence if~$\xi\in \pE$, then~$\xi_{dd}\in
\pF$, and so~$\phi_{\BB}$ is trivial on~$\pE$. Since~$\psi_E$
and~$\phi_{\BB}$ have the same conductor, there exists~$\alpha \in
\oE^{\times}$, such that~$\psi_E(\xi)=\phi_{\BB}(\alpha \xi)$, for
all~$\xi\in E$.

Set~$A(E)=\End_F(E)$, let~$\curlya(E)$ be the hereditary order
corresponding to the lattice chain~$\LL(E)=\{\pE^i:i\in\ZZ\}$, and,
for~$1\le i,j\le d$, let~$\Eins_{ij}\in A(E)$ be given
by~$\Eins_{ij}x_k=\delta_{ik}x_j$, where~$\delta_{ik}$ is the
Kronecker delta. Since~$\BB$ satisfies~(ii) we have~$\Eins_{ij}\in
\curlya$.

Set~$u= 1+ (\alpha-1)\Eins_{dd}\in  \curlya(E)$. We have 
$$
\psi_F(\lambda)=\psi_E(\lambda)=  \phi_{\BB}(\alpha \lambda)=
\psi_F(\alpha_{dd} \lambda), \quad \forall \lambda\in F.
$$ 
Hence~$\alpha_{dd}=1$ and so~$u\in U\cap \curlya(E)$. So the basis~$u
\BB= \{ u x_1, \ldots, u x_d\}$ also satisfies~(i),(ii)
and~(iii). Now
$$
\phi_{u\BB}(\xi)= \phi_{\BB}( u^{-1}\xi u)= \psi_F\left(\sum_{i=1}^d
\xi_{di} \alpha_{id}\right)= \psi_F((\xi\alpha)_{dd})=\phi_{\BB}(\xi
\alpha)=\psi_E(\xi).
$$
Hence~$u\BB$ is~$(U,\chi, \psi_E)$-balanced. 

%Set~$A=\End_F(E)$ and let~$\curlya$ be a hereditary order
%corresponding to the lattice chain~$\{\pE^i:i\in\ZZ\}$ and
%let~$\Eins_{ij}\in A$,~$\Eins_{ij}x_k= 0$, if~$k\neq i$,
%and~$\Eins_{ij} x_i= x_j$. Let~$\psi_A: A\rightarrow  \CC^{\times}$
%be an additive character, given by~$\psi_A(a)=
%\psi_F(\tr_{A/F}(a))$. According to~\cite{bk} (1.3.4) there exists
%an~$(E,E)$-bimodule homomorphism~$s: A\rightarrow E$, such that 
%$$ 
%\psi_A(\xi a)= \phi_{\BB}(\xi s(a)), \quad \forall \xi\in E, \quad
%\forall a\in A.
%$$     
%Since~$\phi_{\BB}(\xi)=\psi_F(\xi_{dd})= \psi_A(\xi
%\Eins_{dd})=\phi_{\BB}(\xi s(\Eins_{dd}))$, for all~$\xi\in E$, we
%have~$s(\Eins_{dd})=1$. Since~$\psi_E$ and~$\phi_{\BB}$ have the same
%conductor, there exists~$\alpha \in \oE^{\times}$, such
%that~$\psi_E(\xi)=\phi_{\BB}(\alpha \xi)$, for all~$\xi\in
%E$. Set~$u= 1+ (\alpha-1)\Eins_{dd}$.  Since~$\BB$ satisfies (ii) we
%have~$\Eins_{dd}\in \curlya$. Hence,~$u\in  \curlya$. We have 
%$$
%\psi_F(\lambda)=\psi_E(\lambda)=  \phi_{\BB}(\alpha \lambda)=
%\psi_F(\alpha_{dd} \lambda), \quad \forall \lambda\in F.
%$$ 
%Hence,~$\alpha_{dd}=1$ and so~$u\in U\cap \curlya$. So the basis~$u
%\BB= \{ u x_1, \ldots, u x_d\}$ also satisfies (i), (ii) and
%(iii). Now 
%$$
%\phi_{u\BB}(\xi)= \phi_{\BB}( u^{-1}\xi u)= \psi_F(\sum_{i=1}^d
%\xi_{di} \alpha_{id})= \psi_F((\xi\alpha)_{dd})=\phi_{\BB}(\xi
%\alpha)=\psi_E(\xi).
%$$
%Hence~$u\BB$ is~$(U,\chi, \psi_E)$-balanced. 

Suppose that~$\BB=\{x_1, \ldots x_d\}$ is~$(U,\chi,\psi_E)$-balanced
and that~$\BB'=\{y_1, \ldots, y_d\}$ satisfies 
Definition~\ref{balba}(i),(iii)
and~(iv). After translating by some~$\lambda\in F^{\times}$, we may
assume that~$x_1=y_1$. Parts~(i),(iii) imply that there exists~$u \in
U$ such that~$\BB'=u\BB$. Let~$(u_{ij})$ be the matrix of~$u$ with
respect to~$\BB$. Then, for all~$\xi\in E$, we have
$$
\psi_E(\xi)=\phi_{u\BB}(\xi)= \psi_F\left((u^{-1}\xi u)_{dd}\right)
=\psi_F\left(\sum_{i=1}^d \xi_{di}u_{id}\right) 
=\psi_A\left(\xi \left(\sum_{i=1}^d u_{id} 
\Eins_{di}\right)\right)\!. 
$$
where $\psi_A=\psi_F\circ\tr_{A(E)/F}$. 
Moreover, $\Eins_{id}=x_ix_d^{-1}\Eins_{dd}$ so
$$
\psi_E(\xi) =\psi_A\left(\xi \left(\sum_{i=1}^d u_{id} x_i
x_d^{-1}\right)\Eins_{dd}\right).
$$
By~\cite{bk}~(1.3.4) there exists an~$(E,E)$-bimodule
homomorphism~$s: A\rightarrow E$ (a \emph{tame corestriction}), such that 
$$ 
\psi_A(\xi a)= \phi_{\BB}(\xi s(a)), 
\quad \forall \xi\in E, \quad\forall a\in A.
$$     
Since~$\phi_{\BB}(\xi)=\psi_F(\xi_{dd})= \psi_A(\xi
\Eins_{dd})=\phi_{\BB}(\xi s(\Eins_{dd}))$, for all~$\xi\in E$, we
obtain that~$s(\Eins_{dd})=1$. Hence we get
$$
\psi_E(\xi)=\phi_{\BB}
\left(\xi\left(\sum_{i=1}^d u_{id} x_i x_d^{-1}\right)\right)
=\psi_E
\left(\xi \left(\sum_{i=1}^d u_{id} x_i x_d^{-1}\right)\right).
$$
Thus~$x_d= \sum_{i=1}^d u_{id} x_i$, which implies that~$u_{id}=0$,
if~$i\neq d$. In particular, we get~$y_d= u x_d = x_d$. 
\end{proof}

%%%%%%%%%%%%%%%%%%%%%%%% A choice of maximal flag %%%%%%%%%%%%%%%%%%%%

\subsection{A choice of maximal flag}  
 
Let~$[\curlya, n, 0, \beta]$ be a simple stratum, such
that~$e(\urlyb|\oE)=1$. Let~$\theta\in \mathcal C (\curlya, 0, \beta)$
and let~$\FF_0$ be a maximal~$E$-stable flag in~$V$,~$U_0$ the
unipotent radical of the~$G$-stabiliser of~$\FF_0$. 

\begin{thm}\label{wavingit} There exist a maximal~$F$-flag~$\FF$
in~$V$, a smooth character~$\chi$ of the unipotent radical~$U=U_{\FF}$
and an element~$b\in X_{\FF}\cap \curlya$ , such that the following
hold:
\begin{enumerate}
\item $\FF_0 \subseteq \FF$,
\item $\theta|_{U\cap H^1}= \chi|_{U\cap H^1},$
\item $\chi$ is trivial on~$U_0$,
\item the character
$$ 
\overline{\psi}_b: (\UU(\urlyb)\cap U)/(\UU^1(\urlyb)\cap U)
\rightarrow \mathbb C^{\times}, \quad u (\UU^1(\urlyb)\cap U)
\mapsto \psi_b(u)
$$
defines a non-degenerate character of a maximal unipotent subgroup
of~$\UU(\urlyb)/\UU^1(\urlyb)$.
\end{enumerate}
\end{thm}

\begin{proof} Let us consider first, as in Step 1 of the proof of
Theorem~\ref{Psikappa} the  case when~$E$ is a maximal subfield
of~$A$. We identify~$V=E$ and~$\curlya=\curlya(E)$ is the hereditary
order associated to the lattice chain $\LL(E)=\{\pE^k: k\in
\ZZ\}$. The parts~(i),(iii) and~(iv) are empty in this case
since~$\FF_0=\{ 0\subseteq V\}$,~$U_0 =\{1\}$,
and~$\UU(\urlyb)=\oE^{\times}$ does not contain non-trivial unipotent
elements, and so~$\UU(\urlyb)\cap U =\{1\}$. Consider the
supercuspidal representation~$\pi=\cIndu{E^{\times}J}{G}{\Lambda}$, as
in Step 1 of the proof of Theorem~\ref{Psikappa}. According
to~\cite{bh}~Proposition~1.6(i), there exist a maximal~$F$-flag~$\FF$
in~$E$, and a smooth character~$\chi$ of the unipotent radical~$U$ of
the~$G$-stabiliser of~$\FF$ such that 
$$ 
\Hom_{U\cap (E^{\times}J)}(\chi, \Lambda)\neq 0.
$$    
Since~$\Lambda|_{H^1}= (\dim\Lambda) \theta$, we obtain 
$$ 
\theta|_{U\cap H^1}= \chi|_{U\cap H^1}.
$$
Hence the theorem holds for~$E$ maximal. Note that, since~$E^{\times}$
normalises~$H^1$ and~$\theta=\theta^x$, for all~$x\in E^{\times}$, we
may replace the pair~$(U,\chi)$ by a conjugate~$(U^x, \chi^x)$,
where~$x\in E^{\times}$.

Now let us consider the general case. Put~$d=[E:F]$. 

\textit{Construction of~$\FF$ and~$\chi$.}  According 
to~\cite{weil}~Proposition~II-3 we may choose
an~$E$-basis~$\BB_E=\{w_1, \ldots, w_r\}$ of~$V$ such that 
$$ 
\FF_0=\{\sum_{i=1}^j E w_i: 1\le j \le r\} ,\quad \textrm{and}\quad
 L_k=\pE^k w_1+\cdots +\pE^k w_r, \quad \forall k \in \mathbb Z.
$$
Let~$P_0$ be the~$G$-stabiliser of~$\FF_0$. Moreover,
put~$G_i=\Aut_F(Ew_i)$, for~$1\le i\le r$, and
$$
M_0\ =\ \prod_{i=1}^r G_i,
$$
a Levi component of the parabolic subgroup~$P_0$ of~$G$. Then~$U_0$ is
the unipotent radical of~$P_0$, so that~$P_0=M_0U_0$. We are in the
situation considered in Step~2 of the proof of Theorem~\ref{Psikappa}. In
particular,
$$ 
H^1(\beta, \curlya)\cap M_0= \prod_{i=1}^r H^1(\beta, \curlya_i)\cong
\prod_{i=1}^r H^1(\beta, \curlya(E)),
$$
where the last isomorphism is induced by identifying~$E w_i$
with~$E$. Moreover, according to~\cite{ihes}~Corollary~10.16, there
exists~$\theta_F\in \mathcal C(\curlya(E), 0, \beta)$ such that
$$
\theta|_{H^1(\beta, \curlya)\cap M_0}= \theta_F\otimes\cdots \otimes
\theta_F
$$
via the above identification. Now, by the case when~$E$ is maximal
considered above, we know that there exist a maximal~$F$-flag~$\FF_1$
in~$E$, and a smooth character~$\chi_1: U_{\FF_1}\rightarrow
\CC^{\times}$ such that  
$$ 
\chi_1(u)= \theta_F(u), \quad  \forall u\in U_{\FF_1}\cap
H^1(\beta,\curlya(E)),
$$
where~$U_{\FF_1}$ is the unipotent radical of~$\Aut_F(E)$-stabiliser
of~$\FF_1$. Proposition~\ref{conju}~(i) implies that~$\chi_1$ is
non-degenerate. Choose an additive character~$\psi_E$ of~$E$, such
that~$\psi_E(x)=\psi_F(x)$, for all~$x\in F$, and~$\psi_E$ trivial
on~$\pE$. Proposition~\ref{balance} gives
a~$(U_{\FF_1},\chi_1,\psi_E)$-balanced~$F$-basis~$\BB_1=\{x_{11},
\ldots x_{1d}\}$ of~$E$. Set~$y=x_{1d} x_{11}^{-1}$ and  for~$2\le j
\le r$,  let~$\mathcal B_j=\{x_{j1},\ldots, x_{jd}\}$ be the basis
of~$E$ over~$F$ given by
$$
 x_{ji}= y^{j-1} x_{1i}, \quad 1\le i \le d.
$$
Note that, in particular,~$x_{j1}=x_{j-1,d}$, for~$2\le j\le r$.

Let~$\FF_j=\{ \sum_{i=1}^k F x_{ji}: 1\le k\le d\}$; then~$\FF_j= y
\FF_{j-1}$ and hence~$U_{\FF_j}= U_{\FF_{j-1}}^{y}$ so we may define 
a character~$\chi_j: U_{\FF_j}\rightarrow \CC^{\times}$ by~$\chi_j=
\chi_{j-1}^{y}$. Since~$y\in E^{\times}$ normalises~$\theta_F$, we
obtain that 
$$ 
\chi_j(u)= \theta_F(u), \quad \forall u\in U_{\FF_j}\cap
H^1(\beta,\curlya(E)), \quad 1\le j\le r.
$$ 
Let~$\FF=\{V_k: 1\le k\le N\}$ be the maximal~$F$-flag in~$V$ given by
$$ 
V_{(i-1)d +j}=F x_{11}w_1+\cdots+ F x_{ji}w_i,\quad 1\le i\le r, \quad
1\le j\le d,
$$
and let~$U_\FF$ be the unipotent radical of the~$G$-stabiliser
of~$\FF$. Since~$U_{\FF}/U_0\cong U_{\FF}\cap M_0\cong \prod_{i=1}^r
U_i$ we may define~$\chi:U_{\FF}\rightarrow \CC^{\times}$ by
$$ 
\chi|_{U_{\FF}\cap M_0}= \chi_1\otimes\cdots \otimes\chi_r,\quad
\chi|_{U_0}=\Eins.
$$  
The Iwahori decomposition implies that 
\begin{equation*}
\begin{split}
U_{\FF}\cap H^1(\beta,\curlya)=&\ \bigl( U_{\FF}\cap M_0\cap
H^1(\beta,\curlya)\bigr)(H^1(\beta,\curlya)\cap U_0) \\ 
\cong &\ \bigl(\prod_{i=1}^r U_{\FF_i} \cap H^1(\beta,
\curlya(E))\bigr) (H^1(\beta, \curlya)\cap U_0) 
\end{split}
\end{equation*}
Since~$\theta$ is trivial on~$U_0\cap H^1(\beta, \curlya)$, it follows
that 
$$
\theta(u)=\chi(u), \quad \forall u\in U\cap H^1(\beta, \curlya).
$$ 

\textit{Construction of~$b$.}
For~$\boldsymbol\mu=(\mu_1,...,\mu_{r-1})\in\oF^{r-1}$, we
define~$b=b(\boldsymbol\mu)\in X_{\FF}$ by 
$$
b(x_{ji}w_j) = 
\begin{cases} \mu_j x_{j+1,1}w_{j+1},
  &\quad\hbox{if }i=d,\ 1\le j\le r-1; \\
0 &\quad\hbox{otherwise.} \end{cases}
$$
We claim that such~$b$ also lies in~$\curlya$. For~$1\le j\le r$ we
have constructed a~$(U_{\FF_j},\chi_j,\psi_E)$-balanced
basis~$\BB_j=\{x_{j1},\ldots, x_{jd}\}$, of~$E$ over~$F$, so there
exist functions~$a_{ji}: \ZZ\rightarrow \ZZ$, for~$1\le i \le d$, such
that  
$$ 
\pE^k= \sum_{i=1}^d \pF^{a_{ji}(k)} x_{ji}, \quad \forall k\in \ZZ.
$$
Hence, 
$$ 
L_k=\sum_{j=1}^r \pE^k w_j=\sum_{j=1}^r \sum_{i=1}^d \pF^{a_{ji}(k)}
x_{ji}w_j, \quad \forall k\in \ZZ.
$$

Since, by construction,~$x_{jd}=x_{j+1,1}$, for~$1\le j \le r-1$, we
have~$a_{jd}(k)=a_{j+1,1}(k)$ for all~$ k \in \ZZ$. Hence 
$$
b L_k= \sum_{j=1}^{r-1}\mu_j\pF^{a_{jd}(k)}x_{j+1,1}w_{j+1}\subseteq
L_k,\quad \forall k\in \ZZ
$$
which implies that~$b\in \curlya$.

Since~$b\in X_{\FF}$, Lemma~\ref{charX} implies that~$\psi_b$ defines
a linear character of~$U$. Since~$b\in \curlya$ we have~$\psi_b(u)=1$
for all~$u\in \UU^1(\curlya)$. Hence the character~$\overline{\psi}_b$
is well defined. Moreover,~$(\UU(\urlyb)\cap U)/(\UU^1(\urlyb)\cap U)$
is the unipotent radical of the~$\Aut_{\kE}(L_0/L_1)$-stabiliser of
the flag~$\{\sum_{j=1}^i \kE (w_j+L_1): 0\le i \le r\}$. This follows
from Lemma~\ref{refine}.

Let~$\boldsymbol\mu=(1,...,1)$ and let~$b=b(\boldsymbol\mu)$ as
above. (Indeed, any~$\boldsymbol\mu\in(\oF^\times)^{r-1}$ would do.)
We claim that~$\overline{\psi}_b$ is non-degenerate. For~$2\le k \le
r$ and~$\xi\in E$ we define~$u_{\xi,k}\in U\cap B^{\times}=U_0\cap
B^{\times}$ by 
$$
u_{\xi,k}(w_j) = 
\begin{cases} w_{k}+\xi w_{k-1},&\quad\hbox{if }k=j; \\
w_j, &\quad\hbox{otherwise.} \end{cases}
$$  
We will compute~$(\chi\psi_b)(u_{\xi,k})=\psi_b(u_{\xi,k})$. Since~$x_{ji}\in E$, it commutes with~$u_{\xi,k}$
and hence
\begin{equation*}
 (u_{\xi,k}-1)(x_{ji}w_j)= \left\{ \begin{array}{ll} 
   \xi x_{ki}w_{k-1} & \textrm{if~$j=k~$,}\\
   0 & \textrm{otherwise.}
   \end {array}\right.
\end{equation*}    
Since~$\mathcal B_{k-1}$ is an~$F$-basis of~$E$ we may write uniquely
$$ 
\xi x_{ki}=\sum_{j=1}^d \lambda_{ji}^{(k)} x_{k-1,j}
$$
where~$\lambda_{ji}^{(k)}\in F$. Then
\begin{equation*}
 (b(u_{\xi,k}-1))(x_{ji}w_j)= \left\{ \begin{array}{ll} 
   \lambda^{(k)}_{di}x_{k1}w_{k} & \textrm{if~$j=k~$,}\\
   0 & \textrm{otherwise.}
   \end {array}\right.
\end{equation*}
Hence, 
$$ 
\trAF( b(u_k-1))=\lambda^{(k)}_{d1}.
$$
Since, by construction,~$x_{k1}=x_{k-1,d}$, we
have~$\lambda^{(k)}_{d1}=\xi^{(k)}_{dd}$, where~$(\xi^{(k)}_{ij})$ is
the matrix of~$\xi\in \End_F(E)$ with respect to~$\BB_{k-1}$. However,
since~$\BB_{k}$ is~$(U_{\FF_k},\chi_k,\psi_E)$-balanced
%= y^{k-1}\BB_1$
%, for~$k\ge 1$, and~$y\in E^{\times}$, 
%so commutes with~$\xi$. Hence,~$(\xi^{(k)}_{ij})=(\xi^{(1)}_{ij})$
%and since~$\BB_1$ was chosen to be~$(U_1,\chi_1,\psi_E)$-balanced, we
%get 
%$$ 
%\psi_b(u_{\xi,k})=\psi_F(\xi^{(1)}_{dd})=\psi_E(\xi).
%$$
$$ 
\psi_b(u_{\xi,k})=\psi_F(\xi^{(k)}_{dd})=\psi_E(\xi).
$$
Now,~$\psi_E$ has conductor~$\pE$, and so~$\overline{\psi}_b$ is
non-degenerate.
\end{proof}

We record the following corollary to the proof of the theorem.

\begin{cor}\label{goodbasis} Fix an additive
character~$\psi_E:E\rightarrow \CC^{\times}$ such
that~$\psi_E(x)=\psi_F(x)$, for all~$x\in F$ and $\psi_E$ is trivial
on~$\pE$. Choose an~$E$-basis~$\BB_E=\{w_1, \ldots, w_r\}$ of~$V$ such
that
$$ 
\FF_0=\{\sum_{i=1}^j E w_i: 1\le j \le r\} ,\quad \textrm{and}\quad
 L_k=\pE^k w_1+\cdots +\pE^k w_r, \quad \forall k \in \mathbb Z.
$$
Let~$\FF$,~$\chi$ and~$\psi_b$ be as in Theorem~\ref{wavingit}. There
exists a basis~$\{x_1, \ldots, x_d\}$ of~$E$ over~$F$, satisfying
Definition~\ref{balba}(iv) with respect to~$\psi_F$ and~$\psi_E$,
such that~$x_1=1$ and such that, if we putt  
$$ 
v_{d(i-1) +j} = x_d^{i-1} x_j w_i, \quad 1\le i\le r, \quad 1\le j \le
d,
$$
then the set~$\BB_F=\{v_1, \ldots, v_N\}$ is an~$F$-basis of~$V$ with
the following properties:
\begin{enumerate} 
\item $\FF=\{\sum_{i=1}^j F v_i : 1\le j\le N\};$
\item there exist functions~$a_{i}: \ZZ\rightarrow \ZZ$, for~$1\le i
\le N$ such that 
$$ 
L_k= \sum_{i=1}^N \pF^{a_{i}(k)} v_{i}, \quad \forall k\in \ZZ;
$$ 
\item if~$u\in U$ has matrix~$(u_{ij})$ with respect to~$\BB_F$, then 
$$
(\chi\psi_b)(u)=\psi_F(\sum_{i=1}^{N-1} u_{i, i+1});
$$
\item if~$u\in U\cap B^{\times}$ has matrix~$(u_{ij})$ with respect
to~$\BB_E$, then
$$
(\chi\psi_b)(u)= \psi_E(\sum_{i=1}^{r-1} u_{i,i+1}).
$$ 
\end{enumerate} 
\end{cor}

If we put~$\beta=0$, so that~$E=F$,~$\oE=\oF$ and~$r=N$ then the proof
of Theorem~\ref{wavingit} formally goes through
with~$\FF=\FF_0$,~$\mathcal B_j=\{1\}$, for~$1\le j \le N$,
and~$b=b(\boldsymbol\mu)$, for any~$\boldsymbol\mu \in
(\oF^{\times})^{N-1}$. We obtain a `level zero version' of
Theorem~\ref{wavingit}:

\begin{cor} \label{version0} Let~$\curlya$ be a maximal~$\oF$-order
in~$A$. Then there exist a maximal~$F$-flag~$\FF$ in~$V$ and an
element~$b\in X_{\FF}\cap \curlya$ such that  
$$ 
\overline{\psi}_b: (\UU(\curlya)\cap U)/(\UU^1(\curlya)\cap U)
\rightarrow \mathbb C^{\times}, \quad u (\UU^1(\curlya)\cap U)
\mapsto \psi_b(u)
$$
defines a non-degenerate character of a maximal unipotent subgroup of
the group~$\UU(\curlya)/\UU^1(\curlya)$.
\end{cor}

%%%%%%%%%%%%%%%%%%%%%%%%%%%%%%%%%%%%%%%%%%%%%%%%%%%%%%%%%%%%%%%%%%%%%%
%%%%%%%%%%%%%%%%%% Supercuspidal representations %%%%%%%%%%%%%%%%%%%%%
%%%%%%%%%%%%%%%%%%%%%%%%%%%%%%%%%%%%%%%%%%%%%%%%%%%%%%%%%%%%%%%%%%%%%%

\section{Supercuspidal representations}\label{super}

Let~$\pi$ be an irreducible supercuspidal representation
of~$G$. By~\cite{bk}~\S6,~$\pi$ is compactly induced: 
$$
\pi\cong\cIndu{\mathbf J}{G}{\Lambda},
$$
from a (rather special) open compact-mod-centre subgroup~$\mathbf J$
of~$G$. More  precisely,~$\mathbf J$ has a unique maximal compact open
subgroup~$J$ and if we put~$\lambda:=\Lambda|_J$ then~$(J,\lambda)$ is
a maximal simple type in the sense of~\cite{bk}~\S6.

\begin{defi}\label{maxtype} The pair~$(J, \lambda)$ is a \emph{maximal
simple type} if one of the following holds:
\begin{enumerate}
\def\labelenumi{{\rm(\alph{enumi})}}
\def\theenumi{\alph{enumi}}
\item $J=J(\beta,\curlya)$ is a subgroup associated to a simple
stratum~$[\curlya,n,0,\beta]$, such that if we write~$E=F[\beta]$
and~$B=\End_E(V)$ then~$\urlyb=\curlya\cap B$ a maximal~$\oE$-order
in~$B$. Moreover, there exists a simple character~$\theta\in\mathcal
C(\curlya,0,\beta, \psi_F)$ such that 
$$
\lambda \cong \kappa\otimes\sigma,
$$ 
where~$\kappa$ is a~$\beta$-extension of the unique irreducible
representation~$\eta$ of~$J^1=J^1(\beta,\curlya)$, which
contains~$\theta$, and~$\sigma$ is the inflation to~$J$ of a cuspidal
representation of~$J/J^1\cong \UU(\urlyb)/\UU^1(\urlyb)\cong
\GL_r(\kE)$. 
\item $(J,\lambda)=(\UU(\curlya),\sigma)$, where~$\curlya$ is a
maximal hereditary~$\oF$-order in~$A$ and~$\sigma$ is an inflation
of a cuspidal representation of~$\UU(\curlya)/\UU^1(\curlya)\cong
\GL_N(\kF)$. 
\end{enumerate}
\end{defi}

In case~(a),~$\mathbf J =E^{\times}J$, and in case~(b),~$\mathbf
J=F^{\times}\UU(\curlya)$. In practice we will treat~(b) as a special
case of~(a), 
with~$\beta=0$,~$E=F$,~$\urlyb=\curlya$,~$J^1=H^1=\UU^1(\curlya)$
and~$\theta$,~$\eta$,~$\kappa$ all trivial. We will refer to~(b) as
the level zero case. 

\medskip

We are going to describe the main result of this
section. Let~$\FF=\{V_i: 1\le i \le N\}$  be \textbf{any}
maximal~$F$-flag in~$V$, let~$U$ be the unipotent radical of
the~$G$-stabiliser of~$\FF$, and let~$\psi_{\alpha}: U\rightarrow
\CC^{\times}$ be \textbf{any} non-degenerate character
of~$U$. Let~$\pi$  be a supercuspidal representation of~$G$. We know
by~\cite{bh}~Proposition~1.6 that there exists~$(\mathbf J, \Lambda)$
as above, such that 
$$ 
\pi\cong \cIndu{\mathbf J}{G}{\Lambda}, \quad \Hom_{U\cap \mathbf
  J}(\psi_{\alpha}, \Lambda)\neq 0.
$$
Moreover we know by Proposition~\ref{conju} that the above properties
determine such a pair~$(\mathbf J, \Lambda)$ up to conjugation
by~$u\in U$. Since~$\Lambda|_{H^1}=(\dim \Lambda) \theta$, we obtain
that  
$$ 
\theta(u)=\psi_{\alpha}(u), \quad \forall u\in U\cap H^1.
$$
Since $J$ normalises $H^1$ and intertwines $\theta$, we may define:

\begin{defi}\label{Psialpha} Let $\Psi_{\alpha}: (J\cap U)H^1 \rightarrow \CC^{\times}$ be the character given by
$$
\Psi_\alpha(uh) = \psi_\alpha(u)\theta(h), \quad\forall u\in J\cap U,\quad \forall h\in H^1.
$$  
\end{defi}

We also define the following subgroups:
\begin{defi}\label{rimdef} Set
$$ \mir_F= \{g\in G: (g-1) V\subseteq V_{N-1}\}, \quad \rim = (\mir_F\cap \UU(\curlya)) \UU^1(\curlya),$$
and 
$$ \GG_F=\{g\in G: g v_1=v_1,\ \forall v_1\in V_1 \}, \quad \GG_{\curlya}= (\GG_F\cap \UU(\curlya)) \UU^1(\curlya).$$
\end{defi}

Let~$\kA=\{g\in G: g^{-1}\curlya g = \curlya\}$ be the~$G$-normaliser 
of~$\curlya$ and let~$\rho$ be the representation of~$\kA$ given by
$$
\rho=\Indu{\mathbf J}{\kA}{\Lambda}.
$$ 

The main result of this Section is the following Theorem.

\begin{thm}\label{mainsup}  
\begin{enumerate} 
\item The restriction~$\Lambda|_{\rim\cap J}$ is an irreducible
representation of~$\rim\cap J$. Moreover,  
$$
\Lambda|_{\rim\cap J} \cong \Indu{(J\cap U) H^1}{\rim\cap J}
{\Psi_{\alpha}}.
$$
\item The restriction~$\rho|_{\rim}$ is an irreducible
representation of~$\rim$. Moreover,  
$$
\rho|_{\rim} \cong \Indu{(J\cap U) H^1}{\rim}{\Psi_{\alpha}}.
$$
\end{enumerate}
Further, both~{\rm (i)} and~{\rm (ii)} hold if we
replace~$\mir_{\curlya}$ with~$\GG_{\curlya}$.
\end{thm}

The strategy is to show that Theorem~\ref{mainsup} holds for a
particular choice of~$U$ and~$\psi_{\alpha}$, constructed from
Theorem~\ref{wavingit}, and then show that the general result may be
obtained by conjugating by some~$g\in \mathbf J$. Before proceeding
with the proof we note that~$\Psi_{\alpha}$ occurs in~$\Lambda$ with
multiplicity~$1$, since~$\psi_{\alpha}$ occurs in~$\Lambda$ with
multiplicity~$1$ and the proof of~\cite{bh}~Lemma~3.1 implies: 

\begin{cor}\label{alphatwo} The character~$\Psi_{\alpha}$ occurs
in~$\pi$ with multiplicity~$1$. 
\end{cor}

We note that the level~$0$ case can be formally recovered from the
general case with~$\beta=0$ and~$\theta$ the trivial character, and is
a well known result of Gel$'$fand~\cite{gelfand}.

%%%%%%%%%%%%%%%%%%%%%%%%%% Some decompositions  %%%%%%%%%%%%%%%%%%%%%%

\subsection{Some decompositions}\label{deecthm}

We will need some decompositions of~$G$, and also of other general linear
groups, so we state the following Theorem for a general field~$\rF$.

\begin{thm}\label{decfield} Let~$\rF$ be any field, let~$\rV$ be
an~$N$-dimensional~$\rF$-vector space, set~$\rG=\Aut_{\rF}(\rV)$. 
Let~$\rK$ be an extension of~$\rF$ of degree~$N$ and suppose that we
are given an embedding of algebras $\iota:\rK\hookrightarrow
\End_{\rF}(\rV)$. Then the following hold: 
\begin{enumerate}
\item Let~$v\in\rV$,~$v\neq 0$, and put~$\GG_v=\{g\in\rG: gv
  =v\}$; then every~$g\in \rG$ can be uniquely decomposed as  
$$ 
g= x h, \quad x\in \iota(\rK^{\times}), \quad h\in \GG_v;
$$ 
\item Let~$\rV'\subset \rV$ be an~$\rF$-subspace of~$\rV$ of
dimension~$N-1$ and set~$\mir= \{g\in \rG: (g-1)\rV\subseteq \rV'\}$; 
then every~$g\in \rG$ can be uniquely decomposed as
$$ 
g=xh, \quad x \in \iota(\rK^{\times}), \quad h\in \mir.
$$
\end{enumerate} 
\end{thm} 

\begin{proof} We can view~$\rV$ as a~$\rK$-vector space,
via~$\iota$. Since,~$[\rK :\rF]= \dim_{\rF} \rV$, we
obtain~$\dim_{\rK} \rV =1$. Hence~$\rK^{\times}$ acts
transitively on the set of non-zero vectors in~$\rV$ and
that~$\rG=\iota(\rK^{\times}) \GG_v$. If~$x\in \GG_v$ then it has
eigenvalue~$1$, and so~$\iota(\rK^{\times})\cap \GG_v =\{1\}$. This
establishes Part~(i).

Choose a basis~$\BB=\{v_1, \ldots, v_N\}$ of~$\rV$ such that~$\rV'=
\rF v_1+\cdots+ \rF v_{N-1}$. Let~$\delta: \rG\rightarrow \rG$ be the
map~$g\mapsto \rw (g^{\top})^{-1} \rw$, where~$\rw\in\rG$ is defined
on~$\BB$ by~$\rw(v_i)= v_{N-i+1}$, and~$g^{\top}$ denotes the
transpose of~$g$ with respect to the basis~$\BB$. We have~$\delta
^2=\id$ and~$\delta(\mir)= \GG_{v_1}$. Part~(ii) follows
from Part~(i) with~$\iota$ replaced with~$\delta \circ \iota$
and~$v=v_1$.
\end{proof}

We apply our decomposition theorem to prove several results on the
intersection with the groups~$\mir_F$ and~$\GG_F$. 
Analogously, we set
$$
\mir_E=\{g\in B^{\times}: (g-1) V \subseteq V_{N-d} \}, \quad
\mir_{\urlyb}=(\mir_E \cap\UU(\urlyb))\UU^1(\urlyb)
$$
and
$$
\GG_E=\{g\in B^{\times}: g w_1 =w_1\}, \quad \GG_{\urlyb}=(\GG_E\cap
\UU(\urlyb)) \UU^1(\urlyb).
$$
Let~$K$ be a maximal unramified extension of~$E$, which
normalises~$\curlya$, so that~$[K:E]=N/d$. Theorem~\ref{decfield}
implies that~$B^{\times}= K^{\times} \mir_E$. Since~$V_{N-d}\subseteq
V_{N-1}$, we obtain that
$$ 
\mir_F\cap B^{\times} =\mir_E.
$$
 
\begin{cor}\label{decem} The following decompositions hold: 
\begin{enumerate}
\item $\kA=K^{\times} (\UU(\curlya)\cap \mir_F), \quad \UU^1(\curlya)=
  (\UU^1(\curlya)\cap \mir_F) (1+\pK)$; 
\item $\mathbf J = K^{\times} (J\cap \mir_F),\quad J^1= (J^1 \cap
  \mir_F) (1+\pK)$; 
\item $\kB= K^{\times} (\UU(\urlyb)\cap \mir_E), \quad \UU^1(\urlyb)=
  (\UU^1(\urlyb)\cap \mir_E) (1+\pK)$. 
\end{enumerate}
Moreover, the statements remain true if we replace~$\mir_F$
with~$\GG_F$ and~$\mir_E$ with~$\GG_E$. 
\end{cor}

\begin{proof} Note, that~$K^{\times}\subset \kB \subset \mathbf
J\subset \kA$. It is enough to prove~(i), since~(ii) and~(iii) are
obtained by intersecting with~$\mathbf J$ and~$\kB$,
respectively. According to Theorem~\ref{decfield} we may 
write~$\kA=K^{\times}(\kA\cap \mir_F)$. 

Let~$\LL=\{L_i:i\in\ZZ\}$ be the~$\oE$-lattice
chain in~$V$ associated to~$\curlya$. 
According to~\cite{weil} Theorem~II-1 there exists a decomposition of~$V=
\sum_{i=1}^N V^i$ into one dimensional subspaces such that~$V_{N-1}=
\sum_{i=1}^{N-1} V^i$ and~$L_j= \sum_{i=1}^N L_j \cap V^i$, for
all~$j\in \ZZ$. If~$g\in \kA\cap \mir_F$ then by projecting to
the~$V^N$ subspace we obtain~$L_{j+v_{\curlya}(g)}\cap V^N = L_j\cap
V^N$, for all~$j\in \ZZ$. This implies that~$v_{\curlya}(g)=0$,
hence~$g\in \UU(\curlya)\cap \mir_F$. Hence~$\kA\cap \mir_F=
\UU(\curlya)\cap \mir_F$.

Let~$g\in \UU^1(\curlya)$. We may write by above~$g= h x$,
where~$h\in \UU(\curlya) \cap \mir_F$ and~$x\in\oK^{\times}$. 
Let~$\bar{g}$,~$\bar{h}$ and~$\bar{x}$ be the images of~$g$,~$h$,~$x$
in~$\Aut_{\kF}(L_m/L_{m+1})$, where~$L_m\in\LL$, such
that~$L_{m+1}+L_m\cap V_{N-1}\neq L_m$. Since~$g\in \UU^1(\curlya)$, 
we have~$\bar{g}=1$. Now   
$$
(\bar{h}-1) (L_m/L_{m+1}) \subseteq  (L_m\cap
V_{N-1}+L_{m+1})/L_{m+1}.
$$
Our assumption on~$L_m$, implies that~$\bar{h}$ has
eigenvalue~$1$. Since~$\bar{h}\bar{x}=1$,~$\bar{x}$ also has
eigenvalue~$1$, hence~$\bar {x}=1$, which implies that~$x\in 1+\pK$,
and so~$h\in \UU^1(\curlya)\cap \mir_F$.  
\end{proof} 
 
\begin{cor}\label{rimdec} We have 
$$ 
\kA= \rim \mathbf J, \quad \rim\cap \mathbf J = (\mir_F \cap \mathbf
J) J^1= \mir_{\urlyb} J^1.
$$
Moreover, the statement is true if we replace~$\rim$
with~$\GG_{\curlya}$;~$\mir_F$ with~$\GG_F$ and~$\mir_{\BB}$
with~$\GG_{\BB}$.
\end{cor}

We end this section with an observation which will prove useful later.

\begin{lem}\label{detv}  Suppose that~$g\in \kA$ then 
$$
v_F(\det g)= \frac{N v_{\curlya}(g)}{e(\curlya|\oF)}.
$$
In particular,~$\mathbf J \cap U= J\cap U.$
\end{lem}

\begin{proof} Let~$e= e(\curlya |\oF)$ and set~$h= g^e
\pif^{-v_{\curlya}(g)}$, then~$h\in \kA$
and~$v_{\curlya}(h)=0$. Hence,~$h\in\UU(\curlya)$. 
Since~$\UU(\curlya)$ is compact, this implies that~$\det h \in
\oF^{\times}$ and this yields the Lemma.    
\end{proof}

%%%%%%%% The proof of Theorem~\ref{mainsup} in a special case %%%%%%%%

\subsection{The proof of Theorem~\ref{mainsup} in a special 
case}\label{specase} 

Let us fix a maximal simple type~$(J,\lambda)$.
Let~$[\curlya,n,0,\beta]$ be an associated simple stratum, and
let~$\theta\in \mathcal C(\curlya, 0, \beta)$ be such
that~$\lambda|_{H^1}=(\dim \lambda) \theta$. Let~$(\FF, \chi, b)$ be a
triple given by Theorem~\ref{wavingit} and let~$U$ be the unipotent radical of
the~$G$-stabiliser of~$\FF$. By construction, a
subflag~$\FF_0:=\{V_{di}: 1\le i\le r\}$ of~$\FF$  is a maximal
stable~$E$-flag in~$V$, and let~$\BB_E=\{w_1,\ldots w_r\}$ denote
an~$E$-basis of~$V$ chosen as in Corollary~\ref{goodbasis}. We choose
some~$a\in X_{\FF}$ such that~$\chi=\psi_a$, and we set  
$$ 
\alpha= a+b.
$$
Then  we have 
$$ 
\psi_{\alpha}(u)= \psi_a(u) \psi_b(u)= \psi_a(u) =\theta(u), \quad
\forall u\in U\cap H^1 .
$$
The first equality is trivial; the second holds, since~$b\in \curlya$
and~$u\in \UU^1(\curlya)$; the third is Theorem~\ref{wavingit}(ii). 
Hence we may define~$\Psi_{\alpha}: (J\cap U)H^1 \rightarrow 
\CC^{\times}$ by~$\Psi_\alpha(uh) = \psi_\alpha(u)\theta(h)$, 
for~$u\in J\cap U$ and~$h\in H^1$, as above.

\begin{lem}\label{welldfn} The function~$\psi_b$ defines a linear
character on~$(J\cap U)J^1$ and 
$$
\Psi_{\alpha}(j)= \Theta(j) \psi_b(j), \quad 
\forall j \in (J\cap U)H^1
$$
where the character~$\Theta$ is given
by~$\Theta(uh)=\psi_a(u)\theta(h)$, for~$u\in (J\cap U)$ and~$h\in
H$, as in Theorem~\ref{Psikappa}. 
\end{lem}

\begin{proof} Since~$b\in X_{\FF}$, by Lemma~\ref{charX},~$\psi_b$
defines a linear character of~$U$. Since~$b\in \curlya$, we
have~$\psi_b(u)=1$, for all~$u\in \UU^1(\curlya)$. This implies that 
$$ 
\psi_b(uj)= \psi_b(u), \quad \forall u\in U\cap J, \quad j\in J^1.
$$
Now~$J$ normalises~$J^1$ and hence~$\psi_b$ is a character on~$(J\cap
U)J^1$. Since~$\alpha=a+ b$, we have an equality of
functions~$\psi_{\alpha}=\psi_{a} \psi_b$ and hence for every~$u\in
J\cap U$ and~$h\in H^1$ we have
$$
\Psi_{\alpha}(uh)=\psi_b(u) \psi_{a}(u) \theta(h)= \psi_b(uh)
\Theta(uh).
$$ 
\end{proof}

%\begin{defi}\label{mir} 
Let~$\LL=\{L_i:i\in\ZZ\}$ be the~$\oE$-lattice
chain in~$V$ associated to~$\curlya$, and put
$$
\mir_{\kE}=\left\{ b\in \UU(\urlyb)/\UU^1(\urlyb): (b-1)
(L_0/L_1)\subseteq \sum_{j=1}^{r-1} \kE (w_j+L_1)\right\},
$$
the image of~$\mir_{\urlyb}$ in~$\UU(\urlyb)/\UU^1(\urlyb)$. 
Similarly, let 
$$ 
\GG_{\kE}=\{ b\in \UU(\urlyb)/\UU^1(\urlyb): b (w_1 +L_1)= w_1 + L_1
\}
$$
be the image of~$\GG_{\urlyb}$ in~$\UU(\urlyb)/\UU^1(\urlyb)$.
%\end{defi}

The group~$\mir_{\kE}$ is known as a mirabolic subgroup. We note
that it contains the image of~$\UU(\urlyb)\cap U$
in~$\UU(\urlyb)/\UU^1(\urlyb)$. Moreover, with respect to the
basis~$w_1+L_1,...,w_r+L_1$ of~$L_0/L_1$, the group~$\mir_{\kE}$ is
identified with the subgroup of matrices of the form 
$$
\begin{pmatrix} *&\cdots&*&*\\ \vdots&\ddots&\vdots&\vdots\\ 
*&\cdots&*&*\\ 0&\cdots&0&1 \end{pmatrix}.
$$
and  the image of~$\UU(\urlyb)\cap U$ in~$\UU(\urlyb)/\UU^1(\urlyb)$
is identified with the subgroup of unipotent upper-triangular
matrices.

\begin{lem}\label{mirabolic} Let~$\sigma$ be a cuspidal representation
of~$\UU(\urlyb)/\UU^1(\urlyb)$. Let~$U_{\kE}$ be the image
of~$\UU(\urlyb)\cap U$ in~$\UU(\urlyb)/\UU^1(\urlyb)$, so that 
$$  
U_{\kE} \cong (\UU(\urlyb)\cap U)/ (\UU^1(\urlyb)\cap U),
$$
and let~$\psi$ be any non-degenerate character of~$U_{\kE}$. Then
$$
\sigma|_{\mir_{\kE}}\cong \Indu{U_{\kE}}{\mir_{\kE}}{\psi}, \quad
\sigma|_{\GG_{\kE}}\cong \Indu{U_{\kE}}{\GG_{\kE}}{\psi}.
$$
Moreover, the representations~$\sigma|_{\mir_{\kE}}$
and~$\sigma|_{\GG_{\kE}}$ are irreducible representations
of~$\mir_{\kE}$ and~$\GG_{\kE}$, respectively.
\end{lem}

\begin{proof} The statement for~$\mir_{\kE}$
is~\cite{gelfand}~Theorem~8. Set~$\bar{w}_i= w_i + L_1$, for~$1\le i
\le r$, and~$\BB_{\kE}=\{\bar{w}_1,\ldots, \bar{w}_r\}$,
a~$\kE$-basis of~$L_0/L_1$. We identify~$\UU(\urlyb)/\UU^1(\urlyb)$
with~$\GL_r(\kE)$, via~$\BB_{\kE}$. Let~$\delta:
\GL_r(\kE)\rightarrow \GL_r(\kE)$ be the automorphism given
by~$\delta(g)=\rw (g^\top)^{-1}\rw$, where~$g^{\top}$ denotes the
transpose of~$g$ and~$\rw$ is given by $\rw(\bar{w}_i)= \bar{w}_{r-i+1}$, for $1\le i\le r$.
Then~$\delta(U_{\kE})=U_{\kE}$, and~$\psi^{\delta}$, given
by~$\psi^{\delta}(u)= \psi(\delta(u))$, for all~$u\in U_{\kE}$, is a
non-degenerate character. The  representation of~$\GL_r(\kE)$, given
by~$\sigma^{\delta}(g)=\sigma(\delta(g))$ is cuspidal. Hence, by the
statement for~$\mir_{\kE}$,
$$
\sigma^{\delta}|_{\mir_{\kE}}\cong
\Indu{U_{\kE}}{\mir_{\kE}}{\psi^{\delta}}.
$$ 
Since~$\delta(\GG_{\kE})=\mir_{\kE}$ and~$\delta^2=\id$, by twisting
by~$\delta$ we obtain 
$$ 
\sigma|_{\GG_{\kE}}\cong \Indu{U_{\kE}}{\GG_{\kE}}{\psi}.
$$
The irreducibility follows from the irreducibility
of~$\sigma^{\delta}|_{\mir_{\kE}}$.
\end{proof}  

\begin{thm}\label{Psilambda} The restriction~$\lambda|_{\mir_{\urlyb}
J^1}$ is an irreducible representation of~$\mir_{\urlyb} J^1$ and  
$$
\lambda|_{\mir_{\urlyb} J^1} \cong \Indu{(J\cap U)H^1}{\mir_{\urlyb}
J^1} \Psi_\alpha.
$$
The statement remains true if we replace~$\mir_{\urlyb}$
with~$\GG_{\urlyb}$.
\end{thm}

\begin{proof} Set~$\mir=\mir_{\urlyb}$. In the level~$0$
case,~$\mir=\mir_{\curlya}$ and the Theorem asserts
that~$\sigma|_{\mir_{\curlya}}$ is irreducible and 
$$
\sigma|_{\mir_{\curlya}}\cong \Indu{(U\cap
  \UU(\curlya))\UU^1(\curlya)}{\mir_{\curlya}}{\psi_b}.
$$
Since~$\sigma$ is an inflation of a cuspidal representation and~$\FF$
and~$b$ were chosen in Corollary~\ref{version0} so
that~$\overline{\psi}_b$ is non-degenerate, the assertion is
Lemma~\ref{mirabolic}.

Let us consider the general case. Since~$\sigma$ is an inflation of a
cuspidal representation of~$\UU(\urlyb)/\UU^1(\urlyb)$ and~$\FF$
and~$b$ were chosen in Theorem~\ref{wavingit}, so that the
character~$\overline{\psi}_b$ is non-degenerate, we may again apply
Lemma~\ref{mirabolic} to obtain 
$$
\lambda|_{\mir J^1}\cong \kappa|_{\mir J^1}\otimes \sigma|_{\mir J^1}
\cong  \kappa|_{\mir J^1}\otimes \Indu{(J\cap U)J^1}{\mir J^1}{\psi_b}
\cong\Indu{(J\cap U)J^1}{\mir J^1}{\kappa\otimes\psi_b}.
$$  
Theorem~\ref{Psikappa} implies that~$\kappa|_{(J\cap U) J^1}\cong
\Indu{(J\cap U)H^1}{(J\cap U)J^1}{\Theta}$  and hence 
$$
\lambda|_{\mir J^1}\cong \Indu{(J\cap U)J^1}{\mir J^1}{\Indu{(J\cap U)H^1}
{(J\cap U)J^1}{\Theta\otimes \psi_b}}\cong\Indu{(J\cap U)H^1}{\mir J^1}
{\Psi_{\alpha}}
$$   
where the last isomorphism is given by Lemma~\ref{welldfn} and the
transitivity of induction. Moreover, since by Lemma~\ref{mirabolic}
the restriction~$\sigma|_{\mir J^1}$ is irreducible, a straightforward
modification of~\cite{bk}~(5.3.2) implies that~$\lambda|_{\mir J^1}$
is an irreducible representation of~$\mir J^1$. The proof
for~$\GG_{\urlyb}$ is analogous.
\end{proof}

Since, by Corollary~\ref{rimdec}, we have $\rim\cap \mathbf J =
\mir_{\urlyb} J^1$, we have now proved Theorem~\ref{mainsup}(i) in our
special case. We also record the following: Since~$\mathbf J \cap U=
J\cap U$, we immediately get

\begin{cor} $\Hom_{U\cap \mathbf J}(\psi_{\alpha}, \Lambda)\neq 0$.
\end{cor}

Finally, part~(ii) of Theorem~\ref{mainsup} is given by:

\begin{prop}\label{Psirho} The restrictions of~$\rho$ to~$\rim$  and
to~$\GG_{\curlya}$ are irreducible representations of~$\rim$
and~$\GG_{\curlya}$, respectively. Moreover,   
$$
\rho|_{\rim}\cong \Indu{(J\cap U)H^1}{\rim}{\Psi_{\alpha}}, \quad
\rho|_{\GG_{\curlya}}\cong \Indu{(J\cap U)H^1}{\GG_{\curlya}}
{\Psi_{\alpha}}.
$$
\end{prop}

\begin{proof} We will prove statement for~$\mir_{\curlya}$, the proof
for~$\GG_{\curlya}$ is analogous. We have
$$
\rho|_{\rim}\cong 
\Indu{\rim \cap \mathbf J}{\rim}{\Lambda|_{\rim\cap \mathbf J}}
\cong \Indu{\mir_{\urlyb} J^1}{\rim}{\lambda|_{\mir_{\urlyb} J^1}}
\cong \Indu{(J\cap U) H^1}{\rim}{\Psi_{\alpha}}
$$  
where the last two isomorphisms follow from Corollary~\ref{rimdec} and
Theorem~\ref{Psilambda}. Since~$\Psi_{\alpha}|_{H^1}=\theta$ we have 
$$
I_G(\Psi_{\alpha})\subseteq I_G(\theta)=J^1 B^{\times} J^1
$$ 
by~\cite{bk}~(3.3.2). Since~$\lambda|_{\mir_{\urlyb} J^1}\cong
\Indu{(J\cap U)H^1}{\mir_{\urlyb} J^1}
     {\Psi_{\alpha}}$,~\cite{bk}~(4.1.1) and~(4.1.5) imply that 
$$
I_G(\lambda|_{\mir_{\urlyb} J^1})=\mir_{\urlyb} J^1 I_G(\Psi_{\alpha}) \mir_{\urlyb} J^1
\subseteq J^1 B^{\times} J^1.
$$  
Hence, 
$$
I_{\rim}(\lambda|_{\mir_{\urlyb} J^1})\subseteq 
\rim \cap (J^1 B^{\times} J^1)=
(\UU(\urlyb)\cap \rim)J^1= \mir_{\urlyb} J^1
$$ 
and hence~$\rho|_{\rim}$ is irreducible. We note that~$\rim$
contains~$\UU^1(\curlya)$, and hence~$J^1$;~$\UU(\curlya)\cap
B^{\times}=\UU(\urlyb)$ and the last equality  above follows from
Corollary~\ref{decem}.
\end{proof}
  
This completes the proof of
Theorem~\ref{mainsup} for our special choice of~$\FF$ and
$\psi_{\alpha}$.

%%%%%%% The proof of Theorem~\ref{mainsup} in the general case %%%%%%%

\subsection{The proof of Theorem~\ref{mainsup} in the general case}
We will prove Theorem \ref{mainsup} in the general case, by showing that 
after conjugation by some $g\in\mathbf J$ we end up in the special case, considered
above.
\begin{proof}
Let~$\FF'=\{V'_i:1\le i\le N\}$ be any maximal~$F$-flag in~$V$, and
let~$U'$ be the unipotent radical of the~$G$-stabiliser of~$\FF'$, and
let $\psi_{\alpha'}$ be any smooth non-degenerate character
of~$U'$. Let~$\pi$ be a supercuspidal representation of~$G$, then
there exists a pair~$(\mathbf J,\Lambda)$, such
that~$\pi\cong\cIndu{\mathbf J}{G}{\Lambda}$ and~$\Hom_{U'\cap \mathbf
  J}(\psi_{\alpha'},\Lambda)\neq 0$, and~$\Lambda|_J\cong \lambda$,
where~$J$ is the maximal compact open subgroup of~$\mathbf J$,
and~$(J,\lambda)$ is a maximal simple type, with the
stratum~$[\curlya,n,0,\beta].$ Moreover,~$\lambda|_{H^1}= (\dim
\lambda) \theta$, where~$\theta\in \mathcal C(\curlya,0,\beta)$. We
define~$\rim'$ as in Definition~\ref{rimdef} 
and~$\Psi_{\alpha'}$ as in Definition~\ref{Psialpha}. 
Let~$\FF=\{V_i:1\le i\le
N\}$,~$U$,~$\psi_{\alpha}$ and~$\rim$ be as in~\S\ref{specase}. 
By Corollary~\ref{goodbasis} the
character~$\psi_{\alpha}$ is non-degenerate, and we know that 
$$
\Hom_{U'\cap \mathbf J}(\psi_{\alpha'}, \Lambda)\neq 0, \quad
\textrm{and} \quad \Hom_{U\cap \mathbf J}(\psi_{\alpha},\Lambda)\neq
0.
$$    
Hence by~\cite{bh}~Proposition~1.6~(ii), there exists~$g\in \mathbf
J$, such that~$U'=U^g$ and~$\psi_{\alpha'}=\psi_{\alpha}^g$. In
particular,~$\FF'=g\FF$, and hence $V'_{N-1}= g V_{N-1}$, which
implies  that~$\rim'= \rim^g$. Since~$g\in \mathbf J$, we
have~$J=J^g$,~$H^1=(H^1)^g$,~$\theta=\theta^g$. Hence~$(J\cap U') H^1
= \bigl((J\cap U)H^1\bigr)^g$ and~$\Psi_{\alpha'}=\Psi_{\alpha}^g$. We
have proved the Theorem for~$\FF$ and~$\psi_{\alpha}$, now twisting
by~$g$, we obtain the result for~$\FF'$  and~$\psi_{\alpha'}$.
\end{proof}

\begin{remar}\label{natural} It follows from the proof that any~$\FF'$
and~$\psi_{\alpha'}$, with the property that~$\Hom_{U'\cap \mathbf
 J}(\psi_{\alpha'},\Lambda)\neq 0$, arise from the construction in
Theorem~\ref{wavingit}, once we replace~$\beta$ by~$g \beta g^{-1}$,
for some~$g\in \mathbf J$, and so the construction in 
Theorem~\ref{wavingit} is a natural one. 
\end{remar}

%%%%%%% A characterisation of~$\Indu{\mathbf J}{\kA}{\Lambda}$ %%%%%%%

\subsection{A characterisation of~$\Indu{\mathbf
 J}{\kA}{\Lambda}$}\label{charind} 

We observe that a result of Gel$'$fand characterising cuspidal
representations of~$\GL_N(\Fq)$ implies a very similar result, for the
representations of~$\kA$ of the form~$\Indu{\mathbf J}{\kA}{\Lambda}$.

\begin{prop}\label{char} Let~$\tau$ be a representation of~$\kA$ such 
that 
$$ 
\tau|_{\rim}\cong \Indu{(J\cap U) H^1}{\rim} \Psi_{\alpha},
$$
then 
$$ 
\tau\cong \Indu{\mathbf J}{\kA}{\Lambda},
$$
for some representation~$\Lambda$ of~$\mathbf J$, such
that~$(J,\Lambda|_J)$ is a maximal simple type, as in 
Definition~\ref{maxtype}.
\end{prop} 

\begin{proof} Since~$ \Indu{(J\cap U) H^1}{\rim} \Psi_{\alpha}$ is
irreducible, so is~$\tau$. Now, 
$$ 
\Indu{(J\cap U)H^1}{J}{\Psi_{\alpha}}\cong \kappa \otimes\Indu{(J\cap
  U) J^1}{J}{\psi_b}\cong \prod_{\sigma}\kappa \otimes \sigma,
$$
where the product runs over all the generic representations~$\sigma$
of~$J/J^1\cong \UU(\urlyb)/\UU^1(\urlyb)\cong
\GL_r(\kE)$. Hence~$\tau|_{J}$ will contain a summand of the
form~$\kappa\otimes \sigma$. It follows from~\cite{bk}~(5.3.2) that
the~$I_{\UU(\curlya)}(\kappa\otimes \sigma)\subseteq  (J B^{\times} J)
\cap \UU(\curlya)= J$. Hence~$\Indu{J}{\UU(\curlya)} {\kappa\otimes
  \sigma}$ is irreducible, and so is isomorphic to~$\tau$. Restricting
to~$\rim$, we obtain that~$\sigma|_{\mir_{\urlyb}J^1}$ is irreducible,
and~\cite{gelfand} implies that~$\sigma$ is cuspidal. Hence,~$\tau|_J$
will contain some simple type~$\lambda$, it follows
from~\cite{bk}~\S6.2 that~$\tau|_{\mathbf J}$ will contain
some~$\Lambda$, and hence~$\tau\cong \Indu{\mathbf J}{\kA}{\Lambda}$.
\end{proof}

\section{Realisation of maximal simple types}\label{realization}

We continue with the situation described in the beginning
of~\S\ref{super}. Let~$U$ be a maximal unipotent subgroup of~$G$,
and~$\psi_{\alpha}$ non-degenerate character of~$U$. Let~$\pi\cong
\cIndu{\mathbf J}{G}{\Lambda}$ be a supercuspidal representation,
and~$\Hom_{U\cap \mathbf J}(\psi_{\alpha}, \Lambda)\neq
0$. Theorem~\ref{mainsup} allows us to use a rather general result of
Alperin and James~\cite{alperin}, and realize the
representation~$\Lambda$ as a concrete space, and describe the action
of~$\mathbf J$ on this space in terms of the character of~$\Lambda$
and~$\Psi_{\alpha}$. This concrete realization enables us to compute a
certain matrix coefficient of~$\pi$, and by integrating it we obtain
an explicit Whittaker function for~$\pi$.

%%%%%%%%%%%%%%%%%%%%%%%%%% Bessel functions %%%%%%%%%%%%%%%%%%%%%%%%%%

\subsection{Bessel functions}\label{besfun}

We will adapt the result of Alperin and James~\cite{alperin} to our
setting. Let~$\mathcal K$ be an open,
compact-modulo-centre subgroup of~$G$ and let~$\tau$ be an irreducible
smooth representation of~$\mathcal K$. 

\begin{ass}\label{assump} 
Suppose that there exists compact open  subgroups~$\mathcal U\subseteq
\mathcal M\subseteq \mathcal K$, and a linear character~$\Psi$
of~$\mathcal U$, such that the following hold:
\begin{enumerate}
\item $\tau|_{\mathcal M}$ is irreducible representation of~$\mathcal M$;
\item $\tau|_{\mathcal M}\cong \Indu{\mathcal U}{\mathcal M} {\Psi}$.
\end{enumerate} 
\end{ass} 

Let~$\mathcal N$ be an open, normal subgroup of~$\mathcal K$ contained
in the~$\Ker \tau$.  Set 
$$
e_{\Psi}= (\mathcal U :\mathcal N)^{-1} \sum_{h\in \mathcal U/\mathcal
  N} \Psi(h) h^{-1}.
$$
Let~$\chi=\chi_{\tau}$ be the (trace) character of~$\tau$ and
let~$\omega=\omega_{\tau}$ be the central character of~$\tau$, so that
$$ 
\chi(xg)= \omega(x) \chi(g), \quad \forall x\in F^{\times}, \quad
\forall g\in \mathcal K.
$$

\begin{defi}\label{defbes} The \emph{Bessel function}~$\bes:\mathcal K
  \rightarrow \mathbb C$ of~$\tau$ is defined by:
$$
\bes(g)= \tr_{\tau}(e_{\Psi} g)=(\mathcal U: \mathcal
N)^{-1}\sum_{h\in \mathcal U/\mathcal N} \Psi(h^{-1})\chi(gh).
$$
\end{defi} 

%some interesting junk:
%
%\begin{lem} Let~$\mathcal H$ be a subgroup of a finite
%  group~$\mathcal G$ and let~$\tau$ be an irreducible representation
%  of~$\mathcal H$ affording character~$\chi_{\tau}$. Let~$g\in
%  \mathcal G$ and let  
%$$e_{\tau}=\dim \tau /|\mathcal H| \sum_{h\in \mathcal H}
%\chi_{\tau}(h^{-1})h.$$
%Suppose that in the group algebra of~$\mathcal G$ we have~$e_{\tau}g
%e_{\tau}=0$ then~$g$ doesnot intertwine~$\tau$. 
%\end{lem}
%\begin{proof}We may write in a unique way
%$$e_{\tau}g e_{\tau}=\sum_{k\in \mathcal G} \mu_k k$$
%where~$\mu_k \in \mathbb C$. Now  
%$$\mu_g=(\dim \tau/|\mathcal H |)^2 \sum_{h\in \mathcal H\cap g
%  \mathcal H g^{-1}} \chi_{\tau}(h^{-1})\chi_{\tau}(g^{-1} h g ).$$ 
%If~$e_{\tau} g e_{\tau}=0$ then~$\mu_g=0$ and hence~$\Hom_{\mathcal H\cap 
 % \mathcal H^ g}(\tau, \tau^g)=0.$
%\end{proof}
 
\begin{prop}\label{proper} The Bessel function~$\bes$ has the
following properties: 
\begin{enumerate}
\item $\bes(1)=1$; 
\item $\bes(xg)=\bes(gx)=\omega(x)\bes(g), \quad \forall x\in
  F^{\times}, \quad \forall g \in \mathcal  K$; 
\item $\bes(hg)=\bes(gh)=\Psi(h)\bes(g), \quad \forall h\in \mathcal
  U, \quad \forall g\in \mathcal  K$;
%\item an element~$g\in \mathbf J$ intertwines~$\Psi_{\alpha}$ if and
% only if~$\bes(g)\neq 0$;
\item if~$\bes(g)\neq 0$ then~$g$ intertwines~$\Psi$; in
particular, if~$m\in \mathcal M$ then~$\bes(m)\neq 0$ if and only
if~$m\in \mathcal U$;
\item for all~$g_1$,~$g_2\in \mathcal K$ we have 
$$ 
\sum_{m\in \mathcal M/\mathcal U} \bes(g_1 m) \bes(m^{-1}g_2)=
\bes(g_1 g_2).
$$
\end{enumerate} 
\end{prop}
 
\begin{proof} We observe that it is enough to prove the Proposition
for a twist of~$\tau$ by an unramified character. Twisting
by~$[g\mapsto (\omega(\pif))^{-v_F(\det(g))/N}]$ we ensure
that~$\pif^{\mathbb Z}\mathcal N$ lies in the kernel
of~$\tau$. Hence,~$\Ker \tau$ is of finite index in~$\mathcal K$ and
we may consider~$\tau$ as a representation of a finite group.  

Part~(i) is a reformulation of the fact that~$\Indu{\mathcal
U}{\mathcal M}{\Psi}$ is irreducible. Since~$\chi$ is defined by
matrix trace, we have~$\chi(g g_1)=\chi(g_1g)$, for all~$g$,~$g_1\in
\mathcal K$, and Parts~(ii) and~(iii) are straightforward
consequences of the definition of~$\bes$.  

Part~(iv): Part~(iii) implies that~$\bes$ is a~$\check{\Psi}$-spherical
function on~$\mathcal  K$, in the sense of~\cite{bk} (4.1),
where~$\check{\Psi}$ is the dual of~$\Psi$. Hence if~$\bes(g)\neq
0$ then according to~\cite{bk}~(4.1.1),~$g$ intertwines~$\Psi$.

Since~$\Indu{\mathcal U}{\mathcal M}{\Psi}$ is irreducible,
the~$\mathcal M$-intertwining of~$\Psi$ is equal to~$\mathcal U$. Now
Parts~(i),(iii) and the argument above finish the proof of Part~(iv).
 
Part~(v) is~\cite{alperin}~Lemma~2, or~\cite{gelfand}~Theorem~9.
\end{proof}

\begin{thm}[cf.~\cite{alperin}]\label{bessel} Let~$\mathcal S$ be the
space of functions from~$\mathcal M$ to~$\mathbb C$ satisfying the
condition 
$$
f(hm)=\Psi(h) f(m), \quad \forall h \in \mathcal U, \quad \forall m\in
\mathcal M.
$$
For each~$g\in \mathcal K$ we define an operator~$L(g)$ 
on~$\mathcal S$ by the formula 
$$ 
[L(g)f](m)= \sum_{m_1\in \mathcal M/\mathcal U} \bes(mg m_1)
f(m_1^{-1}).
$$
Then~$L$ defines a representation of~$\mathcal K$, which is isomorphic
to~$\tau$.
\end{thm}

\begin{proof} Again, it is enough to prove the statement after
twisting by unramified character, and this way we may ensure
that~$\mathcal K/\Ker \tau$ is finite. The assertion now follows from the
main Theorem in~\cite{alperin}. The level~$0$ case 
is~\cite{gelfand}~Theorem~10. 
\end{proof}  

If~$\check{\tau}$ is the dual of~$\tau$,
then~$\check{\tau}$,~$\mathcal M$,~$\mathcal U$ and~$\check{\Psi}$
satisfy Assumption~\ref{assump}. Hence Theorem~\ref{bessel} holds
for~$\check{\tau}$, with~$\check{\Psi}$ instead of~$\Psi$ and 
with~$\check{\mathcal S}$ the space of
functions from~$\mathcal M$ to~$\mathbb C$ satisfying the condition
$$
f(hm)=\check{\Psi}(h) f(m), \quad \forall h \in \mathcal U, \quad
\forall m\in \mathcal M.
$$
Moreover, the Bessel function~$\check{\bes}=\bes_{\check{\tau}}$ satisfies 
$$ 
\check{\bes}(g)= \bes(g^{-1}), \quad \forall g\in \mathcal K.
$$
Let~$(\,,\,)$ be a non-degenerate~$\mathcal K$-invariant pairing
on~$\mathcal S \times \check{\mathcal S}$. Since,~$\tau$ is
irreducible, the pairing is determined up to a scalar
multiple. Let~$\varphi\in \mathcal S$ and~$\check{\varphi} \in
\check{\mathcal S}$ be such that~$\supp \varphi =\supp
\check{\varphi}=\mathcal U$
and~$\varphi(1)=\check{\varphi}(1)=1$. Since~$\varphi$
and~$\check{\varphi}$ span the ~$\Psi$- and~$\check {\Psi}$-isotypical
subspaces in~$\mathcal S$ and~$\check{\mathcal S}$ respectively, we
may normalise~$(\,,\,)$, so that
$$ 
(\varphi,\check{\varphi})=1.
$$
This determines the pairing uniquely.  

\begin{lem}\label{pair} We have 
$$ 
(L(g)\varphi, \check{\varphi})=\bes(g), \quad \forall g\in \mathcal K.
$$
\end{lem} 

\begin{proof} It follows from Theorem~\ref{bessel} that 
$$
L(g)\varphi= \sum_{m\in \mathcal M/\mathcal U} \bes(m^{-1}g)
L(m)\varphi.
$$ 
Since the~$\mathcal M$-intertwining of~$\Psi$ is just~$\mathcal U$,
for~$m\in \mathcal M$ we have~$(L(m)\varphi, \check{\varphi})\neq
0$ if and only if~$m\in \mathcal U$. This implies the Lemma.
\end{proof}

We will apply the preceding results in several situations but,
for now, we observe that Theorem~\ref{mainsup} implies (in the
notation of~\S\ref{super}):

\begin{thm}\label{context} Assumption~\ref{assump} (and hence
Proposition~\ref{proper} and Theorem~\ref{bessel}) holds in the following
contexts:
\begin{enumerate}  
\item $\mathcal K= \kA$,~$\tau=\rho$,~$\mathcal M= \rim$ or~$\mathcal
  M= \GG_{\curlya}$,~$\mathcal U=(J\cap U)H^1$,~$\Psi=\Psi_{\alpha}$;
\item $\mathcal K=\mathbf J$,~$\tau=\Lambda$,~$\mathcal M= \rim\cap J$
  or~$\mathcal M= \GG_{\curlya}\cap J$,~$\mathcal U=(J\cap
  U)H^1$,~$\Psi=\Psi_{\alpha}$.
\end{enumerate}
\end{thm} 
 
If~$\pi$ has level~$0$ then we recover the result of
Gel$'$fand~\cite{gelfand}.

%%%%%%%%%%%%%%%%%%%% Explicit Whittaker functions %%%%%%%%%%%%%%%%%%%%

\subsection{Explicit Whittaker functions}\label{xwf}

Now we argue along the lines of~\cite{bh}~\S3. However, not only do we get 
a uniqueness statement as in~\cite{bh}, but we also obtain explicit 
formulae in terms of the character of~$\Lambda$ and~$\Psi_{\alpha}$.    

Let~$\mathcal V$ be the underlying vector space of~$\pi$ and 
let~$(\check{\pi},\check{\mathcal{V}})$ be the smooth dual 
of~$(\pi, \mathcal V)$. We denote by~$\langle , \rangle$ the pairing 
on~$\mathcal V\times\check{\mathcal V}$ given by the evaluation. 
Let~$v_{\alpha}\in \mathcal V$ and~$\check{v}_{\alpha}\in 
\check{\mathcal{V}}$ be non-zero vectors such that 
$\pi(h)v_{\alpha}=\Psi_{\alpha}(h)v_{\alpha}$ 
and~$\check{\pi}(h)v_{\alpha}=\check{\Psi}_{\alpha}(h) 
\check{v}_{\alpha}$, for all~$h\in (J\cap U)H^1$. 
Corollary~\ref{alphatwo} implies that such vectors exist and they 
are unique up to scalar multiple. We rescale so 
that~$\langle v_{\alpha}, \check{v}_{\alpha}\rangle =1$. 

\begin{prop}\label{unique} The representation~$\pi$ admits a unique 
coefficient function~$f=f_{\alpha, U}$ with the following properties:
\begin{enumerate}
\item $f(1)=1$, and
\item $f(h_1 g h_2)= \Psi_{\alpha}(h_1 h_2) f(g), \quad \forall 
h_1, h_2\in (J\cap U)H^1, \quad \forall g\in G.$
\end{enumerate}
Moreover,~$\supp f \subseteq \mathbf J$ and 
$$
f(g)=\langle \pi(g) v_{\alpha}, \check{v}_{\alpha}\rangle=\bes(g), 
\quad \forall g\in \mathbf J
$$
where~$\bes=\bes_{\Lambda}$ is the Bessel function.
\end{prop}

\begin{proof} For Parts~(i) and~(ii) we argue as in the proof 
of~\cite{bh}~Proposition~3.2. If we set~$f(g)=\langle \pi(g) v_{\alpha}, 
\check{v}_{\alpha} \rangle$, then~$f$ satisfies~(i) and~(ii); 
the uniqueness is implied by Corollary~\ref{alphatwo}. 
 
If~$\langle g v_{\alpha}, \check{v}_{\alpha}\rangle\neq 0$ 
then~$e_{\Psi} \pi(g) v_{\alpha}\neq 0$ and, 
since~$v_{\alpha}\in \VV^{\Lambda}$, Corollary~\ref{alphatwo} implies 
$e_{\Lambda} (\pi(g)\mathcal V^{\Lambda})\neq 0$. Hence~$g$ 
intertwines~$\Lambda$ and so~$g\in \mathbf J$.

If~$g\in \mathbf J$ then Theorem~\ref{bessel} 
and Lemma~\ref{pair} imply that 
$$ 
\langle\pi(g)v_{\alpha},\check{v}_{\alpha}\rangle 
=(L(g)\varphi,\check{\varphi})=\bes(g).
$$
\end{proof}

\begin{thm}\label{whitty} Let~$du$ be an invariant Haar measure 
on~$U$, normalised so that~$\int_{U\cap \mathbf J}du=1$. 
Let~$\Upsilon: \pi \rightarrow \Indu{U}{G}{\psi_{\alpha}}$ be a 
linear map given by
$$ 
v \mapsto [g\mapsto \int_U \psi_{\alpha}(u) \langle \pi(u^{-1}g)v, 
\check{v}_{\alpha}\rangle du].
$$
Then~$\Upsilon$ is non-zero and~$G$-equivariant. Moreover, 
$\supp \Upsilon(v_{\alpha})\subseteq U\mathbf J$
and
$$
[\Upsilon(v_{\alpha})](ug)=\psi_{\alpha}(u) \bes(g), \quad 
\forall u\in U, \quad \forall g\in \mathbf J.
$$
Further,~$\supp \Upsilon(v_{\alpha}) \cap \mir_F = U (H^1\cap \mir_F)$
and
$$
[\Upsilon(v_{\alpha})](uh)=\psi_{\alpha}(u) \theta(h), \quad 
\forall u\in U, \quad \forall h\in H^1\cap \mir_F.
$$
\end{thm}
  
\begin{proof} The first assertion follow directly from 
Proposition~\ref{unique}. Now
$$
\supp \Upsilon(v_{\alpha})\cap \mir_F \subseteq (U \mathbf J)\cap 
\mir_F= U( \mathbf J\cap \mir_F)=U (J\cap \rim)\cap \mir_F,
$$
where the last equality follows from Corollary~\ref{rimdec}. The 
second assertion now follows from Proposition~\ref{proper}(iv). 
\end{proof}

\begin{remar} The Whittaker function~$\Upsilon(v_{\alpha})$ above, 
and the bound on the support, can be obtained by integrating the 
matrix coefficient appearing in~\cite{ihes} (see also~\cite{sally}) and 
this is sufficient for the purposes of~\cite{can}. However, the fact 
that we can realize the representation~$\Lambda$ via Bessel functions 
gives us the precise knowledge of~$\supp \Upsilon(v_{\alpha})\cap 
\mir_F$. We use this in~\S\ref{epsi} to compute, in some cases, 
epsilon factors of pairs.
\end{remar} 

\begin{cor}\label{besrho} Let~$\bes_{\rho}$ be the Bessel function 
of~$\rho$. For~$g\in \kA$,  
\begin{equation*} 
\bes_{\rho}(g)= \left\{ \begin{array}{ll} 
   \bes(g) & \textrm{if~$g\in \mathbf J$,}\\
   0 & \textrm{otherwise.}
   \end {array}\right.
\end{equation*}
\end{cor}

\begin{proof} This is implied by the uniqueness of the matrix 
coefficient in Proposition~\ref{unique}.
\end{proof}

%%%%%%%%%%%%%%%%%%%%%%% Multiplicative property %%%%%%%%%%%%%%%%%%%%%%%

\subsection{Multiplicative property}

Our maximal simple types are of the form~$(J,\lambda)$, 
where~$\lambda=\kappa\otimes\sigma$. In this section, we show that the
Bessel function associated to the extension~$\Lambda$ of~$\lambda$
to~$\mathbf J$ can be split as a product of two Bessel function (see
Proposition~\ref{besmult}).

\begin{lem}\label{Sigma} There exists a representation~$\tilde{\kappa}$ 
of~$\mathbf J$ such that~$\tilde{\kappa}|_J\cong \kappa$. Moreover,
given such a representation~$\tilde{\kappa}$, there
exists a unique representation~$\Sigma$ of~$\mathbf J$, such that
$\Sigma|_J\cong \sigma$ and~$\Lambda\cong \tilde{\kappa} \otimes
\Sigma$. 
\end{lem}

\begin{proof} We may extend  the action
of~$J$ on~$\kappa$  to the action of~$F^{\times}J$, by  making some
uniformiser~$\varpi_F$ act trivially. By definition
of~$\kappa$,~\cite{bk}~(5.2.1),~$E^{\times}$, (and hence~$\mathbf
  J$) intertwines~$\kappa$. Since,~$F^{\times} J$ is normal
in~$\mathbf J$ and the quotient~$\mathbf J /(F^{\times}J) \cong
  E^{\times}/F^{\times} \oE^{\times}$ is a cyclic group, we may extend
the action to~$\mathbf J$.  

Now suppose that we are given a representation~$\tilde{\kappa}$
of~$\mathbf J$ such that~$\tilde{\kappa}|J\cong \kappa$. By the same
argument as above, there exists a representation~$\Sigma$ of~$\mathbf
J$ such that~$\Sigma|_{J}\cong \sigma$. Moreover, we may ensure
that~$\omega_{\tilde{\kappa}}(\varpi_F)\omega_{\Sigma}(\varpi_F)=
\omega_{\Lambda}(\varpi_F)$, where~$\omega$ denotes the central
character. Hence
$$
\Lambda|_{F^{\times}J}\cong (\tilde{\kappa}\otimes
\Sigma)|_{F^{\times} J}.
$$
Thus~$\Lambda$ is a direct summand of 
$$
\Indu{F^{\times}J}{\mathbf J}{\tilde{\kappa}\otimes\Sigma}\cong
\tilde{\kappa}\otimes \Sigma\otimes \Indu{F^{\times}J}{\mathbf J}
{\Eins}\cong \bigoplus_{\chi}\tilde{\kappa}\otimes \Sigma\otimes\chi, 
$$
where~$\chi$ runs over characters of~$\mathbf J/ F^{\times} J\cong
E^{\times}/ F^{\times}\oE^{\times}$. Hence, after replacing~$\Sigma$
by some~$\Sigma\otimes \chi$, we may ensure that~$\Lambda\cong
\tilde{\kappa} \otimes \Sigma$. Moreover, by~\cite{bk}~\S6 we know
that~$\Lambda\cong \Lambda\otimes \chi$ implies that~$\chi$ is the
trivial character. Hence, such~$\Sigma$ is unique. 
\end{proof}

Let us now fix some~$\tilde{\kappa}$ as above and let~$\Sigma$ be the
unique representation of~$\mathbf J$, given by
Lemma~\ref{Sigma}. Let~$(U,\psi_{\alpha})$ be as
in~\S\ref{specase}. In particular, we require that~$U\cap B^{\times}$
is a maximal unipotent subgroup of~$B^{\times}$ and we may
write~$\psi_{\alpha}= \psi_a \psi_b$, such that~$\psi_a$ is trivial
on~$U\cap B^{\times}$;~$\psi_b$ is a non-degenerate character
on~$U\cap B^{\times}$, which descends to a non-degenerate character of
a maximal unipotent subgroup of~$\UU(\urlyb)/\UU^1(\urlyb)$
(see Theorem~\ref{wavingit}(iv)). 
 
\begin{lem}\label{context2} Assumption~\ref{assump} (and hence
Proposition~\ref{proper} and Theorem~\ref{bessel}) holds in the following
contexts:
\begin{enumerate}  
\item $\mathcal K=\mathbf J$,~$\tau=\tilde{\kappa}$,~$\mathcal M=
  J^1$,~$\mathcal U=(J^1\cap U)H^1$,~$\Psi=\Theta$, 
where~$\Theta(uh)=\psi_{\alpha}(u) \theta(h)$, for ~$u\in J^1\cap U$
and~$h\in H^1$.
\item $\mathcal K=\mathbf J$,~$\tau=\Sigma$,~$\mathcal M= \rim\cap J$
  or~$\GG_{\curlya}\cap J$,~$\mathcal U=(J\cap
  U)J^1$,~$\Psi=\psi_{b}$.
\end{enumerate}
\end{lem} 

\begin{proof} Part~(i) is just Theorem~\ref{Psikappa}. Since~$\mathbf
  J/ J^1\cong \kB/\UU^1(\urlyb)$ and~$J^1$ acts
trivially on~$\Sigma$,
part~(ii) is given by Theorem~\ref{context}(i) (together with
Corollary~\ref{rimdec}) applied in the level zero case.
\end{proof}

Theorem~\ref{context} and Lemma~\ref{context2} imply that we may
associate Bessel functions to~$\Lambda$,~$\tilde{\kappa}$
and~$\Sigma$, via Definition~\ref{defbes}.

\begin{prop}\label{besmult} We have 
$$ 
\bes_{\Lambda}(g)= \bes_{\tilde{\kappa}}(g)\bes_{\Sigma}(g), \quad
\forall g\in \mathbf J.
$$
\end{prop}

\begin{proof} We have~$\Lambda\cong \tilde{\kappa} \otimes \Sigma$, 
and~$\Psi_{\alpha}=\Theta \psi_b$, see Lemma~\ref{welldfn}. 
Let~$\VV_{\tilde{\kappa}}$ and~$\VV_{\Sigma}$ be the underlying vector 
spaces of~$\tilde{\kappa}$ and~$\Sigma$, respectively. We claim that 
\begin{equation}\label{mult}
e_{\Psi_{\alpha}} (v\otimes w)= ( e_{\Theta} v) \otimes (e_{\psi_b} w), 
\quad \forall v\in \VV_{\tilde{\kappa}}, \quad \forall w\in \VV_{\Sigma}.
\end{equation}
We choose non-zero vectors~$v_{\Theta}\in e_{\Theta}\VV_{\tilde{\kappa}}$ 
and~$w_{\psi_b}\in e_{\psi_b}\VV_{\Sigma}$. It follows 
from Theorem~\ref{Psikappa} that $\dim e_{\Theta} \VV_{\tilde{\kappa}} =1$ 
and from Lemma~\ref{mirabolic} that~$\dim e_{\psi_b} \VV_{\Sigma} =1$. 
Now~$(J\cap U)H^1$ acts on~$v_{\Theta}\otimes w_{\psi_b}$ 
via~$\Psi_{\alpha}$ and hence Theorem~\ref{Psilambda} implies that the set 
$$
\{ c (v_{\Theta}\otimes w_{\psi_b}): 
\bar{c}\in \mir_{\urlyb}J^1/(J\cap U) H^1\},
$$  
is a basis of~$\VV_{\tilde{\kappa}}\otimes \VV_{\Sigma}$, 
where~$c$ denotes a coset representative of a coset~$\bar{c}$. It is 
enough to show the claim \eqref{mult} holds for the elements of this basis.

If~$g\in \mir_{\urlyb} J^1$ and~$e_{\Psi_{\alpha}} g(v_{\Theta}\otimes 
w_{\psi_b})$ is not equal to zero then~$g$ intertwines~$\Psi_{\alpha}$ 
and hence, by Theorem~\ref{Psilambda}, we obtain that~$g\in (J\cap U) H^1$. 
Conversely, if~$g\in (J\cap U) H^1$ 
then~$e_{\Psi_{\alpha}}g(v_{\Theta}\otimes w_{\psi_b})= 
\Psi_{\alpha}(g) v_{\Theta}\otimes w_{\psi_b}$.  

If~$g\in \mir_{\urlyb} J^1$ and~$(e_{\Theta}g v_{\Theta})\otimes 
(e_{\psi_b} g w_{\psi_b})\neq 0$ then~$g$ intertwines~$\psi_b$ and 
hence, by Lemma~\ref{mirabolic},~$g\in (J\cap U) J^1$. Moreover,~$g$ 
intertwines~$\Theta$ and so by 
Theorem~\ref{Psikappa},~$g\in (J\cap U) H^1$. 
If~$g\in (J\cap U)H^1$ then 
$$
(e_{\Theta} g v_{\Theta})\otimes (e_{\psi_b} g w_{\psi_b}) = 
\Theta(g)\psi_b(g) v_{\Theta} \otimes w_{\psi_b}= 
\Psi_{\alpha}(g)v_{\Theta} \otimes w_{\psi_b}.
$$
Hence we obtain the claim \eqref{mult}. 

Now let~$\Theta'$ be the restriction of~$\Theta$ to~$(J^1\cap U) H^1$. 
According to Theorem~\ref{Psikappa},~$\Theta'$ also occurs 
in~$\kappa$ with multiplicity one and hence~$e_{\Theta'}v= e_{\Theta} v$, 
for all~$v\in \VV_{\tilde{\kappa}}$. Hence,
$$
\bes_{\Lambda}(g)=\tr_{\Lambda}(e_{\Psi_{\alpha}} g)= 
\tr_{\tilde{\kappa}}(e_{\Theta'}  g) \tr_{\Sigma}(e_{\psi_b} g)= 
\bes_{\tilde{\kappa}}(g) \bes_{\Sigma}(g).
$$ 
\end{proof}

%%%%%%%%%%%%%%%%%%%%%%%%%%%%%%%%%%%%%%%%%%%%%%%%%%%%%%%%%%%%%%%%%%%%%%
%%%%%%%%%%%%%%%%%%%%%%%% A numerical invariant %%%%%%%%%%%%%%%%%%%%%%%
%%%%%%%%%%%%%%%%%%%%%%%%%%%%%%%%%%%%%%%%%%%%%%%%%%%%%%%%%%%%%%%%%%%%%%

\section{A numerical invariant}\label{numerical}

In this section, we will define a certain numerical invariant which
appears in our formula for epsilon factors in~\S\ref{epsi}. We
continue with the notation of previous sections.

We suppose that~$E=F[\beta]$ is maximal in $A$ and we 
identify~$V$ with~$E$. Let~$\theta\in \mathcal C(\curlya, 0,\beta)$. 
Let~$\FF=\{ V_i:1\le i\le d\}$ be a maximal~$F$-flag 
in~$E$; let~$U$ be a unipotent radical of~$\Aut_F(E)$-stabiliser 
of~$\FF$; let ~$\chi$ be a smooth, non-degenerate character of~$U$, such that 
$$ 
\chi|_{U\cap H^1}=\theta|_{U\cap H^1};
$$
and let~$\psi_E$ be an additive character of~$E$, trivial on~$\pE$, and 
such that 
$$ 
\psi_E(x)= \psi_F(x), \quad \forall x\in F.
$$ 

\begin{defi}\label{numinv} Choose a basis~$\BB=\{x_1,\ldots, x_d\}$ 
of~$E$ over~$F$, which satisfies Definition~\ref{balba}(i),(iii) 
and~(iv) with~$(U,\chi,\psi_E)$. We define~$\nu \in E^{\times}/ 
(1+\pE)$, by
$$ 
\nu=\nu(\theta,\psi_F,\psi_E)= x_d x_1^{-1} \pmod{\pE}.
$$
\end{defi}

\begin{prop}\label{indi}~$\nu(\theta,\psi_F, \psi_E)$ depends only 
on~$\theta$,~$\psi_F$ and~$\psi_E$.
\end{prop}

\begin{proof} Suppose that we have another 
triple~$U'$,~$\chi'$,~$\BB'$ which satisfy the conditions above. 
Then Proposition~\ref{conju} and the first part of the proof 
of Theorem~\ref{wavingit} imply that there exists~$g\in \mathbf J$ 
such that~$U'= U^g$ and~$\chi'=\chi^g$. Since~$E$ is maximal, we 
may write~$g=x h$, where~$x\in E^{\times}$ and~$h\in J^1$. 

Let~$\xi\in E$ be the unique element such that~$\xi h x_1 = h x_d$. 
Then~$(x_1 x_d^{-1} h^{-1} \xi h) \in \GG \cap \mathbf J$, 
where~$\GG=\{g\in \Aut_F(E); g x_1 =x_1\}$. Now, according 
to Corollary~\ref{decem},~$\GG\cap \mathbf J =\GG \cap J^1$ and 
hence the image of~$ ( x_1 x_d^{-1} h^{-1} \xi h)$ 
in~$\mathbf J / J^1\cong E^{\times}/ 1+\pE$, is equal to~$1$. This 
implies that~$\xi\equiv x_d x_1^{-1} \pmod{\pE}$ and hence we may 
assume that~$U=U'$ and~$\chi=\chi'$. The second part 
of Proposition~\ref{balance} implies that~$\nu$ does not depend on 
the choice of basis~$\BB$. 
\end{proof}  
 
\begin{remar}\label{need} Suppose~$E$ is not necessarily maximal. 
Let~$\theta\in \mathcal C(\curlya,\beta,\psi_F)$ and 
let~$\theta_F\in \mathcal C(\curlya(E),\beta, \psi_F)$ be 
the simple character corresponding to~$\theta$ via the 
correspondence of~\cite{bk}~\S3.6, where~$\curlya(E)$ is the 
hereditary~$\oF$-order in~$\End_F(E)$, corresponding to the lattice 
chain~$\{\pE^i: i\in \ZZ\}$. Let $\{x_1,...,x_d\}$ be the $F$-basis
of~$E$ given by Corollary~\ref{goodbasis}, with $x_1=1$. 
Then it follows from the construction in the proof of 
Theorem~\ref{wavingit}, 
that~$\nu(\theta_F,\psi_F,\psi_E)\equiv x_d\pmod{\pE}$. 
\end{remar}

%%%%%%%%%%%%%%%%%%%% Behaviour under tame lifting %%%%%%%%%%%%%%%%%%%%

\subsection{Behaviour under tame lifting}\label{behaviour}

We continue with the assumption that~$E$ is maximal in~$A$ and let us further assume that $E$ is totally wildly ramified. Let~$K$ 
be a tame extension of~$F$. The algebra~$L=K\otimes_F E$ is a field, which is the compositum 
of~$E$ and~$K$.
Let~$\theta\in \mathcal C(\curlya,0,\beta)$; then Bushnell and 
Henniart in~\cite{ihes} define the \emph{tame lift} of~$\theta$, which 
is a simple character~$\theta^K$ of~$H^1_K=H^1(\beta,\curlya_K)$, 
where~$\curlya_K$ is the hereditary~$\oK$-order in~$\End_K(L)$ 
corresponding to the~$\oK$-lattice chain~$\{\pL^j: j\in \ZZ\}$. We 
will investigate how the invariant~$\nu$ varies with the tame lifting. 

Let~$\psi_L$ and~$\psi_K$ be the additive characters of~$L$ and~$K$,
respectively, given by
$$ 
\psi_L(x)= \psi_E(\tr_{L/E} x), \quad \forall x\in L, \quad 
\psi_K(x)=\psi_F(\tr_{K/F} x), \quad \forall x\in K.
$$ 
Since~$K$ is tame over~$F$,~$\psi_K$ has conductor~$\pK$; 
likewise,~$\psi_L$ has conductor~$\pL$.  

Let~$U$ be a maximal unipotent subgroup of~$G$ and 
let~$\FF=\{V_i:1\le i\le d\}$ be a maximal flag corresponding 
to~$U$. Let~$\chi$ be a smooth non-degenerate character of~$U$, and 
let~$\BB=\{x_1,\ldots, x_d\}$ be an~$F$-basis of~$E$, with respect 
to which~$U$ is the group of unipotent upper-triangular matrices and, 
if~$u\in U$ and~$(u_{ij})\in \MM_d(F)$ is a matrix of~$u$ with respect 
to~$\BB$, then 
$$
\chi(u)=\psi_F(\sum_{i=1}^{d-1} u_{i,i+1}).
$$
Set~$\FF_K=\{V_i\otimes_F K: 1\le i\le d\}$ and let~$U_K$ be the 
unipotent radical of the~$\Aut_K(L)$-stabiliser of~$\FF_K$. 
For~$u\in U_K$, write~$(u_{ij})\in \MM_d(K)$ for the matrix of~$u$ 
with respect to~$\{x_1,\ldots, x_d\}$, and 
let~$\chi^K:U_K\rightarrow \CC^{\times}$ be the character given by 
$$
\chi^K(u)=\psi_K(\sum_{i=1}^{d-1} u_{i,i+1}).
$$

\begin{prop}\label{liftbase} We have~$\theta|_{U\cap H^1}= 
\chi|_{U\cap H^1}$ if and only if 
$$
\theta^K|_{U_K\cap H^1_K}= \chi^K|_{U_K\cap H^1_K}.
$$
\end{prop}

\begin{proof} By~\cite{ihes}~Corollary 9.13(iii), tame lifting is transitive 
in the field extension: if~$K'$ is a subfield of~$K$ containing~$F$, 
then 
$$ 
(\theta^{K'})^{K}= \theta^{K}.%, \quad (\chi^{K'})^{K}= \chi^{K}.
$$
So it is enough to prove the Proposition when~$K/F$ is Galois,  
cyclic, and either unramified or totally tamely ramified, as 
in~\cite{ihes}~(12.2). Let~$\Gamma$ be the Galois group of~$K/F$ and 
fix a generator~$\sigma$ of~$\Gamma$. For~$g\in \Aut_K(L)$, 
let~$N_{\sigma}$ be the cyclic norm map, given 
by~$N_{\sigma}g = g \sigma(g)\cdots \sigma^{l-1}(g)$, where $l=[K:F]$. 
Define~$\mathfrak{H}^1_F$,~$\mathfrak{H}^1_K$,~$\mathfrak{U}_F$,%
~$\mathfrak{U}_K$ by
$$ 
H^1= 1+\mathfrak{H}^1_F, \quad H^1_K= 1+ \mathfrak{H}^1_K, \quad 
U=1+ \mathfrak{U}_F, \quad U_K=1+ \mathfrak{U}_K.
$$  
We observe that the proof of~\cite{ihes}~(12.3)~Proposition (including
the results required from~\cite{ihes}\S11) goes through if we 
replace~$\mathfrak{H}^1_F$ with~$\mathfrak{H}^1_F\cap \mathfrak{U}_F$ 
and~$\mathfrak{H}^1_K$ with~$\mathfrak{H}^1_K\cap \mathfrak{U}_K$. We 
obtain the following:
\begin{enumerate} 
\item For~$x\in H^1_K\cap U_K$, there exists~$u\in H^1_K\cap U_K$ 
such that~$y_x= u x \sigma(u)^{-1}$ satisfies~$N_{\sigma} y_x 
\in H^1\cap U$.
\item The map~$x\mapsto N_{\sigma} y_x$ induces a bijection 
between~$\sigma$-conjugacy classes in~$H^1_K\cap U_K$ and conjugacy 
classes in~$H^1\cap U$.
\end{enumerate}     
Now,~\cite{ihes}~(12.6), and the fact that both~$\theta^K$ and~$\chi^K$ 
are stable under~$\Gamma$, imply that
$$ 
\theta^K(x)=\theta^K(y_x)= \theta(N_{\sigma} y_x), \quad 
\chi^K(x)=\chi^K(y_x)=\chi(N_{\sigma}y_x).
$$  
The above coupled with~(ii) gives the Proposition.
\end{proof}

\begin{remar} 
The above Proposition would follow easily from~\cite{bh}~Lemma~2.10, 
if the gap in its proof were fixed.
\end{remar} 

\begin{cor}\label{change} We have~$\nu(\theta^K,\psi_K,\psi_L)\equiv 
\nu( \theta,\psi_F,\psi_E) \pmod{\pL}$.
\end{cor}

\begin{proof} Let~$(U,\chi)$ be such that~$\theta|_{H^1\cap U}=
\chi|_{H^1\cap U}$, let~$\BB=\{x_1,\ldots,x_d\}$ be an~$F$-basis 
of~$E$, which satisfies Definition~\ref{balba}(i),(iii) and~(iv), 
with respect to~$U$,$\chi$,$\psi_F$ and~$\psi_E$. 
Propositions~\ref{indi} and~\ref{liftbase} imply that it is enough 
to show the following: If~$a \in L$ and~$(a_{ij})$ is the matrix 
of~$a\in \End_K(L)$ with respect to~$\BB$ then~$\psi_L(a)= 
\psi_K(a_{dd})$. 

Since~$\psi_L$ and~$\psi_K$ are additive, and~$L= K\otimes_F E$, it 
is enough to prove this for~$a= c\otimes b$, where~$c\in K$ 
and~$b\in E$. Let~$(b_{ij})$ be a matrix of~$b$ with respect 
to~$\BB$ then~$b_{ij}\in F$ and~$a_{ij}=c b_{ij}$. Hence,
$$ 
\psi_K(a_{dd})= \psi_F( b_{dd} \tr_{K/F}c)=\psi_E(b \tr_{K/F} c)
=\psi_E(\tr_{L/E}(bc))= \psi_L(a).
$$
\end{proof}

%%%%%%%%%%%%%%%%%%%%%%%%%%%%%%%%%%%%%%%%%%%%%%%%%%%%%%%%%%%%%%%%%%%%%%
%%%%%% Application to~$\boldsymbol\varepsilon$-factors of pairs %%%%%%
%%%%%%%%%%%%%%%%%%%%%%%%%%%%%%%%%%%%%%%%%%%%%%%%%%%%%%%%%%%%%%%%%%%%%%

\section{Application to~$\boldsymbol\varepsilon$-factors of 
pairs}\label{epsi}

We will use the Whittaker function constructed in Theorem~\ref{whitty} to 
compute~$\varepsilon$-factor of pairs in the following situation: 

As before, let~$[\curlya, n, 0, \beta]$ be a simple stratum in~$A$, 
such that~$e(\urlyb|\oE)=1$. 
Let~$\theta\in \mathcal C (\curlya, 0, \beta)$%, \psi_F)$
 be a simple 
character and let~$\kappa$ and~$\eta$ be representations 
of~$J=J(\beta,\curlya)$ and~$J^1= J^1(\beta,\curlya)$, as 
in Definition~\ref{maxtype}. Let~$\sigma_1$ and~$\sigma_2$ be lifts 
of cuspidal representations of~$J/J^1\cong \UU(\urlyb)/
\UU^1(\urlyb)\cong\GL_r(\kE)$ to~$J$. We allow the 
case~$\sigma_1\cong\sigma_2$. For~$i=1,2$, 
set~$\lambda_i= \kappa\otimes \sigma_i$ and let~$\pi_i$ be a
supercuspidal representation of~$G$ such that~$\pi_i|_J$ 
contains~$\lambda_i$. According to~\cite{bk}~\S6, there exists an 
irreducible representation~$\Lambda_i$ of~$\mathbf J= E^{\times} J$, 
such that~$\Lambda_i|_J\cong \lambda_i$ and 
$$ 
\pi_i\cong \cIndu{\mathbf J} {G} {\Lambda_i}.
$$ 
We fix some extension~$\tilde{\kappa}$ of~$\kappa$ to~$\mathbf J$, as 
in Lemma~\ref{Sigma}. For~$i=1,2$, let~$\Sigma_i$ be the unique 
representation of~$\mathbf J$, also given by Lemma~\ref{Sigma}, such 
that~$\Lambda_i\cong \tilde{\kappa}\otimes \Sigma_i$  
and~$\Sigma_i|_J\cong \sigma_i$; we view~$\Sigma_i$ as a representation 
of~$\kB=E^{\times}\UU(\urlyb)$ and set 
$$ 
\tau_i= \cIndu{\kB} {B^{\times}} {\Sigma_i}.
$$  
Then~$\tau_1$ and~$\tau_2$ are level zero supercuspidal 
representations of~$B^{\times}\cong \GL_r(E)$. Let~$ \curlya(E)$ be 
the hereditary order in~$\End_F(E)$, corresponding to the lattice 
chain~$\{\pE^i: i\in \ZZ\}$. Let~$\theta_F\in 
\mathcal C(\curlya(E),n,0,\beta)$%, \psi_F)$
 be the simple character 
corresponding to~$\theta$ via the correspondence of~\cite{bk}~\S3.6.

\begin{thm}\label{apl} Choose an additive, unitary character~$\psi_E: 
E\rightarrow \CC^{\times}$, such that~$\psi_E$ is trivial on~$\pE$ 
and~$\psi_E(x)=\psi_F(x)$, for all~$x\in F$. Then 
$$ 
\varepsilon(\pi_1\times \check{\pi}_2, s , \psi_F)= 
\zeta\omega_{\tau_1}(\nu^{-r})
\omega_{\tau_2}(\nu^{r})q^{(s-1/2)r v_E(\nu)N/e} 
\varepsilon(\tau_1\times \check{\tau}_2, s,\psi_E),
$$
where:~$\nu=\nu(\theta_F, \psi_F,\psi_E)\in E^{\times} /(1+\pE)$ is 
the invariant defined in Definition~\ref{numinv};~$r=\dim_E(V)$;%
~$\check{\pi}$ denotes the contragredient of~$\pi$;~$q=q_F$ is the 
cardinality of~$\kF$; and~$\zeta= \omega_{\tau_2}(-1)^{r-1} 
\omega_{\pi_2}(-1)^{N-1}$. 
\end{thm}

We remark that, although the representation~$\tau_1$ and~$\tau_2$
depend on the choice of~$\beta$-extension~$\kappa$, and the choice
of~$\tilde{\kappa}$, the~$\varepsilon$-factor in Theorem~\ref{apl} does
not. For a different choice of~$\tilde{\kappa}$ would twist~$\tau_1$
and~$\tau_2$ by the same tamely ramified character~$\chi$ and we have
$$
\varepsilon(\tau_1\chi\times \check{\tau}_2\chi^{-1}, s,\psi_E)=
\varepsilon(\tau_1\times \check{\tau}_2, s,\psi_E).
$$

In \S\ref{twists}, we use our invariant~$\nu$ to recover 
(and generalise) 
certain results in~\cite{can} on the behaviour of~$\varepsilon$-factors 
of pairs under twists by tamely ramified characters.

%%%%%%%%%%%%%%%%%%%%%%%%%%%% Preparation %%%%%%%%%%%%%%%%%%%%%%%%%%%%%

\subsection{Preparation}

Let~$U$ be a maximal unipotent subgroup of~$G$, which is
the~$G$-stabiliser of the maximal flag~$\FF=\{V_i:1\le i\le N\}$, and 
let~$\psi_{\alpha}=\chi\psi_b$ be a smooth non-degenerate character 
of~$U$, as constructed in Theorem~\ref{wavingit}. 
Let~$\LL=\{L_i: i\in \ZZ\}$ be the lattice chain in~$V$ corresponding 
to~$\curlya$. Recall that, in Corollary~\ref{goodbasis}, we constructed 
an~$F$-basis~$\BB_F=\{v_1,\ldots, v_N\}$ of~$V$, with respect to 
which~$U$ is the group of upper triangular matrices and~$\psi_{\alpha}$ 
is the `standard' character. Moreover, the vectors~$v_i$ are of the form 
$$
v_{d(i-1)+j}=  x_d^{i-1} x_j w_i, \quad 1\le i \le r, 
\quad 1\le j \le d,
$$
where~$\BB_E=\{w_1,\ldots,w_r\}$ is an~$E$-basis of~$V$ such 
that~$L_0=\sum_{i=1}^r \oE w_i$, and~$\{x_1,\ldots, x_d\}$ is 
an~$F$-basis of~$E$, which depends on~$\psi_E$ (see 
Corollary~\ref{goodbasis}). Further,~$x_1=1$. Whenever it is required 
of us we will identify~$G$ with~$\GL_N(F)$ via~$\BB_F$. 
Let~$\rw\in G$ be the element defined on the basis~$\BB_F$ by
$$
\rw(v_i)= v_{N-i+1}, \quad 1\le i \le N.
$$
We also define an involution~$\delta:G\rightarrow G$, by 
$$ 
\delta(g)= \rw g^{\top-1} \rw,
$$
where~$g^{\top}$ is the transpose of~$g$ with respect to~$\BB_F$ and $g^{\top-1}= (g^{\top})^{-1}= (g^{-1})^{\top}$. 

We will briefly recall the definition of~$\varepsilon$-factors of 
pairs, using the the formulation of Jacquet, Piatetskii-Shapiro and 
Shalika~\cite{jpss}, rather than Shahidi~\cite{shahidi}.

Let~$\WW_1= \WW(\pi, \psi_{\alpha})$ and~$\WW_2=\WW(\check{\pi}_2, 
\overline{\psi}_{\alpha})$ be the Whittaker models of~$\pi_1$ 
and~$\check{\pi}_2$ respectively. Let~$\mathcal S(F^N)$ be the set 
of compactly-supported, locally constant functions~$\phi:
F^N\rightarrow \CC$. We denote by~$\ee_1= (1,0,\ldots, 0), \ldots ,
\ee_N= (0, \ldots ,0, 1)$ the standard basis of~$F^N$. Given~$W_1\in 
\WW_1$,~$W_2\in \WW_2$ and~$\Phi\in \mathcal S(F^N)$, we define the 
zeta function
$$
Z(W_1,W_2, \Phi,s)= \int_{U\backslash G} W_1(g) W_2(g) 
\Phi(\ee_N g)|\det g|^s dg,
$$
where~$dg$ is a~$G$-equivariant measure on~$U\backslash G$. Note that, 
under our identification of~$G$ with~$\GL_N(F)$, via~$\BB_F$, the 
term~$\ee_N g$ is the~$N^{\rm th}$ row of the matrix of~$g$ with respect 
to~$\BB_F$. The integral converges absolutely for~$\Re(s)$ sufficiently 
large, and is a rational function of~$q^{-s}$. This zeta function  
satisfies a functional equation,~\cite{jpss}(2.7):
$$ 
\frac{Z(\widetilde{W}_1, \widetilde{W}_2,\hat{\Phi}, 1-s)}
{L(\check{\pi}_1\times \pi_2, 1-s)}= \omega_{\pi_2}(-1)^{N-1} 
\varepsilon(\pi_1\times \check{\pi}_2, s, \psi_F)
\frac{Z(W_1, W_2, \Phi, s)}{L(\pi_1\times \check{\pi}_2, s)},
$$
where, for~$i=1,2$,~$\widetilde{W}_i(g)= W_i(\rw g^{\top-1})$;%
~$\hat{\Phi}$ is the Fourier transform of~$\Phi$, given by
$$ 
\hat{\Phi}(\mathrm y)= \int_{F^N}\Phi(\mathrm x) 
\psi_F(\mathrm x%\centerdot 
\mathrm y^{\top}) d \mathrm x, \quad 
\forall \mathrm y\in F^N,
$$
where~$d\mathrm x$ is normalised so 
that~$\hat{\hat {\Phi}}(\mathrm x)= \Phi(-\mathrm x)$; 
and~$L(\pi_1\times\check{\pi}_2,s)$ is the~$L$-function. 
Since~$\pi_1$ and~$\pi_2$ are supercuspidal, it is enough for our 
purposes to know that 
$$ 
L(\pi_1\times\check{\pi}_2,s)= \prod_{\chi} L(\chi,s),
$$
where the product is taken over all the unramified 
characters~$\chi:F^{\times}\rightarrow \CC^{\times}$ such 
that~$\pi_1\cong \pi_2\otimes \chi\circ \det$, and~$L(\chi,s)$ is 
as in Tate's thesis~\cite{tate}. In particular, 
if~$\sigma_1\not\cong \sigma_2$ then the product is taken over an 
empty set and so $L(\pi_1\times\check{\pi}_2,s)=1$. 
If~$\sigma_1\cong \sigma_2$ then there exists~$\chi$ as above such 
that~$\pi_1\cong \pi_2\otimes\chi\circ \det$. In this case, it follows 
from~\cite{bk}~(6.2.5) that 
$$
L(\pi_1\times \check{\pi}_2, s)= (
1- \chi(\varpi_F)^{-N/e} q^{-s N/e})^{-1},
$$
where~$e=e(\curlya|\oF)= e(E|F)$.

%%%%%%%%%%%%%%%%%%%%%%%%%%%% Computation %%%%%%%%%%%%%%%%%%%%%%%%%%%%%

\subsection{Computation}

Let~$W_1\in \WW(\pi_1, \psi_{\alpha})$ and~$W_2\in \WW(\check{\pi}_2, 
\overline{\psi}_{\alpha})$ be the Whittaker functions constructed 
in Theorem~\ref{whitty}. Then~$\supp W_1\subseteq U\mathbf J$,%
~$\supp W_2\subseteq U\mathbf J$ and 
$$ 
W_1(ug)= \psi_{\alpha}(u)\bes_{\Lambda_1}(g), \quad W_2(ug)= 
\overline{\psi}_{\alpha}(u) \bes_{\check{\Lambda}_2}(g), 
\quad \forall u\in U, \quad \forall g\in \mathbf J.
$$
Set~$\bes_1= \bes_{\Lambda_1}$ and~$\bes_2=\bes_{\check{\Lambda}_2}$, and 
let~$\Phi\in \mathcal S(F^N)$ be the indicator function on the 
set~$\ee_N J^1$. We are going to compute the zeta functions on both 
sides of the functional equation for this particular choice 
of~$W_1$,~$W_2$ and~$\Phi$. This will give us Theorem~\ref{apl}. 

%\marginpar{\tiny Have I got this volume right here? See also~$\volF$ below.}
For~$X$ a subset of~$G$ which is a union of right~$U$-cosets, we
write~$\volU(X)$ for the volume of~$U\backslash X$ with respect to the
measure~$du$ on~$U\backslash G$.

\begin{prop}\label{right}Let~$\rF: G\rightarrow \CC$ be the function 
given by
$$
\rF(g)= W_1(g) W_2(g) \Phi(\ee_N g).
$$
Then~$\rF$ is an indicator function on the set~$U H^1$. In particular, 
$$
Z(W_1,W_2,\Phi, s)=\volU (U H^1).
$$
\end{prop} 

\begin{proof} We have 
$$
\supp W_1\subseteq U\mathbf J, \quad  \supp W_2 \subseteq U\mathbf J, 
\quad \supp [g\mapsto \Phi(\ee_N g)] = \mir_F J^1,
$$
where~$\mir_F=\{g\in G: (g-1) V\subseteq V_{N-1}\}$ is the mirabolic
subgroup, as in~\S\ref{super}. Hence,
$$ 
\supp \rF \subseteq U\mathbf J \cap \mir_F J^1= 
U (\mathbf J \cap \mir_F) J^1= U \mir_{\urlyb}J^1,
$$
where the last equality is given by Corollary~\ref{rimdec}. We have 
$$ 
\rF(ug)= W_1(g)W_2(g)\Phi(\ee_N g)= \bes_1(g) \bes_2(g), 
\quad \forall u\in U, \forall g\in \mir_{\urlyb}J^1.
$$
It follows from Proposition~\ref{proper}(ii) and~(iv) 
that~$\bes_1(g)\bes_2(g)=\Psi_{\alpha}(g)\overline{\Psi_{\alpha}}(g)=1$, 
if~$g\in (J\cap U)H^1$, and~$\bes_1(g) \bes_2(g)=0$, otherwise. 
Hence~$\rF$ is an indicator function on the set~$U H^1$. Since~$H^1$ 
is compact and~$U$ is unipotent, we obtain that~$|\det g |=1$, for 
all~$g\in U H^1$. Hence~$Z(W_1, W_2, s,\Phi)= \volU (UH^1)$.  
\end{proof}

We write $\GG_F=\{g\in G:gv_1=v_1\}$, as in~\S\ref{super}.
\begin{lem}\label{ccos} For all  $g_1\in \GG$, $h\in G$ and $ g_2\in (J\cap \mir_F)J^1$ we have:
$$\hat{\Phi}(\ee_1 (g_1h g_2)^{\top-1})=\hat{\Phi}(\ee_1 h^{\top-1}).$$
\end{lem}
\begin{proof} Since $\ee_1 g_1^{\top-1} =\ee_1$, we obtain 
$\hat{\Phi}(\ee_1 (g_1 h)^{\top-1})= \hat{\Phi}(\ee_1 h^{\top-1})$. Since $\Phi$ is an indicator 
function on the set $\ee_N J^1= \ee_N (J\cap \mir_F)J^1$, we have $g_2 \Phi =\Phi$, hence 
\begin{equation}\notag
\begin{split}
 \hat{\Phi}(\ee_1 (hg_2)^{\top-1})=&\int_{F^N} \Phi(\mathrm x) \psi_F(\mathrm x  (g_2^{-1} h^{-1}\ee_1^{\top}))d\mathrm x\\
                                =&\int_{F^N} [g_2\Phi](\mathrm x) \psi_F(\mathrm x  (h^{-1}\ee_1 ^{\top}))d\mathrm x=\hat{\Phi}(\ee_1 h^{\top-1}).
\end{split}
\end{equation} 
\end{proof}

For~$L$ a lattice in~$F^N$, we write~$\volF(L)$ for the volume of~$L$
with respect to the measure~$d{\mathrm x}$ on~$F^N$.

\begin{lem}\label{volume} Let~$i,j\in \ZZ$; 
then~$\volF(\ee_N \curlyp^i)= q^{(j-i)N/e} \volF(\ee_N \curlyp^j)$.
\end{lem}

\begin{proof} Since there exists~$\gamma\in \kA$ such that~$\ee_N=
\ee_1 \gamma$, and then $\ee_N \curlyp^i= \ee_1 \gamma \curlyp^i =\ee_1 
\curlyp^{i+v_{\curlya}(\gamma)}$, it is enough to prove that 
$$
\volF(\ee_1 \curlyp^i)=q^{(j-i)N/e} \volF(\ee_1 \curlyp^j).
$$ 
Since our basis~$\BB_F$ splits the lattice chain we may write
$$ 
\curlyp^{k}= \bigoplus_{1\le i,j\le N} \pF^{c_{ij}(k)} \Eins_{ij},
$$
where~$\Eins_{ij}\in A$ are the projections given 
by~$\Eins_{ij}(v_k)= \delta_{ik} v_j$, for~$1\le k \le N$, 
and~$\delta_{ik}$ is the Kronecker delta.
Hence,~$\ee_1 \curlyp^k= \sum_{j=1}^N \pF^{c_{1j}(k)} \ee_j$. According 
to~\cite{bk}~(1.1.4) we have 
$$
\curlyp^{1-k} =\{a\in A: \psi_F( \tr_{A/F}( xa))=1, 
\forall x\in\curlyp^k\}, \quad \forall k\in \ZZ,
$$
which, since~$\psi_F$ has conductor~$\pF$, implies that~$1-c_{ij}(k)=
c_{ji}(1-k)$, for all~$k\in \ZZ$.  
Since we have chosen~$v_1$, so that~$v_1\in L_0$,~$v_1\not\in L_1$ and 
the lattice chain is principal, we have~$\curlyp^k v_1= L_k$, for 
all~$k\in \ZZ$. Now,~$\curlyp^k v_1 = \sum_{j=1}^N \pF^{c_{j1}(k)} v_j= 
\sum_{j=1}^N \pF^{1- c_{1j}(1-k)} v_j.$ Hence
$$
(\ee_1 \curlyp^i : \ee_1 \curlyp^j)= (L_{1-j} : L_{1-i}),
$$
where the brackets denote the generalised index. Since~$(L_i:
L_{i+e})=q^N$ and~$\LL$ is principal, we have~$(L_i:
L_{i+1})=q^{N/e}$, for all~$i\in \ZZ$, and the result follows. 
\end{proof}

We write~$q_{\curlya}=q^{N/e}$, so Lemma~\ref{volume} says
that~$\volF(\ee_N \curlyp^i)= q_{\curlya}^{j-i} \volF(\ee_N \curlyp^j)$.

Let~$\rw_E \in \UU(\urlyb)$ be the element defined by its action 
on the basis~$\BB_E$ by 
$$ 
\rw_E(w_i)= w_{r-i+1}, \quad 1\le i \le r.
$$
From our construction of the bases~$\BB_E$ and~$\BB_F$, we  
have~$x_d^r \rw_E v_1 = x_d^r w_r = v_N$. In terms of matrices with 
respect to~$\BB_F$ we can rephrase this as
$$ 
(x_d^r \rw_E) \ee_1^{\top}= \ee_N^{\top}.
$$

\begin{lem}\label{addcond} Let~$\phi: A\rightarrow \CC^{\times}$ be the 
function 
$$
\phi: a\mapsto \psi_F((\ee_N a)%\centerdot 
\ee_1^{\top})=\psi_F(a_{N1}),
$$
where~$(a_{ij})$ is the matrix of~$a\in A$ with respect to~$\BB_F$. 
Then 
$$ 
\phi( a x_d^r \rw_E)= \psi_F(a_{NN}), \quad \forall a\in A.
$$
Hence~$\phi$ defines an additive character on~$A$, which is trivial 
on~$\curlyp^{1+r v_{E}(x_d)}$, and non-trivial on~$\curlyp^{rv_{E}(x_d)}$.
\end{lem}

\begin{lem}\label{condE} Let~$b\in B$, let~$(b_{ij})$ be the matrix 
of~$b$ with respect to~$\BB_E$, and define~$\phi: A\rightarrow
\CC^{\times}$ as in Lemma~\ref{addcond}; then
$$
\phi(b)= \psi_E(x_d^{-r} b_{r1}).
$$
\end{lem}

\begin{proof} Set~$a=b x_d^{-r} \rw_E$, let~$(a^F_{ij})$ be the 
matrix of~$a$ with respect to~$\BB_F$ and let~$(a^E_{ij})$ be the 
matrix of~$a$ with respect to~$\BB_E$. According 
to Lemma~\ref{addcond}, we have~$\phi(b)= \psi_F( a^F_{NN})$. We 
have~$ a w_r + V_{N-r}=  a^E_{rr} w_{r}+ V_{N-r}$ and, 
since~$x_d\in E$, we obtain~$ a(x_d^r w_r) +V_{N-r} = 
a_{rr}^E(x_d^r w_r) + V_{N-r}.$ Since~$a^E_{rr}\in E$, we may consider 
it as~$a^E_{rr}\in \End_F(E)$. Let~$(\alpha_{ij})$ be the matrix 
of~$a^E_{rr}$ with respect to the basis~$\{x_1, \ldots, x_d\}$. 
Since~$v_{d(i-1)+j}=  x_d^{i-1} x_j w_i$, for~$1\le i \le r$ 
and~$1\le j \le d$, and in particular~$v_N= x_d^r w_r$, we obtain that  
$$ 
a v_N + V_{N-1} = \alpha_{dd} v_N + V_{N-1}.
$$
In particular,~$\alpha_{dd}= a_{NN}$. Recall that the 
basis~$\{x_1, \ldots x_d\}$ was chosen so 
that~$\psi_E(a^E_{rr})=\psi_F(\alpha_{dd})$, see 
Definition~\ref{balba}(iv). Hence,~$\phi(b)= \psi_E(a^E_{rr})$. 
Since~$x_d\in E$, we obtain that~$a^E_{rr}= x_d^{-r} (b \rw_E )_{rr} 
= x_d^{-r} b_{r1}$. 
\end{proof}

To ease the notation, we set 
$$ 
c=\volF(\ee_N \curlyp^{1+rv_{E}(x_d)}).
$$
As in~\S\ref{deecthm} let~$K$ be a maximal unramified extension 
of~$E$, such that~$K^{\times}$ normalises~$\curlya$.

\begin{lem}\label{phihat} Let~$h\in \kA$ and 
set~$j= v_{\curlya}(h)- v_E(x_d^{-r})$; then 
\begin{equation*} \hat{\Phi}(\ee_1 h^{\top-1})= 
\left\{ \begin{array}{ll}
   0 & \textrm{ if }j>0;\\
   c|\det h|\phi(h^{-1})  & \textrm{ if } j=0;\\    
   c|\det h| q_{\curlya}^{j} & \textrm{ if } j<0.
        \end{array} \right. 
\end{equation*}
\end{lem}

\begin{proof} Corollary~\ref{decem} implies that~$\ee_N \UU^1(\curlya)=
\ee_N J^1=\ee_N (1+\pK)$. Hence, 
\begin{equation*}
\begin{split}
\hat{\Phi}(\ee_1 h^{\top-1})
&=\int_{\ee_N \UU^1(\curlya)} \psi_F( (\mathrm x  h^{-1})
       %\centerdot 
       \ee_1^{\top})d\mathrm x\\
&=\psi_F( (\ee_N h^{-1})%\centerdot 
                        \ee_1^{\top}) \int_{\ee_N \curlyp} 
        \psi_F((\mathrm x h^{-1})%\centerdot 
                                 \ee_1^{\top})d \mathrm x\\  
&=|\det h| \psi_F( (\ee_N h^{-1})%\centerdot 
                                 \ee_1^{\top})
        \int_{\ee_N \curlyp^{1-v_{\curlya}(h)}}
        \psi_F(\mathrm x %\centerdot 
                         \ee_1^{\top})d \mathrm x.
\end{split}
\end{equation*}
Lemmas~\ref{volume},~\ref{addcond} and the orthogonality of characters 
imply the Lemma. 
\end{proof}

\begin{lem}\label{voldelta}~$\volU (U \delta(H^1))= \volU( U H^1)$.
\end{lem} 

\begin{proof} Let~$\mathcal K_m=\{g\in \GL_N(\oF): 
g \equiv 1 \pmod{\pF^m}\}$. Since~$\mathcal K_m$, for~$m\ge 1$, 
form a basis of neighbourhoods of~$1$ in~$G$, there exists~$m$, such 
that~$\mathcal K_m\subseteq H^1\cap \delta(H^1)$. 
Since~$\delta(U)=U$,~$\delta(\mathcal K_m)=\mathcal K_m$ and the 
measure on~$U\backslash G$ is~$G$-invariant, we obtain 
that~$\volU (U \delta(H^1))= \volU( U H^1)$.
\end{proof}  

Let~$\widetilde{\rF}: G\rightarrow \CC$ be the function given by 
$$ 
\widetilde{\rF}(g)= 
\widetilde{W}_1(g)\widetilde{W}_2(g) \hat{\Phi}(\ee_Ng).
$$
We have~$\supp \widetilde{\rF}\subseteq \supp \widetilde{W}_1 
= \delta(\supp W_1) \rw\subseteq U \delta(\mathbf J) \rw$, and 
$$ 
\widetilde{\rF}(u \delta(g)\rw)=\bes_1(g)\bes_2(g) 
\hat{\Phi}(\ee_1 g^{\top-1}), \quad 
\forall g\in\mathbf J, \forall u\in U.
$$
For~$x\in E^{\times}$ we define~$S(x)$, by 
$$ 
S(x)=\int_{U\delta(Jx)\rw} \widetilde{\rF}(g) |\det g|^{1-s} dg.
$$
Then~$S(x)$ depends only on~$v_E(x)$. For~$\Re(-s)$ sufficiently large, 
%\marginpar{\tiny Does this mean $\Re(1-s)$ sufficiently large?}
$$
Z(\widetilde{W}_1,\widetilde{W}_2,\hat{\Phi}, 1-s)=
\sum_{x\in \oE^{\times} \backslash E^{\times}} S(x).
$$
Since~$H^1$ is normal in~$\mathbf J$, Proposition~\ref{proper}(iii), 
and Lemma~\ref{ccos} imply that 
$$ 
\widetilde{\rF}(\delta(h)g)= \widetilde{\rF}(g), \quad 
\forall h\in H^1, \forall g\in G.
$$
Hence, 
$$ 
S(x)=\volU(U \delta(H^1)) |\det x|^{s-1}
\sum_{h\in (J\cap U) H^1\backslash J} 
\bes_1(hx) \bes_2(hx)\hat{\Phi}(\ee_1 (hx)^{\top-1}).
$$
We forget the volume term, by using Lemma~\ref{voldelta} and 
normalising the 
measure on~$U\backslash G$, so that~$\volU(U\delta(H^1))=\volU(U
H^1)=1$.

As in~\S\ref{deecthm}, we put
$$
\GG_E=\{g\in B^{\times}: g w_1 =w_1\}, \quad \GG_{\urlyb}=(\GG_E\cap
\UU(\urlyb)) \UU^1(\urlyb),
$$
where~$\BB_E=\{w_1,...,w_r\}$ is our~$E$-basis of~$V$. 
Corollary~\ref{decem} implies that we 
have~$J= (\GG_{\urlyb} J^1) \oK^{\times}$. 
Then, using Lemma~\ref{ccos}, we obtain:

\begin{lem}\label{state}
$$ 
S(x)= |\det x|^{s-1} \sum_{y\in (1+\pK)\backslash \oK^{\times}} 
\hat{\Phi}(\ee_1 (yx)^{\top-1})
\sum_{h\in(J\cap U)H^1\backslash \GG_{\urlyb} J^1}  
      \bes_1(h yx) \bes_2(hyx).
$$
\end{lem}

%%%%%%%%%%%%%%%%%% The case~$\sigma_1\cong\sigma_2$ %%%%%%%%%%%%%%%%%%

\subsection{The case~$\sigma_1=\sigma_2$}

Suppose that~$\sigma_1\cong \sigma_2$; then it follows 
from~\cite{bk}~(6.2.3) that there exists an unramified 
quasi-character~$\chi: F^{\times}\rightarrow \CC^{\times}$, such 
that~$\Lambda_1\cong \Lambda_2\otimes \chi\circ\det$, and 
hence~$\pi_1\cong \pi_2\otimes \chi\circ\det$. 
Then~$\bes_2(g)=\bes_1(g^{-1})\chi(\det g)$, for all~$g\in \mathbf J$. 
Hence, for all~$g\in \mathbf J$, we have  
$$ 
\sum_h \bes_1(h g) \bes_2(h g)=
\chi(\det g)\sum_{h\in(J\cap U)H^1\backslash \GG_{\urlyb} J^1 }
\bes_1(h g) \bes_1(g^{-1} h^{-1})= \chi(\det g),
$$
where the last equalities follow from Proposition~\ref{proper}(v),(i). 
It follows from Corollary~\ref{decem} 
that~$(J: \GG_{\urlyb} J^1)= (\oK^{\times}: 1+\pK)= 
q^{N/e}-1=q_{\curlya}-1$. There exists~$a\in \CC$ such 
that~$\chi(x)= |x|_F^a$, for all~$x\in F$. If~$x\in E$ then 
Lemma~\ref{detv} 
implies that~$\chi(\det x)= q_{\curlya}^{-a v_E(x)}$. 
Let~$x\in E^{\times}$ and set~$j=v_E(x)-v_E(x_d^{-r})$, 
Lemmas~\ref{detv} and~\ref{phihat} imply that
\begin{equation*} 
S(x)= \left\{ \begin{array}{ll}
  0 & \textrm{if }j> 0;\\
  cq_{\curlya}^{ rv_{E}(x_d)(s-a)}\sum_{y} \phi(y x_d^r)  & 
      \textrm{if } j=0;\\    
  cq_{\curlya}^{-(s-1-a)j +(s-a) rv_{E}(x_d)}(q_{\curlya}-1)  & 
      \textrm{if } j<0.
\end{array} \right. 
\end{equation*}
where, in the sum,~$y$ runs over the cosets~$\oK^{\times}/1+\pK$ and~$\phi$ is 
defined in Lemma~\ref{addcond}. It follows from~\ref{addcond} that~$\phi$ 
restricted to~$K$ defines an additive character, which is trivial 
on~$\pK^{rv_E(x_d)+1}$ and non-trivial on~$\pK^{rv_E(x_d)}$. Hence
$$
\sum_{y\in \oK^{\times}/1+\pK} \phi(y x_d^r)= -\phi(0)=-1.
$$  
%\marginpar{\tiny Again, does this mean $-\Re(s)$ sufficiently large?}
Set~$\widetilde{Z}=Z(\widetilde{W}_1,\widetilde{W}_2,\hat{\Phi},1-s)$; 
then, for~$\Re(-s)$ sufficiently large, we obtain that 
\begin{equation*}
\begin{split}
\widetilde{Z}&= c q_{\curlya}^{rv_{E}(x_d)(s-a)} (-1 + (q_{\curlya}-1) 
  \sum_{k\ge 1} q_{\curlya}^{(s-1-a)k})\\
&= c q_{\curlya}^{rv_{E}(x_d)(s-a)} 
  \frac{q_{\curlya}^{s-a }-1}{1-q_{\curlya}^{s-1-a}}
=c q_{\curlya}^{(r v_{E}(x_d)+1)(s-a)}
\frac{L(\check{\pi}_1\times \pi_2, 1-s)}{L(\pi_1\times\check{\pi}_2, s)}.  
\end{split}
\end{equation*}
It follows from the functional equation and Proposition~\ref{right} that 
$$
\varepsilon(\pi_1\times \check{\pi}_2,\psi_F, s)= 
\omega_{\pi_2}(-1)^{N-1} c q_{\curlya}^{(s-a) (rv_{E}(x_d)+1)}.
$$
Following~\cite{calc}, we observe that the symmetry in the functional 
equation implies 
that~$\varepsilon(\pi_1\times \check{\pi}_1,\psi_F, 1/2)^2=1$. Hence, 
$$ 
c=\volF(\ee_N \curlyp^{1+r v_{E}(x_d)})= q_{\curlya}^{-(rv_{E}(x_d)+1)/2}.
$$ 

We will now prove Theorem~\ref{apl} in the case~$\sigma_1\cong\sigma_2$.
\begin{proof}
Note that~$\Lambda_1\cong \Lambda_2\otimes \chi\circ\det_A$ 
implies~$\Sigma_1\cong \Sigma_2\otimes \chi_E\circ \det_B$, 
where~$\chi_E=\chi\circ \mathrm N_{E/F}$ and $\mathrm N_{E/F}$ denotes
the field norm. Moreover, Remark~\ref{need} 
implies that~$v_E(x_d)= v_E(\nu)$, 
where~$\nu=\nu(\theta_F,\psi_F,\psi_E)$. Hence, 
$$
q^{av_E(x_d^{-r})N/e}= \chi(\det \nu^{-r})= 
\omega_{\tau_1}(\nu^{-r})\omega_{\tau_2}(\nu^r).
$$
If we compute~$\varepsilon(\tau_1\times\check{\tau}_2,s,\psi_E)$ by 
the same recipe we obtain that 
$$ 
\varepsilon(\tau_1\times\check{\tau}_2,s,\psi_E)= 
\omega_{\tau_2}(-1)^{r-1}q_{\urlyb}^{-1/2} q_{\urlyb}^{s-a},
$$
where~$q_{\urlyb}= q_E^{r}= q_F^{N/e}= q_{\curlya}$. Hence,
$$ 
\varepsilon(\pi_1\times \check{\pi}_2, s , \psi_F)= 
\zeta\omega_{\tau_1}(\nu^{-r})\omega_{\tau_2}(\nu^{r})
q^{(s-1/2)r v_E(\nu)N/e} 
\varepsilon(\tau_1\times \check{\tau}_2, s,\psi_E),
$$
where~$\zeta= \omega_{\pi_2}(-1)^{N-1}\omega_{\tau_2}(-1)^{r-1}$.
\end{proof}

%%%%%%%%%%%%%%%% The case~$\sigma_1\not\cong\sigma_2$ %%%%%%%%%%%%%%%%

\subsection{The case~$\sigma_1\not=\sigma_2$}

Now let us suppose that~$\sigma_1\not\cong \sigma_2$. 

\begin{lem}\label{vanish} Let~$x\in E^{\times}$. 
If~$v_{E}(x)\neq v_{E}(x_d^{-r})$ then~$S(x)=0$.
\end{lem}

\begin{proof} Set~$j= v_E(x)-v_E(x_d^{-r})$. If~$j>0$ then 
Lemma~\ref{phihat} gives us~$S(x)=0$. If~$j< 0$ then 
Lemma~\ref{phihat} implies that 
$$
S(x)=c|\det x|^s q_{\curlya}^{j}\sum_{h\in (J\cap U)H^1 \backslash J} 
\bes_1(hx) \bes_2(hx).
$$
Recall from~\S\ref{besfun} that we have the
idempotent~$e_{\Psi_{\alpha}}$ given by
$$
e_{\Psi_{\alpha}}= Q^{-1} 
\sum_{h\in \UU^{n+1}(\curlya)\backslash (J\cap U) H^1} \Psi_{\alpha}(h) h^{-1},
$$
where~$Q=((J\cap U)H^1:\UU^{n+1}(\curlya))$, 
and similarly~$e_{\overline{\Psi}_{\alpha}}$. Now, 
\begin{equation*}
\begin{split}
\sum_{h\in (J\cap U)H^1 \backslash J} \bes_1(hx) \bes_2(hx)
&= \sum_{h\in (J\cap U) H^1 \backslash J} 
  \tr_{\Lambda_1}(x e_{\Psi_{\alpha}} h) 
  \tr_{\check{\Lambda}_2}(x e_{\overline{\Psi}_{\alpha}} h)\\
&=Q^{-1}\tr_{\Lambda_1\otimes\check{\Lambda}_2}
   \biggr ( x (e_{\Psi_{\alpha}}\otimes e_{\overline{\Psi}_{\alpha}}) 
   \sum_{h\in \UU^{n+1}(\curlya)\backslash J} h \biggl ), 
\end{split}
\end{equation*}   
Since~$\sigma_1\not\cong \sigma_2$ we  
have~$\Lambda_1\not\cong \Lambda_2$, 
hence~$\Hom_{J}( \Eins, \Lambda_1\otimes\check{\Lambda}_2)=0$. This 
implies that 
$$ 
\sum_{h\in \UU^{n+1}(\curlya)\backslash J} 
\Lambda_1(h)\otimes \check{\Lambda}_2(h)=0,
$$
and so~$S(x)=0$.
\end{proof}

\begin{lem}\label{cancel} Let~$b\in \kB$; then 
$$ 
\sum_{h\in (J^1\cap U)H^1\backslash J^1} \bes_1(hb)\bes_2(hb)= 
\bes_{\Sigma_1}(b) \bes_{\check{\Sigma}_2}(b).
$$
\end{lem}

\begin{proof} According to Proposition~\ref{besmult}, we 
have~$\bes_{\Lambda_i}(g)=\bes_{\tilde{\kappa}}(g) \bes_{\Sigma_i}(g)$, 
for~$i=1,2$ and for all~$g\in \mathbf J$. The assertion follows from 
the fact that~$J^1$ acts trivially on~$\Sigma_i$ and 
Proposition~\ref{proper}(v),(i) applied to~$\tilde{\kappa}$, 
via Lemma~\ref{context2}.
\end{proof}  

We now prove Theorem~\ref{apl} in the case 
when~$\sigma_1\not\cong \sigma_2$.
\begin{proof}
Set~$\widetilde{Z}=
Z(\widetilde{W}_1,\widetilde{W}_2,\hat{\Phi}, 1-s)$. It follows from 
Lemma~\ref{vanish} that, for~$\Re(-s)$ sufficiently  
%\marginpar{\tiny Does this mean $-\Re(s)$ sufficiently large again?}
large,~$\widetilde{Z}=S(x_d^{-r})$. Lemmas~\ref{state},~\ref{detv} 
and~\ref{phihat} imply that
$$
\widetilde{Z}= c q_{\curlya}^{-r v_E(x_d)s}
\sum_{y\in (1+\pK)\backslash \oK^{\times}} 
\phi(y^{-1} x_d^{r})\sum_{h\in(J\cap U)H^1\backslash \GG_{\urlyb} J^1} 
\bes_1(h yx_d^{-r}) \bes_2(hyx_d^{-r}),
$$
where~$c= q_{\curlya}^{(-r v_E(x_d)-1)/2}$. Lemma~\ref{cancel} implies 
that 
$$ 
\widetilde{Z}= c q_{\curlya}^{r v_E(x_d)s}
\sum_{y\in (1+\pK)\backslash \oK^{\times}} \phi(y^{-1} x_d^{r})
\sum_{h\in \UU^1(\urlyb_m)\backslash \GG_{\urlyb}} 
\bes_{\Sigma_1}(hyx_d^{-r}) \bes_{\check{\Sigma}_2}(hyx_d^{-r}),
$$
where~$\UU^1(\urlyb_m)= (U\cap \UU(\urlyb)) \UU^1(\urlyb)$. 
Now~$x_d\in E$ so we can use Proposition~\ref{proper}(ii) and 
Lemma~\ref{ccos} to obtain 
$$ 
\widetilde{Z}= c q_{\curlya}^{r v_E(x_d)s} 
\omega_{\Sigma_1}(x_d^{-r})\omega_{\Sigma_2}(x_d^r)
\sum_{h\in \UU^1(\urlyb_m)\backslash \UU(\urlyb)} \phi(h x_d^r)
\bes_{\Sigma_1}(h^{-1}) \bes_{\check{\Sigma}_2}(h^{-1}).
$$
Lemma~\ref{condE} implies that 
$$ 
\widetilde{Z}= c q_{\curlya}^{r v_E(x_d)s} 
\omega_{\Sigma_1}(x_d^{-r})\omega_{\Sigma_2}(x_d^r)
\sum_{h\in \UU^1(\urlyb_m)\backslash \UU(\urlyb)} \psi_E(h_{r1})
\bes_{\Sigma_1}(h^{-1})\bes_{\check{\Sigma}_2}(h^{-1}),
$$
where~$h_{r1}$ is the~$r1$-coefficient of the matrix of~$h$ with 
respect to the basis~$\BB_E$. It now follows from the functional equation 
and Proposition~\ref{right} that~$\varepsilon(\pi_1\times\check{\pi}_2,
s,\psi_F)= \omega_{\pi_2}(-1)^{N-1} \widetilde{Z}$. 

If we 
compute~$\varepsilon(\tau_1\times\check{\tau_2},s, \psi_E)$ by the same 
recipe, we obtain that
$$ 
\varepsilon(\tau_1\times\check{\tau}_2,s, \psi_E)= 
q_{\urlyb}^{-1/2} \omega_{\tau_2}(-1)^{r-1}
\sum_{h\in \UU^1(\urlyb_m)\backslash \UU(\urlyb)} \psi_E(h_{r1})
\bes_{\Sigma_1}(h^{-1})\bes_{\check{\Sigma}_2}(h^{-1}),
$$
where~$q_{\urlyb}= q_E^{r}= q_F^{N/e}= q_{\curlya}$. Hence,
$$
\varepsilon(\pi_1\times \check{\pi}_2, s, \psi_F)= 
\zeta \omega_{\tau_1}(x_d^{-r})\omega_{\tau_2}(x_d^r)  
q^{rv_E(x_d)(s-1/2)N/e}
\varepsilon(\tau_1\times\check{\tau}_2,s, \psi_E),
$$
where~$\zeta= \omega_{\tau_2}(-1)^{r-1} \omega_{\pi_2}(-1)^{N-1}$. 
Now~$\omega_{\tau_1}$ and~$\omega_{\tau_2}$ are trivial on~$1+\pE$ 
and Remark~\ref{need} finishes the proof, as in the
case~$\sigma_1\cong\sigma_2$. 
\end{proof}

%%%%%%%%%%%%%%% Twists by tamely ramified quasi-characters %%%%%%%%%%%

\subsection{Twists by tamely ramified quasi-characters}\label{twists}

We continue in the same situation as above. Theorem~\ref{apl} 
immediately implies the following:

\begin{cor}\label{twistcor}
 Let~$\chi:F^{\times}\rightarrow \CC^{\times}$ be a tamely 
ramified quasi-character and put~$\chi_E=\chi\circ \mathrm{N}_{E/F}$; then  
$$
\frac{\varepsilon(\pi_1\chi\times \check{\pi}_2,s,\psi_F)}
{\varepsilon(\pi_1\times \check{\pi}_2,s,\psi_F)}  =
\chi(\mathrm{N}_{E/F}(\nu^{-r^2})) 
\frac{\varepsilon(\tau_1\chi_E\times \check{\tau}_2,s,\psi_E)}
{\varepsilon(\tau_1\times \check{\tau}_2,s, \psi_E)},
$$
where~$\nu=\nu(\theta_F,\psi_F,\psi_E)$.
\end{cor}

In the case~$E$ is maximal, totally ramified over~$F$
and~$\pi_1=\pi_2$ we recover~\cite{can}\S6.1~Corollaire~2, with (in
the notation of~\cite{can})~$c(\pi_1,\check{\pi}_1,\psi_F)=\mathrm{N}_{E/F}(\nu)$. 
Moreover, Corollary~\ref{change} implies~\cite{can}\S7.1~Th\'eor\`eme,
which describes how the constant~$c(\pi_1,\check{\pi}_1,\psi_F)$ changes under the tame lifting
operation.

%%%%%%%%%%%%%%%%%%%%%%%%%%%%%%%%%%%%%%%%%%%%%%%%%%%%%%%%
%%%%%%%%%%%%%%%%%%%%%%%% References %%%%%%%%%%%%%%%%%%%%
%%%%%%%%%%%%%%%%%%%%%%%%%%%%%%%%%%%%%%%%%%%%%%%%%%%%%%%%

\bigskip

%%%%%%%%%%%%%%%%%%%%%%%%%%%%%%%%%%%%%%%%%%%%%%%%%%%%%%%%
%%%%%%%%%%%%%%%%%%%%%%%%% Addresses %%%%%%%%%%%%%%%%%%%%
%%%%%%%%%%%%%%%%%%%%%%%%%%%%%%%%%%%%%%%%%%%%%%%%%%%%%%%%

\small

Vytautas Paskunas\hfill Shaun Stevens \break \indent
Fakult\"at f\"ur Mathematik%
\hfill School of Mathematics \break\indent
Universit\"at Bielefeld%
\hfill University of East Anglia \break\indent
Postfach 100131\hfill Norwich NR4 7TJ \break\indent
D-33501 Bielefeld\hfill United Kingdom \break\indent
Germany

\smallskip
paskunas@math.uni-bielefeld.de\hfill ginnyshaun@bigfoot.com \break\indent

\end{document}